\documentclass{amsart}
\usepackage{amsfonts}
\usepackage{amsmath}
\usepackage{graphics}
\usepackage{subfigure}
\usepackage[all]{xy}
\usepackage[colorlinks=true]{hyperref}
\usepackage{xspace}

\newtheorem{theorem}{Theorem}[section]
\newtheorem{corollary}[theorem]{Corollary} 
\newtheorem{lemma}[theorem]{Lemma}
\newtheorem{proposition}[theorem]{Proposition}

\theoremstyle{definition} 
\newtheorem{definition}[theorem]{Definition}
\newtheorem{remark}[theorem]{Remark} 
\newtheorem{example}[theorem]{Example} 

\newcommand{\ti}[1]{\tilde{{#1}}}
\newcommand{\baseRing}[1]{\ensuremath{\mathbb{#1}}}
\newcommand{\R}{\baseRing{R}} 
\newcommand{\BF}{\baseRing{F}} 
\newcommand{\BB}{\baseRing{B}} 
\newcommand{\N}{\baseRing{N}} 
\newcommand{\NZ}{\ensuremath{\baseRing{N}\cup\{0\}}} 
\newcommand{\Z}{\baseRing{Z}} 
\newcommand{\jdef}[1]{\emph{#1}}
\newcommand{\jgsg}{\ensuremath{\mathfrak{g}}\xspace}
\newcommand{\jgm}{\ensuremath{\mathfrak{m}}\xspace}
\newcommand{\stext}[1]{\ensuremath{\quad \text{{#1}} \quad}}
\newcommand{\del}{\partial}
\newcommand{\CC}{\ensuremath{{\mathcal{A}}}}
\newcommand{\CD}{\ensuremath{{\mathcal{A}_d}}}
\newcommand{\CDp}[1]{\ensuremath{{\mathcal{A}_d^{#1}}}}
\newcommand{\BM}{\ensuremath{{\mathcal{B}_m}}}
\newcommand{\BMp}{\ensuremath{{\mathcal{B}_m^+}}}
\newcommand{\BMm}{\ensuremath{{\mathcal{B}_m^-}}}
\newcommand{\BMR}{\ensuremath{{\hat{\mathcal{B}}_m}}}
\newcommand{\BMRp}{\ensuremath{{\hat{\mathcal{B}}_m^+}}}
\newcommand{\BMRm}{\ensuremath{{\hat{\mathcal{B}}_m^-}}}
\newcommand{\HLds}[1]{\ensuremath{\overline{h_d^{{#1}}}}}
\newcommand{\HLd}[1]{\ensuremath{{h_d^{{#1}}}}}

\newcommand{\HLc}[1]{\ensuremath{{h^{{#1}}}}}
\newcommand{\SG}{\ensuremath{G}}
\newcommand{\RS}{\ensuremath{\ti{\SG}\times (Q/\SG)}}
\newcommand{\RSsec}{\ensuremath{Q^{(2)}_\SG}}
\newcommand{\RSsecC}{\ensuremath{\check{Q}^{(2)}_{\SG}}}
\newcommand{\trint}{\;\;\makebox[0pt]{$\top$}\makebox[0pt]{$\cap$}\;\;}
\newcommand{\jiff}{\ensuremath{\Leftrightarrow}\xspace}
\newcommand{\imp}{\ensuremath{\Rightarrow}\xspace}
\newcommand{\annihilator}{\ensuremath{\circ}}

\newenvironment{roman-enumerate}{ \begin{enumerate}}{\end{enumerate}}

\DeclareMathOperator{\im}{Im}


\email{jfernand@ib.edu.ar}
\email{cora@mate.unlp.edu.ar}
\email{marce@mate.unlp.edu.ar}

\thanks{This research was partially supported by grants from the
  Universidad Nacional de Cuyo and Universidad Nacional de La
  Plata. C.T. was partially supported by a fellowship from CONICET.}

\begin{document}

\bibliographystyle{amsplain}

\title[Reduction of discrete mechanical systems]{Lagrangian reduction of
  nonholonomic discrete mechanical systems}
\author[Javier Fernandez, Cora Tori and Marcela Zuccalli]{}

\maketitle


\centerline{\scshape Javier Fernandez$^1$, Cora Tori$^{2}$ and
  Marcela Zuccalli$^2$}

\medskip {\footnotesize 

  \centerline{$^1$Instituto Balseiro, Universidad Nacional de Cuyo --
    C.N.E.A.}  

  \centerline{ Av. Bustillo 9500, San Carlos de
    Bariloche, R8402AGP, Rep\'ublica Argentina} }
\medskip
{\footnotesize
  \centerline{$^2$Departamento de Matem\'atica, Facultad de Ciencias Exactas,
    Universidad Nacional de La Plata} 
  \centerline{50 y 115, La Plata, Buenos Aires, 1900, Rep\'ublica
    Argentina} 
}


\begin{abstract}
  In this paper we propose a process of lagrangian reduction and
  reconstruction for nonholonomic discrete mechanical systems where
  the action of a continuous symmetry group makes the configuration
  space a principal bundle. The result of the reduction process is a
  discrete dynamical system that we call the discrete reduced
  system. We illustrate the techniques by analyzing two types of
  discrete symmetric systems where it is possible to go further and
  obtain (forced) discrete mechanical systems that determine the
  dynamics of the discrete reduced system.
\end{abstract}


\section{Introduction}
\label{sec:introduction}

The elimination of degrees of freedom of a symmetric mechanical
system, the basic goal of reduction theory, is an old subject that
dates back to the mid-nineteen century. The work of Routh in the
context of abelian symmetries of classical mechanical systems was
extended by many others in an effort to explore different aspects of
the reduction process. Currently, there are well developed theories of
reduction in the Hamiltonian setting, where the emphasis is in the
reduction of the symplectic and Poisson structures, and the lagrangian
setting, mostly focused on the reduction of variational
principles. The literature in this area is vast; modern references
are, for
instance,~\cite{ar:koiller-reduction_of_some_classical_nonholonomic_systems_with_symmetry,ar:bates_sniatycki-nonholonomic_reduction,bo:marsden_misiolek_ortega_perlmutter_ratiu-hamiltonian_reduction_by_stages,ar:cendra_marsden_ratiu-geometric_mechanics_lagrangian_reduction_and_nonholonomic_systems}.

Discrete time mechanical systems have been considered in the
literature since the 1960s, mostly as a way to approximate and model
the behavior of (continuous) mechanical systems
(see~\cite{ar:marsden_west-discrete_mechanics_and_variational_integrators}
and the references therein). The dynamics of nonholonomic discrete
mechanical systems has been introduced by J. Cort\'es and S. Martinez
in~\cite{ar:cortes_martinez-non_holonomic_integrators}. Many
characteristics of mechanical systems have a discrete analogue. One
important feature whose discrete analogue has been exposed only
partially in the literature is the reduction of symmetries. The
purpose of the present work is to describe a reduction and
reconstruction process for discrete time mechanical systems with
nonholonomic constraints, in the lagrangian setting.

There are many reasons for being interested in the reduction of a
symmetric system, be it continuous or discrete. On the one hand, the
reduced system has less degrees of freedom, which may be an advantage
when trying to solve the equations of motion. On the other, in many
cases, the reduced system encodes the core dynamics of the system of
interest. A good example of the advantage of working on the reduced
system in the discrete case is presented by S. Jalnapurkar et
al. in~\cite{ar:jalnapurkar_leok_marsden_west-discrete_routh_reducion}
where they show how the reduced system is free from geometric phases'
effects allowing them to observe interesting dynamical structures,
which were hard to separate in the unreduced system.

There are already several results on the reduction and reconstruction
of symmetric nonholonomic discrete mechanical systems available in the
literature. The case where the configuration space of the system $Q$
is a Lie group has been studied by Y. Fedorov and D. Zenkov
in~\cite{ar:fedorov_zenkov-dynamics_of_the_discrete_chaplygin_sleigh,ar:fedorov_zenkov-discrete_nonholonomic_ll_systems_on_lie_groups}
and by R. McLachlan and M. Perlmutter
in~\cite{ar:mclachlan_perlmutter-integrators_for_nonholonomic_mechanical_systems}. They
obtain reduced equations of motion on the Lie algebra of the Lie group
of symmetries which coincide, in the unconstrained case, with those
obtained by A. Bobenko and Y. Suris
in~\cite{ar:bobenko_suris-discrete_lagrangian_reduction_discrete_euler_poincare_equations_and_semidirect_products}. The
case of Chaplygin systems, that is, when the symmetry group acts in
such a way that no symmetry direction is compatible with the given
constraints is analyzed by J. Cort\'es
in~\cite{bo:cortes-non_holonomic} and, in context of groupoids, by
D. Iglesias et
al. in~\cite{ar:iglesias_marrero_martin_martines-discrete_nonholonomic_lagrangian_systems_on_lie_groupoids}. Also,
the unconstrained case has been treated by S. Jalnapurkar et
al. in~\cite{ar:jalnapurkar_leok_marsden_west-discrete_routh_reducion}
and, using discrete connections, by M. Leok et
al. in~\cite{ar:leok_marsden_weinstein-a_discrete_theory_of_connections_on_princip%
  al_bundles}. A different approach using groupoids is given
in~\cite{ar:iglesias_marrero_martin_martines-discrete_nonholonomic_lagrangian_systems_on_lie_groupoids};
in this case, the abstract reduction theory of groupoids produces a
general reduced system on a groupoid; our approach differs from
theirs, ours being more elementary ---we stay in the ``pair groupoid''
case--- and, at the same time, more explicit because we incorporate
additional information (connections) to our models.

The reduction results that we present here are modeled on the results
of H. Cendra et
al. in~\cite{ar:cendra_marsden_ratiu-geometric_mechanics_lagrangian_reduction_and_nonholonomic_systems}
for continuous systems. They consider the action of a Lie group $\SG$
on a configuration manifold $Q$ so that the quotient map
$\pi:Q\rightarrow Q/G$ is a principal bundle and the other data is
$\SG$-invariant. They relate the variational principle that determines
the dynamics of the original system with a reduced variational
principle, that determines the dynamics of the reduced
system. Eventually, equations of motion are derived from both
variational principles. The introduction of a connection on the
principal bundle $\pi$ serves them to, first, construct an isomorphic
model for the natural reduced space $TQ/\SG$ and, also, to split the
reduced variational principle and equations of motion into horizontal
and vertical parts. Appropriate choice of a connection can lead to a
simplified analysis of a specific mechanical system.

At a philosophical level, a common approach to discrete mechanics
consists of replacing the tangent bundle $TQ$ by $Q\times Q$, with the
idea that infinitesimal displacement ---velocities--- are replaced by
finite displacement ---pairs of points. Even though this is a powerful
idea, there are some very important differences between $TQ$ and
$Q\times Q$, that make it difficult to transfer techniques developed
for continuous systems to discrete systems, as we will see later. The
construction of the reduced space
in~\cite{ar:cendra_marsden_ratiu-geometric_mechanics_lagrangian_reduction_and_nonholonomic_systems}
relies on an isomorphism defined using a connection on the principal
bundle $Q\rightarrow Q/G$, which is seen as a $\SG$-invariant
splitting of $TQ$, that descends to a splitting of $TQ/\SG$. In the
discrete case, we follow the same path, but using what we call an
affine discrete connection, a minor generalization of the discrete
connections introduced by M. Leok et
al. in~\cite{ar:leok_marsden_weinstein-a_discrete_theory_of_connections_on_principal_bundles}
which allows us to split $(Q\times Q)/\SG$. Using, in addition, a
connection, we are able to derive a reduced variational principle as
well as reduced equations of motion that split in horizontal and
vertical parts.

The vertical part of the variational principle and equations turn out
to be equivalent to what is known as the discrete nonholonomic
momentum evolution
equation~\cite{ar:cortes_martinez-non_holonomic_integrators} that
determines how the discrete nonholonomic momentum mapping, defined in
terms of the symmetry directions that are compatible with the
constraints, evolves for discrete mechanical system. The so called
horizontal symmetries have the property of preserving the discrete
nonholonomic momentum.

We specialize the general reduction process in two different settings:
the Chaplygin and horizontal systems, where we can go further to
obtain reduced systems that can be described easily. The first case is
already present in the literature and we derive the known results from
our approach. The other hasn't been considered before in the discrete
setting to our knowledge; horizontal systems have been considered in
the continuous setting by Cort\'es in~\cite{bo:cortes-non_holonomic}.
Our approach can also be used in the case where the configuration
space is a Lie Group to re derive the results
of~\cite{ar:mclachlan_perlmutter-integrators_for_nonholonomic_mechanical_systems}
mentioned above.

Constraints in a mechanical system on a configuration space $Q$ are
usually given by two distributions: $\mathcal{D}$ that describes the
allowed variations and $\mathcal{C}_K$ that limits the allowed
trajectories. The D'Alembert Principle used to determine the dynamics
of such systems requires that $\mathcal{D}=\mathcal{C}_K$. More
general situations, called generalized nonholonomic systems, where no
connection is made between $\mathcal{D}$ and $\mathcal{C}_K$, have
been considered
in~\cite{ar:marle-kinematic_and_geometric_constraints_servomechanism_and_control_of_mechanical_systems,
  ar:cendra_ibort_de_leon_martin-a_generalization_of_chetaevs_principle_for_a_class_of_higher_order_nonholonomic_constraints,
  ar:cendra_grillo-generalized_nonholonomic_mechanics_servomechanisms_and_related_brackets}
. Discrete mechanical systems as introduced
in~\cite{ar:cortes_martinez-non_holonomic_integrators} should be
called generalized in the same sense as above, since the setup for
constrained discrete mechanical systems consists of $\mathcal{D}$,
with the same meaning as in the continuous case, and a submanifold
$\mathcal{D}_d\subset Q\times Q$ which restricts the allowable
discrete trajectories, and no relation between the two is assumed. 

\vskip 0.3cm

The layout of the paper is as follows. In
Section~\ref{sec:reduction_of_continuous_systems}, we sketch very
roughly the lagrangian reduction of classical mechanical systems,
which serves the dual purpose of introducing the results in the
continuous case as well as the variational approach that serves as
motivation for our work in the discrete case.
Section~\ref{sec:discrete_mechanical_systems_and_symmetries}
introduces discrete mechanical systems and their symmetries. In
Section~\ref{sec:discrete_tools} we introduce some tools, including
the affine discrete connections, that will be useful in the analysis
of the reduction process to be carried out in
Section~\ref{sec:variations_and_reduced_variations} and condensed as
Theorem~\ref{th:4_points-general}. The statement of this result is
written in terms that are not obviously defined on the reduced system;
the purpose of
Section~\ref{sec:intrinsic_version_of_the_reduced_equations_of_motion}
is to give an intrinsic version of that result, which we achieve with
Corollary~\ref{cor:4_points_general-intrinsic}. The reconstruction of
the original dynamics starting from that of the reduced system is
explored in Section~\ref {sec:reconstruction}. Finally, in
Section~\ref{sec:nonholonomic_discrete_momentum} we study the discrete
nonholonomic momentum mapping, showing that the vertical part of the
reduced variational principle is equivalent to the discrete
nonholonomic momentum equation.
The rest of the paper deals with the specialization of the general
theory developed so far to particular situations. In
Section~\ref{sec:reduced_equations_of_motion_for_trivial_bundles} we
obtain the equations of motion for systems where $Q\rightarrow Q/\SG$
is a trivial principal bundle. In
Sections~\ref{sec:reduced_equations_of_motion-chaplygin}
and~\ref{sec:reduced_equations_of_motion-horizontal} we study systems
with Chaplygin and horizontal symmetries respectively; in the
Chaplygin case we obtain intrinsic versions that specialize to the
results of~\cite{bo:cortes-non_holonomic}, while in the horizontal
case we find the discrete version of a type of symmetry whose
continuous counterpart had been studied
in~\cite{bo:cortes-non_holonomic}.

\vskip 0.3cm
Last, we wish to thank Hern\'an Cendra for his interest and valuable
comments on this work.


\section{Reduction of classical mechanical systems}
\label{sec:reduction_of_continuous_systems}

In this section we recall some basic facts of the lagrangian reduction
theory of (generalized) nonholonomic mechanical systems in the
presence of symmetry. We refer
to~\cite{ar:bloch-nonholonomic_symmetry} and
~\cite{ar:cendra_ferraro_grillo-lagrangian_reduction_of_generalized_NHS}
for further details.


\subsection{Generalized nonholonomic mechanical systems}
\label{sec:generalized_non_holonomic_mechanical_systems}

\begin{definition}\label{def:continuous_mechanical_system}
  A \jdef{generalized nonholonomic mechanical system} is a quadruple
  $(Q,L,\mathcal{D},C_K)$ where $Q$ is a differentiable manifold, the
  \jdef{configuration space}, $L:TQ\rightarrow\R$ is a smooth function on
  the tangent bundle of $Q$, the \jdef{lagrangian}, $\mathcal{D}$ is a
  subbundle of $TQ$, the \jdef{variational constraints} or
  \jdef{virtual displacements}, and $C_K\subset TQ$ is a submanifold, the
  \jdef{kinematic constraints}.
\end{definition}

For every such system, the \jdef{action functional} is defined by
\begin{equation*}
  S(q):=\int_{t_0}^{t_1} L(q(t),\dot{q}(t)) dt,
\end{equation*}
where $q:[t_0,t_1]\rightarrow Q$ is a smooth curve in $Q$ and
$\dot{q}:[t_0,t_1]\rightarrow TQ$ is its velocity. An
\jdef{infinitesimal variation} of $q$ is a smooth curve $\delta
q:[t_0,t_1]\rightarrow TQ$ . An infinitesimal
variation is said to have vanishing end points if $\delta q(t_0)=0$ and
$\delta q(t_1)=0$.

The dynamics of a generalized nonholonomic mechanical system is
determined by the following Principle.

\begin{definition}[Lagrange--D'Alembert Principle]
  A trajectory of $(Q,L,\mathcal{D},C_K)$ is a curve
  $q:[t_0,t_1]\rightarrow Q$ which
  \begin{itemize}
  \item satisfies the kinematic constraints: $(q(t),\dot{q}(t))\in
    C_K$ for all $t\in [t_0,t_1]$ and
  \item is a critical point of $S$ for the admissible
    variations: $dS(q)(\delta q)=0$ for all infinitesimal
    variations $\delta q$ of $q$ with vanishing end points and
    such that $\delta q(t) \in \mathcal{D}_{q(t)}$ for all $t\in
    [t_0,t_1]$.
  \end{itemize}
\end{definition}

As usual, the Lagrange--D'Alembert Principle gives rise to a set of
equations called the generalized Lagrange--D'Alembert
equations for $(Q,L,\mathcal{D},C_K)$. These equations of motion can
be written in a coordinate-free way if an affine connection $\nabla$
on $Q$ is chosen. Assuming this additional
datum,~\cite{ar:cendra_ferraro_grillo-lagrangian_reduction_of_generalized_NHS}
proves the following result.

\begin{theorem}
  A smooth curve $q$ in $Q$ is a trajectory of $(Q,L,\mathcal{D},C_K)$
  if and only if $(q(t),\dot{q}(t))\in C_K$ for all $t\in [t_0,t_1]$ and
  \begin{equation}
    \label{eq:generalized_lagrange_dalambert}
    -\frac{D}{Dt}\BF L(q(t),\dot{q}(t)) + \BB L (q(t),\dot{q}(t)) \in
    (\mathcal{D}_{q(t)})^\annihilator \stext{ for all } t\in [t_0,t_1],
  \end{equation}
  where $\BF L$ and $\BB L$ denote the fiber and base derivatives of
  $L$, as defined in~\cite[Definition
  3]{ar:cendra_ferraro_grillo-lagrangian_reduction_of_generalized_NHS}
  and $(\mathcal{D}_{q(t)})^\annihilator\subset T^*_{q(t)}Q$ is the annihilator of
  $\mathcal{D}_{q(t)}\subset T_{q(t)}Q$.
\end{theorem}

Condition~\eqref{eq:generalized_lagrange_dalambert} is known as the
(generalized) Lagange--D'Alembert equation.


\subsection{Symmetric generalized nonholonomic mechanical systems}
\label{sec:generalized_non_holonomic_mechanical_systems_with_symmetry}

In what follows, we assume that $\SG$ is a Lie group that acts (on the
left) on $Q$ freely and properly, so that the quotient
$\pi:Q\rightarrow Q/\SG$ is a principal bundle with structure
group $\SG$. This action will be denoted by $l_g^Q(q)$
for all $g\in \SG$ and $q\in Q$; there is an induced action of $\SG$
on $TQ$ given by $l^{TQ}_g(v_q):=dl^Q_g(q)(v_q)$ for all $g\in \SG$
and $v_q\in T_qQ$, that is called the \jdef{lifted action}.

\begin{definition}
  $\SG$ is a \jdef{symmetry group} of $(Q,L,\mathcal{D},C_K)$ if, in
  addition to the general assumptions stated above, $L$, $\mathcal{D}$
  and $C_K$ are $\SG$-invariant, that is, if $L\circ l^{TQ}_g = L$,
  $l^{TQ}_g(\mathcal{D})\subset\mathcal{D}$ and $l^{TQ}_g(C_K)\subset
  C_K$ for all $g\in \SG$.
\end{definition}

The \jdef{vertical bundle} $\mathcal{V}^\SG$ over $Q$ is a subbundle
of $TQ$ with fibers $\mathcal{V}^\SG_q:=T_q(l^Q_\SG(\{q\}))$. We
assume that the distribution $\mathcal{S}$ with fibers
$\mathcal{S}_q:=\mathcal{V}^\SG_q\cap \mathcal{D}_q$ has locally
constant rank, so that it is a subbundle of $TQ$. This condition,
which applies to many interesting examples, is verified when, for
instance, $T_qQ = \mathcal{D}_q + \mathcal{V}_q$ for all $q\in Q$
---the ``dimension assumption''--- holds, as is usually considered
in much of the literature.

Assume that $\SG$-invariant subbundles of $TQ$, $\mathcal{H}$,
$\mathcal{U}$ and $\mathcal{W}$ can be chosen in such a way that
$\mathcal{H}$ and $\mathcal{U}$ are direct complements of
$\mathcal{S}$ in $\mathcal{D}$ and $\mathcal{V}^\SG$ respectively, and
$\mathcal{W}$ is a direct complement of $\mathcal{D}+\mathcal{V}^\SG$
in $TQ$. All together, we have the decomposition
\begin{equation}
  \label{eq:TQ_decomposition} 
  TQ = \mathcal{W} \oplus  \mathcal{U} \oplus \mathcal{S} \oplus \mathcal{H},
\end{equation}
where $\mathcal{D} = \mathcal{S} \oplus \mathcal{H}$ and
$\mathcal{V}^\SG = \mathcal{S} \oplus \mathcal{U}$.  One way to
construct the complementary bundles $\mathcal{H}$, $\mathcal{U}$ and
$\mathcal{W}$ is as orthogonal complements for some $\SG$-invariant
inner product on $TQ$. This type of inner product is usually available
as the kinetic energy of symmetric mechanical systems.

Following~\cite{ar:bloch-nonholonomic_symmetry}
and~\cite{ar:cendra_ferraro_grillo-lagrangian_reduction_of_generalized_NHS},
we associate a connection to the
decomposition~\eqref{eq:TQ_decomposition}.
\begin{definition}\label{def:non_holonomic_connection}
  The unique connection $\CC$ on the principal bundle
  $\pi:Q\rightarrow Q/\SG$ whose horizontal space is
  $Hor_{\CC}=\mathcal{W}\oplus \mathcal{H}$ is called the
  \jdef{generalized nonholonomic connection} associated to the system
  $(Q,L,\mathcal{D},C_K)$ and the
  splitting~\eqref{eq:TQ_decomposition}.
\end{definition}


\subsection{Reduction of symmetry}
\label{sec:reduction_of_symmetry_of_a_generalized_non_holonomic_mechanical_system}

When $\SG$ is a symmetry group of $(Q,L,\mathcal{D},C_K)$
we define the \jdef{reduced lagrangian} $\ell:TQ/\SG\rightarrow \R$ by
$\ell([(q,\dot{q})]_\SG) := L(q,\dot{q})$.

In order to establish a reduced variational principle and study the
reduced equations of
motion~\cite{bo:cendra_marsden_ratiu-lagrangian_reduction_by_stages}
and~\cite{ar:cendra_marsden_ratiu-geometric_mechanics_lagrangian_reduction_and_nonholonomic_systems}
find it more convenient to work on a diffeomorphic model of $TQ/\SG$
as follows. Given a principal connection $\CC$ on $\pi:Q\rightarrow
Q/\SG$ there is an isomorphism of vector bundles over $Q/\SG$
\begin{equation*}
  \alpha_\CC:TQ/\SG \rightarrow T(Q/\SG)\oplus\ti{\jgsg},
\end{equation*}
where $\jgsg$ is the Lie algebra of $\SG$, $\ti{\jgsg}$ is the adjoint
vector bundle (see~\cite{bo:kobayashi_nomizu-foundations-v1}) of $\pi$
and $\alpha_\CC([q,\dot{q}]_\SG) := d\pi(q,\dot{q})\oplus
[(q,\CC(q,\dot{q}))]_\SG$. The connection of choice used to define
$\alpha_\CC$ is the generalized nonholonomic connection $\CC$
introduced above. The reduced lagrangian $\ell$ can be transported to
this new model space as
$\hat{L}:T(Q/\SG)\oplus\ti{\jgsg}\rightarrow\R$, defined by
\begin{equation*}
  \hat{L}(x,\dot{x},\bar{v}) := 
  \ell\big((\alpha_\CC)^{-1}(x,\dot{x},\bar{v})\big) = 
  L(q,\dot{q}),
\end{equation*}
where $(x,\dot{x}) = d\pi(q,\dot{q})$ and $\bar{v} =
[(q,\CC(q,\dot{q}))]_\SG$. 

The reduced dynamics of the system is defined
in~\cite{ar:cendra_ferraro_grillo-lagrangian_reduction_of_generalized_NHS}
by the following principle.

\begin{definition}[Generalized Lagrange--D'Alembert--Poincar\'e Principle]
  A trajectory of the reduced system determined by the symmetry group
  $\SG$ of $(Q,L,\mathcal{D},C_K)$ is a curve $\mu:= (x,\dot{x})\oplus
  \bar{v}: [t_0,t_1]\rightarrow T(Q/\SG)\oplus\ti{\jgsg}$ which
  \begin{itemize}
  \item satisfies the reduced kinematic constraints:
    $\mu(t)\in\hat{C}_K:=\alpha_\CC(C_K/\SG)$ for all $t\in
    [t_0,t_1]$.
  \item is a critical point of the reduced action
    \begin{equation*}
      \int_{t_0}^{t_1} \hat{L}(x(t),\dot{x}(t),\bar{v}(t)) dt
    \end{equation*}
    for some infinitesimal variations $\delta\mu=\delta x\oplus
    \delta^\CC\bar{v}$ of $\mu$ such that, if $q$ is a lift of $x$ to
    $Q$, $\delta x(t)\in
    \hat{\mathcal{D}}^h_{x(t)}:=d\pi(q(t))(\mathcal{D}_{q(t)})$ and
    $\delta^\CC\bar{v}(t)\in
    \hat{\mathcal{D}}^v:=\alpha_\CC(\mathcal{S}/\SG)$ with all
    variations vanishing at the end points. The precise description of
    the set of infinitesimal variations where the criticality
    condition applies requires a special form for the
    $\delta^\CC\bar{v}$ involving the curvature of $\CC$ (we refer
    to~\cite{ar:cendra_ferraro_grillo-lagrangian_reduction_of_generalized_NHS}
    for further details).
  \end{itemize}
\end{definition}

The dynamics generated by this principle can be described using
equations of motion, called the reduced generalized
Lagrange--D'Alembert--Poincar\'e equations, which are equivalent to
the generalized nonholonomic Lagrange--D'Alembert equations of the
system $(Q,L,\mathcal{D},C_K)$, as shown by Theorem 9
of~\cite{ar:cendra_ferraro_grillo-lagrangian_reduction_of_generalized_NHS}.

\begin{theorem}\label{th:generalized_reduction-continuous}
  Let $\SG$ be a symmetry group of $(Q,L,\mathcal{D},C_K)$, $\CC$ be
  the generalized nonholonomic connection on $\pi:Q\rightarrow Q/\SG$
  associated to a splitting of $TQ$ and $q:[t_0,t_1]\rightarrow Q$ a
  curve on $Q$. Then the following statements are equivalent.
  \begin{itemize}
  \item The curve $q$ satisfies $(q(t),\dot{q}(t))\in C_K$ for all
    $t\in [t_0,t_1]$ and~\eqref{eq:generalized_lagrange_dalambert}
    holds.
  \item The curve $\mu:[t_0,t_1]\rightarrow T(Q/\SG)\oplus \ti{\jgsg}$
    given by $\mu:=(x,\dot{x})\oplus \bar{v}
    =\alpha_\CC([(q,\dot{q})]_\SG)$ satisfies $\mu(t)\in \hat{C}_K$
     and
    \begin{gather*}
      -\frac{D}{Dt}\frac{\partial \hat{L}}{\partial \bar{v}}(\mu(t)) +
      ad^*_{\bar{v}}\frac{\partial \hat{L}}{\partial \bar{v}}(\mu(t)) 
      \in (\hat{\mathcal{D}}_{x(t)}^v)^\annihilator\\
      -\frac{D}{Dt} \frac{\partial \hat{L}}{\partial \dot{x}}(\mu(t))
      + \frac{\partial \hat{L}}{\partial x}(\mu(t)) - \big\langle
      \frac{\partial\hat{L}}{\partial\bar{v}},i_{x(t)}\ti{B}\big\rangle
      \in (\hat{\mathcal{D}}_{x(t)}^h)^\annihilator
    \end{gather*}
    for all $t\in [t_0,t_1]$, where
    $\ti{B}$ is the reduced curvature of $\CC$.
  \end{itemize}

\end{theorem}


\section{Discrete mechanical systems and symmetries}
\label{sec:discrete_mechanical_systems_and_symmetries}

In this section we review the notion of discrete mechanical system
with nonholonomic constraints, that is the discrete time analogue of
the (generalized) nonholonomic mechanical systems considered in
Section~\ref{sec:reduction_of_continuous_systems}. We also consider
symmetries of such systems.


\subsection{Discrete mechanical systems}
\label{sec:discrete_mechanical_systems}

\begin{definition}
  A \emph{nonholonomic discrete mechanical system} consists of a
  quadruple $(Q,L_d,\mathcal{D},\mathcal{D}_d)$ where $Q$ and
  $\mathcal{D}$ are as in
  Definition~\ref{def:continuous_mechanical_system}, $L_d:Q\times Q
  \rightarrow\R$ is a smooth map, the \emph{discrete lagrangian}, and
  $\mathcal{D}_d\subset Q\times Q$ is a submanifold, the
  \emph{discrete kinematic constraints}.
\end{definition}

\begin{remark}
  The discrete mechanical systems defined above are slightly more
  general than those considered
  in~\cite{ar:cortes_martinez-non_holonomic_integrators} because we
  are not requiring that the diagonal of $Q\times Q$ be contained in
  $\mathcal{D}_d$. Still, in order to construct a dynamical system, we
  will assume that $\mathcal{D}_d$ contains the graph of a smooth map
  $Q\rightarrow Q$, with the case
  of~\cite{ar:cortes_martinez-non_holonomic_integrators} corresponding
  to the identity.
\end{remark}

Following~\cite{ar:cortes_martinez-non_holonomic_integrators}, discrete
mechanical systems define discrete dynamical systems using a discrete
Lagrange--D'Alembert Principle, roughly saying that trajectories of
the dynamical system are critical points of the \emph{discrete action
  functional}
\begin{equation*}
  S_d(q_\cdot ) := \sum_{k=0}^{N-1}L_d(q_{k},q_{k+1})
\end{equation*}
that satisfy the constraints. The following definition makes this
notion precise.

\begin{definition}[Discrete Lagrange--D'Alembert Principle]
  \label{def:discrete_lagrange_dalembert_principle}
  A \emph{discrete curve} in $Q$ is a map $q_\cdot : \{0,1,\ldots, N\}
  \rightarrow Q$ and a \emph{variation} of a discrete curve $q_\cdot$
  consists of a map $\delta q_\cdot:\{0,1,\ldots N\}\rightarrow TQ$
  such that $\delta q_k\in T_{q_k}Q$ for all $k$. A variation is said
  to have \emph{vanishing end points} if $\delta q_0 = 0$ and $\delta
  q_N = 0$. A \emph{trajectory} of $(Q,L_d,\mathcal{D},\mathcal{D}_d)$
  is a discrete curve $q_\cdot$ which
  \begin{itemize}
  \item satisfies the kinematic constraints $(q_k,q_{k+1})\in
    \mathcal{D}_d$ for all $k$ and
  \item is a critical point of $S_d$ for all admissible variations
    $\delta q_\cdot$ of $q_\cdot$: $dS_d(q_\cdot)(\delta q_\cdot)=0$
    for all infinitesimal variations with vanishing end points $\delta
    q_k\in \mathcal{D}_{q_k}$ for all $k$.
  \end{itemize}
\end{definition}

The discrete Lagrange--D'Alembert Principle leads to a set of
equations. Indeed, if $D_j$ denotes differentiation with respect to
the $j$-th component of a Cartesian product,
\begin{equation*}
  \begin{split}
    dS_d(q_\cdot)(\delta q_\cdot) =& \sum\nolimits_{k=0}^{N-1}
    dL_d(q_k,q_{k+1})(\delta
    q_k ,\delta q_{k+1}) \\
    =& \sum\nolimits_{k=1}^{N-1} \big(D_1L_d(q_k,q_{k+1}) +
    D_2L_d(q_{k-1},q_k)\big)
    (\delta q_k) \\
    & + D_1L_d(q_0,q_1)(\delta q_0) + D_2L_d(q_{N-1},q_N)(\delta q_N),
  \end{split}
\end{equation*}
so that $q_\cdot$ is a trajectory of
$(Q,L_d,\mathcal{D},\mathcal{D}_d)$ if and only if, for all $k$,
\begin{equation} \label{eq:dla-eqs} D_1L_d(q_k,q_{k-1}) +
  D_2L_d(q_{k-1},q_k) \in \mathcal{D}_{q_k}^\annihilator \stext{ and }
  (q_k,q_{k+1})\in\mathcal{D}_d
\end{equation}
or, if $\chi_\kappa(q_k,q_{k+1})=0$ are local equations defining
$\mathcal{D}_d$,
\begin{equation*}
  D_1L_d(q_k,q_{k-1}) + D_2L_d(q_{k-1},q_k) = \sum\nolimits_{a=1}^M \lambda_{a,k}\,
    \omega^a(q_k) \stext{ and }
    \chi_\kappa(q_k,q_{k+1})=0,
\end{equation*}
for some constants $\lambda_{a,k}\in\R $ and
$\mathcal{D}_{q_k}^\annihilator = \langle \omega^1(q_k),\ldots,
\omega^M(q_k)\rangle$.

It can be shown that under sufficient regularity of the data,
equation~\eqref{eq:dla-eqs} has solutions, that are unique given
sufficiently closely spaced initial data (see~\cite[Proposition 3]
{ar:mclachlan_perlmutter-integrators_for_nonholonomic_mechanical_systems}
and~\cite{ar:cortes_martinez-non_holonomic_integrators}).

Sometimes it is necessary to consider systems that are forced. Such
will be the case below when we consider reduced systems, even when the
unreduced one is not forced. Thus, it is convenient to add the
following notion.

\begin{definition}
  A \emph{forced discrete mechanical system} consists of a discrete
  mechanical system $(Q,L_d,\mathcal{D},\mathcal{D}_d)$ together with
  a $1$-form $f_d$ on $Q\times Q$. We often write $f_d(q_0,q_1)(\delta
  q_0,\delta q_1) := f_d^-(q_0,q_1)(\delta q_0) +
  f_d^+(q_0,q_1)(\delta q_1)$ were $f_d^+: p_2^*(TQ)\rightarrow \R $
  and $f_d^-: p_1^*(TQ)\rightarrow \R $, where $p_j:Q\times
  Q\rightarrow Q$ denotes the projection on the $j$-th component and
  $p_j^*(TQ)$ is the pullback vector bundle
  (see~\cite{bo:AM-mechanics}) over $Q\times Q$ .
\end{definition}

The dynamics of a forced discrete mechanical system is given by the
appropriately modified discrete Lagrange--D'Alembert principle as
follows.

\begin{definition}
  A discrete curve $q_\cdot$ is a \jdef{trajectory of the forced
    discrete lagrangian mechanical system}
  $(Q,L_d,\mathcal{D},\mathcal{D}_d,f_d)$ if
  \begin{itemize}
  \item the curve satisfies the kinematic constraints:
    $(q_k,q_{k+1})\in \mathcal{D}_d$ for all $k$ and
  \item for all variations $\delta q_\cdot$ of $q_\cdot$ that have
    vanishing end points and $\delta q_k\in\mathcal{D}_{q_k}$,
    \begin{equation*}
      dS_d(q_\cdot)(\delta q_\cdot) 
      +\sum_{k=0}^{N-1} f_d(q_k,q_{k+1})(\delta q_k, \delta q_{k+1})
      = 0.
    \end{equation*}
  \end{itemize}
\end{definition}

Just as we found a characterization in terms of equations in the
unforced situation, $q_\cdot$ is a trajectory of a forced system if
and only if, for all $k$,
\begin{equation*} 
  \begin{cases}
    D_1L_d(q_k,q_{k+1}) + D_2L_d(q_{k-1},q_k) + f_d^-(q_k,q_{k+1}) +
    f_d^+(q_{k-1},q_k) \in \mathcal{D}_{q_k}^\annihilator\\
    (q_k,q_{k+1})\in\mathcal{D}_d.
  \end{cases}
\end{equation*}


\subsection{Symmetric discrete mechanical systems}
\label{sec:symmetric_dicsrete_mechanical_systems}

The main subject of this paper is the analysis of discrete mechanical
systems with continuous symmetry groups. Here we continue under the
assumption that the Lie group $\SG$ acts on $Q$ by $l^Q$ making the
quotient map $\pi:Q\rightarrow Q/\SG$ a principal bundle. In
addition, we consider the \jdef{diagonal action} of $\SG$ on $Q\times
Q$ defined by $l^{Q\times Q}_g(q_0,q_1):=(l^Q_g(q_0),l^Q_g(q_1))$.

\begin{definition}
  $\SG$ is a \emph{symmetry group} of
  $(Q,L_d,\mathcal{D},\mathcal{D}_d)$ if, in addition to what was
  stated in the previous paragraph, $L_d$ and $\mathcal{D}_d$ are
  invariant by $l^{Q\times Q}$, and $\mathcal{D}$ is $\SG$-invariant
  by $l^{TQ}$. If the system is forced, we also require that the
  discrete force $f_d$ be $\SG$-equivariant.
\end{definition}

Several structures introduced in
Section~\ref{sec:generalized_non_holonomic_mechanical_systems_with_symmetry}
remain useful in the current context. In particular, we have the
subbundles
$\mathcal{S},\mathcal{V}^\SG,\mathcal{H},\mathcal{U},\mathcal{W}
\subset TQ$ and the decomposition of $TQ$ given
by~\eqref{eq:TQ_decomposition}. We also have a nonholonomic
connection $\CC$ associated to that decomposition. We denote the
horizontal lift associated with ${\CC}$ by $\HLc{}: Q\times_{Q/\SG}
T(Q/\SG) \rightarrow TQ$, that is $\HLc{q}(w_{\pi(q)}) = v_q$ if $v_q
\in Hor_{\CC}(q)$ and $d\pi(q)(v_q)=w_{\pi(q)}$. The following result
is straightforward.

\begin{lemma}
  Considering the $\SG$-actions given by the lifted action $l^{TQ}$
  and by $l_g^{Q\times_{Q/\SG} T(Q/\SG)}(q,t_{\pi(q)}) :=
  (l_g^Q(q),t_{\pi(q)})$, the map $\HLc{}$ is $\SG$-equivariant. 
\end{lemma}

For each $q\in Q$, the elements $\xi\in\jgsg$ such that $\xi_Q(q)\in
\mathcal{D}_q$ form a subspace in $\jgsg$, which may depend on
$q$. For that reason, it is convenient to define the space
\begin{equation*}
  \jgsg^\mathcal{D} := \{ (q,\xi) \in Q\times \jgsg : 
  \xi_Q(q)\in \mathcal{D}_q\}.
\end{equation*}
Using the projection on the first variable, $\jgsg^\mathcal{D}$ is a
vector bundle on $Q$. It is easy to see that the differential in the
group direction of $l^Q$ establishes an isomorphism of vector bundles
$\jgsg^\mathcal{D} \simeq \mathcal{S}$.


\section{Some discrete tools}
\label{sec:discrete_tools}

An important step in the reduction process for classical systems is
the passage from $TQ/\SG$ to a model space $T(Q/\SG)\oplus\ti{\jgsg}$,
as seen in Section~\ref{sec:reduction_of_continuous_systems}. This is
achieved using the nonholonomic connection. In the discrete case, we
follow the same philosophy but connections are not the right tool for
the task. In this Section we introduce the notion of affine discrete
connection and use it to construct an isomorphism that, eventually,
will play the role that $\alpha_\CC$ played in
Section~\ref{sec:reduction_of_symmetry_of_a_generalized_non_holonomic_mechanical_system}.


\subsection{Affine discrete connections}
\label{sec:affine_discrete_connections}

\begin{definition}
  The \emph{discrete vertical bundle} for the $l^Q$ action of $\SG$ is
  the submanifold
  \begin{equation*}
    \mathcal{V}^\SG_d:= \{(q,l^Q_g(q)) \in Q\times Q: q\in Q, g\in \SG\}.
  \end{equation*}
  For $q\in Q$ we define $\mathcal{V}^\SG_d(q) := \mathcal{V}^\SG_d \cap
  (\{q\}\times Q)\subset Q\times Q$.
\end{definition}

\begin{definition}
  \label{def:composition} 
  The \emph{composition} of vertical and arbitrary elements of
  $Q\times Q$ with the same first element is defined by $\cdot :
  \mathcal{V}^\SG_d\times_Q (Q\times Q) \rightarrow Q\times Q$ with
  \begin{equation*}
    (q_0,l^Q_g(q_0))\cdot (q_0,q_1) := (q_0, l^Q_g(q_1)).
  \end{equation*}
\end{definition}

We denote by $l^\SG$ the conjugation action of $\SG$ on itself, that
is, $l^\SG_g(h):=g h g^{-1}$ for all $g,h\in\SG$.

\begin{definition}
  \label{def:affine_discrete_connection} 
  Let $\gamma:Q\rightarrow \SG$ be a smooth $\SG$-equivariant map with
  respect to $l^Q$ and $l^\SG$, $\Gamma:=\{(q,l^Q_{\gamma(q)}(q)):
  q\in Q\}$ and $Hor\subset Q\times Q$ be a $\SG$-invariant
  submanifold such that $%
  \Gamma\subset Hor$. For each $q\in Q$ let $Hor(q):=Hor \cap
  (\{q\}\times Q)$ and $Hor^2(q) := p_2(Hor(q))$.

  We say that $Hor$ defines the \emph{affine discrete connection}
  $\CD$ on the principal bundle $\pi:Q\rightarrow Q/\SG$ if, for
  each $q\in Q$ and all $q_1\in l^Q_\SG(\{q\})\cap Hor^2(q)$, 
  \begin{equation}\label{eq:affine_discrete_connection_transversality-0}
    Hor^2(q)\subset Q \text{ is a submanifold} \stext{and} 
    T_{q_1}Q = T_{q_1} (l^Q_\SG(\{q\})) \oplus T_{q_1} Hor^2(q)
  \end{equation}
  (see Figure~\ref{fig:orbit_and_horizontal-a}).  We denote $Hor$ by
  $Hor_{\CD}$, and call $\gamma$ (or even $\Gamma$) the \emph{level}
  of $\CD$.
\end{definition}

\begin{remark}
  Recalling the notion of transversality of submanifolds
  (see~\cite{bo:Guillemin-Pollack-differential_topology}),
  condition~\eqref{eq:affine_discrete_connection_transversality-0} is
  equivalent to requiring that for each $q\in Q$,
  \begin{equation*}
    \dim(Hor^2(q))=\dim Q-\dim \SG \stext{and} l^Q_\SG(\{q\})\trint Hor^2(q).
  \end{equation*}
\end{remark}

\begin{remark}
  The notion of discrete connection introduced
  in~\cite{ar:leok_marsden_weinstein-a_discrete_theory_of_connections_on_principal_bundles},
  coincides with that of an affine discrete connection where the level
  $\Gamma$ is the diagonal of $Q\times Q$. Since the diagonal of
  $Q\times Q$ plays the role of a ``null element'' for composition,
  affine discrete connections need not contain the ``null element'' in
  their horizontal space, just like affine spaces need not contain the
  null element of a vector space.
\end{remark}

\begin{remark}
  The idea behind the introduction of the previous definitions is that
  $Q\times Q$ should be a discrete version of $TQ$. Even though this
  is a powerful idea, there are several important differences between
  those spaces. Using the projection on the first factor and the
  standard projection, both spaces are fibered over $Q$, but $TQ$ is a
  vector bundle, while $Q\times Q$ is usually not one. In particular,
  tangent vectors at the same base point can be added, whereas there
  is nothing similar for elements of $Q\times Q$. A partial fix for
  this problem is the composition operation introduced in
  Definition~\ref{def:composition}: it provides a way of combining
  vertical elements of $Q\times Q$ with arbitrary elements of $Q\times
  Q$ based at the same point (that is, with the same first
  component). Even though this is not a complete analogue of addition,
  it is enough to handle discrete connections.

  Connections are an important tool in differential geometry.
  Essentially, for principal bundles $\pi:Q\rightarrow Q/\SG$,
  they provide compatible splittings
  \begin{equation}
    \label{eq:decomposition_TQ-CC} 
    T_qQ = \mathcal{V}^\SG_q  \oplus Hor(q) \stext{for all} q\in Q.
  \end{equation}
  In the discrete setting, we want to be able to split $\{q\}\times Q$
  in vertical and (some) complementary space. Not having an addition
  operation we choose a more geometric view of the problem and define
  affine discrete connections in terms of spaces that are
  complementary to $\mathcal{V}^\SG_d$ in the precise sense of
  Definition~\ref{def:affine_discrete_connection}.  Furthermore, this
  definition is equivalent to being able to decompose $\{q\}\times Q =
  \mathcal{V}^\SG_d(q)\cdot Hor(q)$ (at least in a neighborhood of
  $\Gamma$), which is the composition-analogue
  of~\eqref{eq:decomposition_TQ-CC}.
\end{remark}

\begin{proposition}\label{prop:CD_imp_pair_decomposition}
  Let $\CD$ be an affine discrete connection of level $\gamma$.  Then,
  there exists $U$, a $\SG$-invariant open neighborhood of $\Gamma$ in
  $Q\times Q$, such that, for all $(q_0,q_1)\in Q\times Q$ with
  $\pi(q_1)$ sufficiently close to $\pi(q_0)$, there is a unique $g\in
  \SG$ such that
  \begin{equation}  \label{eq:pair_decomposition_Ad}
    (q_0,q_1) = (q_0,l^Q_g(q_0))\cdot (q_0,l^Q_{g^{-1}}(q_1)),
  \end{equation}
  with $(q_0,l^Q_{g^{-1}}(q_1))\in Hor_{\CD}(q_0)\cap U$.
\end{proposition}

\begin{proof}
  Since $l^Q_{\gamma(q_0)}(q_0)\in l_\SG^Q(\{q_0\})\cap
  Hor^2_{\CD}(q_0)$, the condition $l^Q_\SG(\{q_0\})\trint
  Hor_{\CD}^2(q_0)$ implies that $Hor^2_{\CD}(q_0)$ intersects all
  orbits of $\SG$ that are close to the one through $q_0$ (see
  Figure~\ref{fig:orbit_and_horizontal-b}). Therefore, for any $q_1$
  with $\pi(q_1)$ sufficiently close to $\pi(q_0)$, there are $g\in
  \SG$ such that $l^Q_{g^{-1}}(q_1)\in Hor^2_{\CD}(q_0)$. Furthermore, again
  by the transversality condition, there is an open set
  $V_{q_0}\subset Q$ where the intersection $l^Q_\SG(\{q_1\})\cap
  Hor^2_{\CD}(q_0) \cap V_{q_0}$ consists of a single point. Hence,
  there is a unique $g\in \SG$ such that $l^Q_{g^{-1}}(q_1)\in
  Hor^2_{\CD}(q_0) \cap V_{q_0}$. Finally, taking into account the
  smoothness of $Hor$, $U$ is constructed by gluing the sets
  $V_{q_0}$.
  \begin{figure}[htbp]
    \centering \subfigure[$\SG$ orbit ($l^Q_\SG(q_0)$) and
    $Hor_{\CD}^2(q_0)$\label{fig:orbit_and_horizontal-a}]{\input{orbita_y_horizontal-2-a.pstex_t}}\hspace{10em}
    \subfigure[$\SG$ orbit through $q_0$ and nearby
    orbits\label{fig:orbit_and_horizontal-b}]{\input{orbita_y_horizontal-2-b.pstex_t}}
    \caption{$\SG$ orbit and $Hor_{\CD}^2(q_0)$ in $Q$}
    \label{fig:orbit_and_horizontal}
  \end{figure}
\end{proof}

\begin{remark}
  The existence of a level function as a datum for an affine discrete
  connection may seem a bit strange at first. However, if a
  decomposition like~\eqref{eq:pair_decomposition_Ad} is
  expected, by taking $q_1=q_0$, we see that $(q_0,l^Q_{g^{-1}}(q_0))
  \in Hor_{\CD}(q_0)$, so $\gamma(q_0):=g^{-1}$ is a level function for
  $\CD$.
\end{remark}

\begin{definition}
  Given an affine discrete connection $\CD$ we define its
  \emph{discrete connection $1$-form} $\CD:Q\times Q\rightarrow \SG$
  by $\CD(q_0,q_1):=g$ where $g$ is the element of $\SG$ that appears
  in decomposition~\eqref{eq:pair_decomposition_Ad}.
\end{definition}

\begin{remark}
  Notice that we will use the same letter to name both the connection
  $\CD$ and its discrete connection $1$-form. Also, calling $\CD$ a
  $1$-form can be misleading since it is not a $1$-form in the usual
  sense but, rather, a function. Still, the name is coming from the fact
  that when Leok et al. introduced it
  in~\cite{ar:leok_marsden_weinstein-a_discrete_theory_of_connections_on_principal_bundles}
  they also introduce the notion of ``discrete $k$-form'', and $\CD$
  is a discrete $1$-form in that sense.
\end{remark}

\begin{remark}\label{rem:CD_not_globally_defined}
  Due to Proposition~\ref{prop:CD_imp_pair_decomposition}, $\CD$ is
  only defined in a $\SG$-invariant open set $U\subset Q\times Q$. In
  what follows, we will abuse the notation and pretend that $\CD$ is
  defined everywhere not to make the notation too cumbersome.
\end{remark}

\begin{proposition}\label{prop:discrete_connection_vs_one_form}
  Let $\CD$ be an affine discrete connection on the principal
  bundle $\pi:Q\rightarrow Q/\SG$. Then, for all $(q_0,q_1)\in Q\times
  Q$ and $g_0,g_1\in \SG$,
  \begin{equation} 
    \label{eq:Ad_GxG} 
    \CD(l^Q_{g_0}(q_0),l^Q_{g_1}(q_1)) = g_1 \CD(q_0,q_1) g_0^{-1},
  \end{equation}
  as long as both sides are defined. Conversely, given a smooth
  function $\mathcal{A}:Q\times Q\rightarrow \SG$ such
  that~\eqref{eq:Ad_GxG} holds (with $\CD$ replaced by $\mathcal{A}$),
  then $Hor:=\{(q_0,q_1)\in Q\times Q: \mathcal{A}(q_0,q_1)=e\}$
  defines an affine discrete connection $\CD$ with level
  $\gamma(q):=\mathcal{A}(q,q)^{-1}$ and whose discrete connection
  $1$-form is $\mathcal{A}$.
\end{proposition}

\begin{proof}
  Writing $\CD(q_0,q_1)=g$ and $\CD(l^Q_{g_0}(q_0),l^Q_{g_1}(q_1))
  =\ti{g}$, by definition,
  \begin{gather*}
    (q_0,q_1) = (q_0,l^Q_g(q_0))\cdot
    \underbrace{(q_0,l^Q_{g^{-1}}(q_1))}_{\in
      Hor_{\CD}(q_0)}\\
    (l^Q_{g_0}(q_0),l^Q_{g_1}(q_1)) =
    (l^Q_{g_0}(q_0),l^Q_{\ti{g}}(l^Q_{g_0}(q_0))) \cdot
    \underbrace{(l^Q_{g_0}(q_0), l^Q_{\ti{g}^{-1}}(l^Q_{ g_1}(q_1)))}_{\in
      Hor_{\CD}(l^Q_{g_0}(q_0))}.
  \end{gather*}
  By the $\SG$-invariance of $Hor_{\CD}$, since
  $(q_0,l^Q_{g^{-1}}(q_1))\in Hor_{\CD}(q_0)$, it follows that
  $(l^Q_{g_0}(q_0),l^Q_{g_0 g^{-1}}(q_1))\in Hor_{\CD}(l^Q_{g_0}(q_0))$. We
  can write
 \begin{equation*}
    (l^Q_{g_0}(q_0),l^Q_{g_1}(q_1)) = (l^Q_{g_0}(q_0),
    l^Q_{g_1 g g_0^{-1}}(l^Q_{g_0}(q_0))) \cdot
    \underbrace{(l^Q_{g_0}(q_0),
      l^Q_{g_0 g^{-1}}(q_1))}_{\in Hor_{\CD}(l^Q_{g_0}(q_0))}.
  \end{equation*}
  Identity~\eqref{eq:Ad_GxG} follows, then, from the uniqueness of the
  decomposition~\eqref{eq:pair_decomposition_Ad}.

  The converse result is, mostly, routine checking and we will omit
  the details. That $Hor$ is a submanifold follows from the easily
  proved fact that $e\in \SG$ is a regular value of $\mathcal{A}$. The
  same proof shows that $p_2(Hor(q_0))$ is a submanifold of $Q$ of the
  right dimension. The $\SG$-equivariance of $\mathcal{A}$ shows that
  $Hor$ is $\SG$-invariant. An application of~\eqref{eq:Ad_GxG} leads to
  establishing the transversality condition for $Hor$. Finally, taking
  $\gamma(q):=\mathcal{A}(q,q)^{-1}$, one concludes that $Hor$ defines
  an affine discrete connection of level $\gamma$.
\end{proof}

Affine discrete connections have discrete horizontal lifts, just as
regular connections do.

\begin{definition}
  Let $\CD$ be an affine discrete connection on the principal bundle
  $\pi:Q\rightarrow Q/\SG$. The \emph{discrete horizontal lift}
  $\HLd{}:Q\times Q/\SG\rightarrow Q\times Q$ is given by
  \begin{equation*}
    \HLd{q_0}(r_1) := (q_0,q_1) \jiff (q_0,q_1)\in
    Hor_{\CD} \stext{ and } \pi(q_1)=r_1.
  \end{equation*}
  We define $\HLds{q_0} :=p_2\circ \HLd{q_0}$.
\end{definition}

\begin{remark}\label{rem:HLD_not_globally_defined}
  The map $\HLd{q_0}(r_1)$ is well defined, provided that $r_1$ is
  sufficiently close to $\pi(q_0)$, since for any
  $q_1\in\pi^{-1}(r_1)$,
  \begin{equation} \label{eq:horizontal_lift_using_Ad}
    \HLd{q_0}(r_1):= (q_0,l^Q_{\CD (q_0,q_1)^{-1}}(q_1)).
  \end{equation}
\end{remark}

\begin{lemma}
  Let $\CD$ be an affine discrete connection on $\pi: Q\rightarrow
  Q/\SG$ and consider the natural actions of $\SG$ on $Q\times
  (Q/\SG)$ (extending $l^Q$ trivially on the second factor) and the
  diagonal action on $Q\times Q$. Then, $\HLd{}$ is
  $\SG$-equivariant. 
\end{lemma}

\begin{proof}
  It is a direct computation
  using~\eqref{eq:horizontal_lift_using_Ad}.
\end{proof}

For completeness we mention that having a smooth and $\SG$-equivariant
map $\eta:Q\times (Q/\SG)\rightarrow Q\times Q$ such that $\pi\circ
p_2\circ \eta^q = id|_{Q/\SG}$ for all $q\in Q$ is equivalent to
having an affine discrete connection whose discrete horizontal lift is
$\eta$.

\begin{remark}
  Discrete connection forms $\CD$ and discrete horizontal lifts may
  not be defined everywhere. Indeed, if $\HLd{}:Q\times
  (Q/\SG)\rightarrow Q\times Q$ is defined everywhere, then for any
  $q\in Q$, the map $r\mapsto \HLds{q}(r)$ is a global section of the
  principal bundle $\pi:Q\rightarrow Q/\SG$, so that the bundle is
  trivial. Hence, for nontrivial principal bundles $\HLd{}$ and,
  consequently, $\CD$ can only be defined in some open set of the total
  space.

  On the other hand, if $Q$ has a $\SG$-invariant riemannian metric
  ---as it usually happens for mechanical systems--- the following
  construction provides a discrete affine connection. Let $\gamma$ be
  a level function and define
  \begin{equation*}
    Hor :=\{ (q_0,\exp_{\ti{\gamma}(q_0)}(v_1))\in Q\times Q  
    \stext{ for some } v_1\in (\mathcal{V}^{\SG}_{\ti{\gamma}(q_0)})^\perp 
    \subset T_{\ti{\gamma}(q_0)}Q\},
  \end{equation*}
  where $\exp$ and $\perp$ are those of the riemannian metric and
  $\ti{\gamma}(q_0) := l^Q_{\gamma(q_0)}(q_0)$. It can be checked that
  $Hor$ defines a discrete affine connection of level $\gamma$. Notice
  that the domain of the associated discrete affine connection form
  and discrete horizontal lifts are limited by the domain of the
  exponential mapping.
\end{remark}

\begin{example}\label{ex:continued_example-affine_discrete_connection}
  Consider $Q:=\R^2$ with the action of $\SG:=\R$ given by $l^Q_g(q)
  := (x,y+g)$, where $q=(x,y)$. Clearly, $p_1=\pi:Q\rightarrow Q/\SG$
  is a (trivial) principal bundle with structure group $\SG$. In order
  to define an affine discrete connection on this bundle we need to
  find a $\SG$-invariant manifold $Hor$ that is ``complementary'' to
  $\mathcal{V}^\SG_d = \{(q_0,q_1)\in Q\times Q : x_0=x_1\}$ in
  $Q\times Q$. We consider complements to $l^Q_{\SG}(\{q_0\})$ for
  $q_0\in Q$. A $G$-invariant family of curves that is complementary
  to the orbits is $Hor^{(2)}(q_0) = \{ q_1\in Q : y_1-y_0 =
  b(x_1+x_0)(x_1-x_0)/2\}$ for a parameter $b\in\R$. Since
  $l^Q_{\SG}(\{q_0\}) \cap Hor^{(2)}(q_0) =\{q_0\}$, the only thing we
  have to prove in order to see that $Hor = \{(q_0,q_1)\in Q\times Q:
  y_1-y_0 = b(x_1+x_0) (x_1-x_0)/2\}$ defines an affine discrete
  connection is that $T_{q_0}Q = T_{q_0}l^Q_{\SG}(\{q_0\}) \oplus
  T_{q_0}Hor^{(2)}(q_0)$, which is evident because
  $T_{q_0}l^Q_{\SG}(\{q_0\}) = \langle \del_y\big|_{q_0}\rangle$ and
  $T_{q_0}Hor^{(2)}(q_0) = \langle \del_x\big|_{q_0} + b x_0
  \del_y\big|_{q_0}\rangle$. We denote the discrete connection defined
  by $Hor$ with $\CDp{b}$.

  Since
  \begin{equation*}
    (q_0,q_1) = (q_0,l^Q_g(q_0)) \cdot \underbrace{(q_0,
      l^Q_{g^{-1}}(q_1))}_{\in Hor_{\CDp{b}}} \stext{ for } 
    g=y_1-y_0-\frac{1}{2}b(x_1+x_0) (x_1-x_0),
  \end{equation*}
  we have 
  \begin{equation}\label{eq:continued_example-discrete_connection_formulas}
    \begin{split}
      \CDp{b}(q_0,q_1) =& \;y_1-y_0-b(x_1+x_0)(x_1-x_0)/2\\
      \HLd{q_0}(r_1) =& \;
      (q_0,(r_1,y_0+b(r_1+x_0)(r_1-x_0)/2)).
    \end{split}
  \end{equation}
  The fact that $\CDp{b}$ and $\HLd{q_0}$ are defined everywhere is
  compatible with $Q\rightarrow Q/\SG$ being a trivial principal
  bundle.
\end{example}


\subsection{Isomorphisms associated to an affine discrete connection}
\label{sec:isomorphisms}

When working with symmetric discrete mechanical systems, one is led to
consider the space $(Q\times Q)/\SG$. It is convenient to have a
different model for this space. In this section we construct such a
model associated to an affine discrete connection $\CD$.

We start with a special case of a general construction called the
\emph{associated bundle} of a principal bundle (see Chap. 1,
Sect. 5 of~\cite{bo:kobayashi_nomizu-foundations-v1}). 

\begin{definition}
  Let $\pi:Q\rightarrow Q/\SG$ be a principal bundle and consider the
  action of $\SG$ on $Q\times \SG$ defined by $l_g^{Q\times
    \SG}(q,w):=(l_g^Q(q),l_g^\SG(w))$, with $l_g^\SG(w):=g w
  g^{-1}$. Being $\pi$ a principal bundle, the quotient
  $\ti{\SG}:=(Q\times \SG)/\SG$ by this action is a manifold, called
  the \emph{conjugate associated bundle}. The quotient map is denoted
  by $\rho:Q\times \SG\rightarrow \ti{\SG}$. The projections onto each
  of the two components of $Q\times \SG$ induce smooth maps
  $p^{Q/\SG}:\ti{G}\rightarrow Q/\SG$ and
  $p^{\SG/\SG}:\ti{\SG}\rightarrow \SG/\SG$. The first, $p^{Q/\SG}$,
  turns $\ti{\SG}$ into a bundle over $Q/\SG$ with fiber $\SG$.
\end{definition}

It is convenient to define
\begin{equation*}
  \ti{F}_1:Q\times \SG\times
  (Q/\SG)\rightarrow Q \stext{ with } 
  \ti{F}_1(q_0,w_0,r_1):=l^Q_{w_0}\big(\HLds{q_0}(r_1)\big).
\end{equation*}
Extending the $\SG$ action $l^{Q\times\SG}$ to $Q\times \SG\times
(Q/\SG)$ trivially on the last component and considering the $\SG$
action on $Q$, a simple computation using the $\SG$-equivariance of
$\HLd{}$ shows that $\ti{F}_1$ is $\SG$-equivariant.

\begin{proposition}\label{prop:isomorphisms_of_reduced_spaces}
  Given an affine discrete connection $\CD$ on $\pi:Q\rightarrow Q/\SG$,
  let $\tilde{{\Phi}}_{\CD}:Q\times Q \rightarrow Q\times \SG \times
  (Q/\SG)$ and $\tilde{{\Psi}}_{\CD} : Q\times \SG \times (Q/\SG)\rightarrow
  Q\times Q$ be defined by
  \begin{equation*}
    \begin{split}
      &\tilde{{\Phi}}_{\CD}(q_0,q_1) := (q_0,\CD(q_0,q_1), \pi(q_1)) \\
      &\tilde{{\Psi}}_{\CD}(q_0,w_0,r_1) :=
      (q_0,\ti{F}_1(q_0,w_0,r_1)).
    \end{split}
  \end{equation*}
  Then, $\tilde{{\Phi}}_{\CD}$ and $\tilde{{\Psi}}_{\CD}$ are smooth
  and mutually inverses (when restricted to the domain of
  $\CD$). Furthermore, considering the diagonal action of $\SG$ on
  $Q\times Q$ and $l^{Q\times \SG \times
    (Q/\SG)}_g(q_0,w_0,r_1):=(l^Q_g(q_0),l^\SG_g(w_0),r_1)$, both
  maps are $\SG$-equivariant, so that they induce diffeomorphisms
  $\Phi_{\CD}:(Q\times Q)/\SG \rightarrow \RS $ and $\Psi_{\CD} : \RS
  \rightarrow (Q\times Q)/\SG$.
\end{proposition}

\begin{proof}
  That $\ti{\Phi}_{\CD}$ and $\ti{\Psi}_{\CD}$ are smooth is clear
  from the definition. Checking that they are inverses and the
  $\SG$-equivariance is done by direct computation.
\end{proof}

\begin{remark}
  In fact, when $\CD$ and $\HLds{}$ are not globally defined the maps
  $\Phi_{\CD}$ and $\Psi_{\CD}$ are diffeomorphisms between open
  neighborhoods of of $\Gamma/\SG$ and $(Q\times\{e\})/\SG \times
  (Q/\SG)$.
\end{remark}

The commutative diagram~\eqref{eq:spaces_and_maps_reduction-1}, where
$\tilde{{\pi}}$ is the quotient map and $\Upsilon:= \Phi_{\CD}\circ
\tilde{{\pi}}$, shows some of the spaces and maps we have introduced
and that will be used later, when analyzing the behavior of symmetric
discrete mechanical systems.
\begin{equation} \label{eq:spaces_and_maps_reduction-1}
  \xymatrix{{Q\times Q} \ar[d]_{\ti{\pi}}
    \ar[r]^(.4){\ti{\Phi}_{\CD}}_(.4){\sim} \ar[dr]^{\Upsilon} & {Q\times
      \SG\times (Q/\SG)}\ar[d]^{\rho\times id} \\ {(Q\times Q)/\SG}
    \ar[r]_{\Phi_{\CD}}^{\sim} & {\RS}}
\end{equation}
.

\begin{lemma}
  \label{le:dgamma_computation} 
  Let $(q_0,q_1)\in Q\times Q$ and $(\delta q_0,\delta q_1)\in
  T_{(q_0,q_1)}(Q\times Q)$. Then,
  \begin{equation*}
    \begin{split}
      d\Upsilon(q_0,q_1)(\delta q_0, \delta q_1) =
      (d\rho(q_0,\CD(q_0,q_1))(\delta q_0,d\CD(q_0,q_1)(\delta
      q_0,\delta q_1)),d\pi(q_1)(\delta q_1)).
    \end{split}
  \end{equation*}
\end{lemma}

\begin{proof}
  Compute $d\Upsilon$, using $\Upsilon=(\rho\times id) \circ
  \ti{\Phi}_{\CD}$.
\end{proof}

\begin{example}\label{ex:continued_example-isomorphisms}
  In the context of
  Example~\ref{ex:continued_example-affine_discrete_connection}, we
  have that $Q=(Q/\SG) \times \SG$ with the $\SG$-action on $Q$
  corresponding to left multiplication on $\SG$. Then $\ti{\SG} \simeq
  (Q/\SG)\times \SG$ with $\rho((r_0,h_0),w_0) \mapsto (r_0,w_0)$. The
  map $\rho$ has a section $s(r_0,w_0) := ((r_0,0),w_0)$ (where we are
  using implicitly the isomorphism just defined). The isomorphisms
  that appear in Proposition~\ref{prop:isomorphisms_of_reduced_spaces}
  are
  \begin{gather*}
    \ti{\Phi}_{\CDp{b}}(q_0,q_1) =
    ((x_0,y_0),y_1-y_0-b(x_1+x_0)(x_1-x_0)/2,x_1)\\
    \ti{\Phi}_{\CDp{b}}^{-1}((r_0,h_0),w_0,r_1) =
    ((r_0,h_0),(r_1,h_0+b(r_1+r_0)(r_1-r_0)/2+w_0)).
  \end{gather*}
\end{example}

\begin{remark}
  The model spaces $Q\times \SG\times (Q/\SG)$ and $\RS$ are by no means
  unique. In fact, given an affine discrete connection $\CD$ it is
  possible to consider the isomorphism $\alpha_{\CD}:(Q\times
  Q)/\SG\rightarrow (Q/\SG\times Q/\SG) \times_{Q/\SG} \ti{\SG}$ defined by
  \begin{equation*}
    \alpha_{\CD}(\tilde{{\pi}}(q_0,q_1)) = ((\pi(q_0),\pi(q_1)),
    \rho(q_0,\CD(q_0,q_1))).
  \end{equation*}
  This isomorphism and the corresponding model space are very close to
  the data used in the analysis of the continuous case in
  Section~\ref{sec:reduction_of_symmetry_of_a_generalized_non_holonomic_mechanical_system}
  and is introduced
  by~\cite{ar:leok_marsden_weinstein-a_discrete_theory_of_connections_on_principal_bundles}
  in their proposed reduction of discrete unconstrained
  systems. Still, the presence of the fibered product in the model
  space makes it harder to work with compared to the simpler Cartesian
  product that appears in $\RS$.
\end{remark}


\section{Variations and reduced variations}
\label{sec:variations_and_reduced_variations}

In this section we analyze the relationship between the dynamics of a
symmetric discrete mechanical system on $Q\times Q$, which induces a
dynamics on (an open neighborhood of $Q\times\{e\}\times (Q/\SG)$ in)
$(Q\times Q)/\SG$ and the dynamics on (an open neighborhood of
$(Q\times\{e\})/\SG \times (Q/\SG)$ in) the isomorphic model $\RS$ of
$(Q\times Q)/\SG$. In order to describe the dynamics on this last
space we start by introducing some relevant notions.


\subsection{Reduced lagrangians}
\label{sec:reduced_lagrangians}

Given $(Q,L_d,\mathcal{D},\mathcal{D}_d)$ with a symmetry group $\SG$
and an affine discrete connection $\CD$ we construct a dynamical system
on $Q\times \SG\times (Q/\SG)$ using the isomorphism
$\tilde{{\Phi}}_{\CD}$. In particular, the dynamics will be determined
using the function $\check{L}_d:=L_d\circ \tilde{{\Psi}}_{\CD}$.
Since $L_d$ is $\SG$-invariant and, taking into account the symmetry
properties of the system, $\check{L}_d$ is $\SG$-invariant, $L_d$ and
$\check{L}_d$ induce maps on the corresponding quotient spaces. We
denote the map induced by $\check{L}_d$ on $\RS $ by
$\hat{L}_d$. Therefore,
\begin{equation*}
  \check{L}_d(q_0,w_0,r_1) = L_d(q_0,\ti{F}_1(q_0,w_0,r_1))
  \stext{and} \hat{L}_d(\rho(q_0,w_0),r_1) = \check{L}_d(q_0,w_0,r_1).
\end{equation*}
The following diagram shows all the relevant maps introduced so far.
\begin{equation*}
  \xymatrix{{} & {} & {\R} \\ {Q\times Q} \ar[urr]^{L_d} \ar[d]_{\ti{\pi}}
    \ar[r]_{\ti{\Phi}_{\CD}} \ar[dr]_{\Upsilon} & {Q\times \SG\times (Q/\SG)}
    \ar[ur]^(.3){\check{L}_d} \ar[d]^{\rho\times id} & {}\\ {(Q\times Q)/\SG}
    \ar[r]_{\Phi_{\CD}} & {\RS} \ar@/_/[uur]_{\hat{L}_d} & {}}
\end{equation*}
Associated to the reduced discrete lagrangian $\hat{L}_d$ we define a
reduced discrete action $\hat{S}_d$ by $\hat{S}_d(v_\cdot,r_\cdot) :=
\sum_k \hat{L}_d(v_k,r_{k+1})$.

\begin{example}\label{ex:continued_example-discrete_mechanical_system}
  Consider the discrete mechanical system
  $(Q,L_d,\mathcal{D},\mathcal{D}_d)$, where
  \begin{gather*}
    Q:=\R^2 \stext{ with points } q=(x,y),\\
    L_d(q_0,q_1) := m\big( (x_1-x_0)^2 +
    (y_1-y_0)^2\big)/2,\\
    \mathcal{D}_q:= \{\dot{x}\del_x\big|_q + \dot{y}\del_y\big|_q \in
    T_qQ : \dot{y} = x \dot{x}\} = \langle \del_x\big|_q - x
    \del_y\big|_q \rangle \subset T_qQ,\\
    \mathcal{D}_d:=\{(q_0,q_1)\in Q\times Q : y_1-y_0 =
    (x_1+x_0)(x_1-x_0)/2\},
  \end{gather*}
  that originates as a discretization of a free particle in $\R^2$
  subject to a nonholonomic constraint. This system can be readily
  solved using the discrete Lagrange--D'Alembert
  equations~\eqref{eq:dla-eqs}. Rather than following that path we
  will make a detailed study of its reduction in the following
  sections in order to illustrate some of the techniques developed so
  far.

  The group $\SG:=\R$ acts on $Q$ as considered in
  Example~\ref{ex:continued_example-affine_discrete_connection},
  turning $\SG$ into a symmetry group of
  $(Q,L_d,\mathcal{D},\mathcal{D}_d)$. Furthermore, for arbitrary
  $b\in\R$, $\CDp{b}$ defines isomorphisms $\ti{\Phi}_{\CDp{b}}$ and
  $\Phi_{\CDp{b}}$ that can be used to study the system. Since, by
  Example~\ref{ex:continued_example-isomorphisms}, $\RS\simeq
  (Q/\SG)\times \SG\times (Q/\SG)$, the relevant lagrangians are
  \begin{equation*}
    \begin{split}
      \check{L}_d((r_0,h_0),w_0,r_1) =& \;
      m\big( (r_1-r_0)^2 +
      (w_0+b(r_1+r_0)(r_1-r_0)/2)^2\big)/2\\
      \hat{L}_d(r_0,w_0,r_1) =& \;\check{L}_d((r_0,0),w_0,r_1).
    \end{split}
  \end{equation*}
\end{example}

Back in the general setting, we notice that, since $\hat{L}_d$ is
$\SG$-invariant and $Q\times \SG\times (Q/\SG)$ is a Cartesian product
(so that differentials on that space decompose in terms of
differentials on each of the spaces),
\begin{equation}  
  \label{eq:decomposition_reduced_extended_lagrangian}
  \begin{split}
    d\hat{L}_d(v_0,r_1)(\delta v_0,\delta r_1) =&
    d\hat{L}_d(\rho(q_0,w_0),r_1)(d\rho(q_0,w_0)(\delta q_0,\delta
    w_0),
    \delta r_1) \\
    =& d\check{L}_d(q_0,w_0,r_1)(\delta q_0,\delta w_0, \delta r_1) \\
    =& D_1\check{L}_d(q_0,w_0,r_1)(\delta q_0) +
    D_2\check{L}_d(q_0,w_0,r_1)(\delta w_0) \\
    & + D_3\check{L}_d(q_0,w_0,r_1)(\delta r_1).
  \end{split}
\end{equation}
In Section~\ref{sec:mixed_curvature_and_reduced_forces}, we study more
carefully the second term in this last equality.

\begin{lemma}\label{le:equivariance_of_DcL_d}
  Considering the $\SG$ actions $l^{Q\times \SG\times (Q/\SG)}$,
  $l^Q$, $l^\SG$ and their lifted actions to the corresponding tangent
  spaces, the following morphisms of vector bundles over $Q\times
  \SG\times (Q/\SG)$
  \begin{equation*}
    D_1\check{L}_d : p_1^*TQ\rightarrow \R \stext{and} 
    D_2\check{L}_d : p_2^*T\SG\rightarrow \R
  \end{equation*}
  are $\SG$-equivariant, where $p_1$ and $p_2$ are the projections of
  $Q\times \SG\times (Q/\SG)$ onto the corresponding factors.
\end{lemma}
\begin{proof}
  The $\SG$-invariance of $\check{L}_d$ leads to the
  $\SG$-equivariance of $d\check{L}_d:T(Q\times \SG\times
  (Q/\SG))\rightarrow\R$.  The statement follows by decomposing
  $T(Q\times \SG\times (Q/\SG)) = p_1^*TQ\oplus p_2^*T\SG\oplus
  p_3^*T(Q/\SG)$.
\end{proof}


\subsection{Mixed curvature and reduced forces}
\label{sec:mixed_curvature_and_reduced_forces}

Given a principal bundle $\pi:Q\rightarrow Q/\SG$, a connection
$\CC$ and an affine discrete connection $\CD$ we construct an object
that, somehow, compares the notions of continuous and discrete
horizontality introduced by $\CC$ and $\CD$.

\begin{definition}
  The \emph{mixed curvature} $\BM$ of $\CC$ and $\CD$ is the morphism
  of vector bundles over $Q\times Q$, $\BM: T(Q\times Q) \rightarrow
  \CD^*(T\SG)$ defined by
  \begin{equation*} 
    \BM(q_0,q_1)(\delta q_0, \delta q_1)
    := d\CD(q_0,q_1)(Hor_\CC(\delta q_0),Hor_\CC(\delta q_1)),
  \end{equation*}
  where $Hor_\CC(\delta q)$ denotes the $\CC$-horizontal part of
  $\delta q$. Taking into account the isomorphism $T(Q\times Q) \simeq
  p_1^*(TQ) \oplus p_2^*(TQ)$, the mixed curvature decomposes as
  $\BM=\BMp + \BMm $, where $\BMp:p_2^*TQ\rightarrow \CD^*(T\SG)$ and
  $\BMm:p_1^*TQ\rightarrow \CD^*(T\SG)$ are 
  \begin{equation*}
    \begin{split}
      \BMp(q_0,q_1)(\delta q_1) :=& \BM(q_0,q_1)(0,\delta q_1) =
      D_2\CD(q_0,q_1)(Hor_\CC(\delta q_1)) \\
      \BMm(q_0,q_1)(\delta q_0) :=& \BM(q_0,q_1)(\delta q_0,0) =
      D_1\CD(q_0,q_1)(Hor_\CC(\delta q_0)).
    \end{split}
  \end{equation*}
\end{definition}

\begin{remark}
  An obvious difference between the mixed curvature and the curvature
  of a connection is that the former can be represented as a $1$-form,
  while the latter is represented as a $2$-form.
\end{remark}

\begin{remark}
  The mixed curvature measures how $\CC$-horizontal deformations at
  $(q_0,q_1)$ depart from the level manifolds of $\CD$.
\end{remark}

\begin{example}\label{ex:continued_example-mixed_curvature}
  We continue in the context of
  Example~\ref{ex:continued_example-discrete_mechanical_system}. In
  order to compute the mixed curvature (and continue analyzing the
  reduction of the system) we introduce a connection on the principal
  bundle $\pi:Q\rightarrow Q/\SG$. For that purpose, the
  splitting~\eqref{eq:TQ_decomposition} is
  \begin{equation*}
    T_qQ = \underbrace{\{0\}}_{=\mathcal{W}_q}\oplus 
    \underbrace{\langle \del_y\big|_q\rangle }_{=\mathcal{U}_q}\oplus 
    \underbrace{\{0\}}_{=\mathcal{S}_q}\oplus 
    \underbrace{\langle \del_x\big|_q + x\del_y\big|_q\rangle.
    }_{=\mathcal{H}_q}
  \end{equation*}
  The connection $\CC$ with horizontal space
  $Hor_{\CC}=\mathcal{W}\oplus\mathcal{H}$ has connection $1$-form and
  horizontal lift given by
  \begin{equation*}
    \CC(\dot{x}\del_x\big|_q+\dot{y}\del_y\big|_q) =
    \dot{y}-\dot{x} x \stext{ and } 
    \HLc{q}(\del_r\big|_r) = \del_x\big|_q+r\del_y\big|_q.
  \end{equation*}
  We can now compute the mixed curvature $\BM$ associated to $\CC$
  and $\CDp{b}$:
  \begin{equation*}
    \begin{split}
      \BM(q_0,q_1)(\delta q_0,\delta q_1) =&
      d\CDp{b}(q_0,q_1)(Hor_{\CC}(\delta q_0),
      Hor_{\CC}(\delta q_1)) \\=&
      \del_{w_0}\big|_{\CDp{b}(q_0,q_1)} \otimes \big(dy_1\big|_{q_1}
      - dy_0\big|_{q_0} - b x_1 dx_1\big|_{q_1} + b x_0
      dx_0\big|_{q_0}\big)\\ & \phantom{==}\big( c_0
      (\del_{x_0}\big|_{q_0} + x_0\del_{y_0}\big|_{q_0}) + c_1
      (\del_{x_1}\big|_{q_1} + x_1\del_{y_1}\big|_{q_1}) \big) \\=&
      (1-b)(c_1 x_1 - c_0 x_0) \del_{w_0}\big|_{\CDp{b}(q_0,q_1)}.
    \end{split}
  \end{equation*}
  Notice that depending on the ``relative slopes'' $1$ and $b$ of
  $\CC$ and $\CDp{b}$, the mixed curvature vanishes or not.
\end{example}

Just as regular curvature, the mixed curvature can be seen as a
``reduced object'' defined on $(Q\times Q)/\SG$ or, using the
isomorphism $\Phi_{\CD}$, on $\RS$. We study this object next.

Since $\CD:Q\times Q\rightarrow \SG$ is $\SG$-equivariant for
$l^{Q\times Q}$ and $l^\SG$, considering the lifted actions on
$T(Q\times Q)$ and $TG$, $d\CD$ induces a morphism $\widetilde{d\CD}$
on the quotient vector bundles:
\begin{equation*}
  \xymatrix{ {T(Q\times Q)/\SG} \ar[d] \ar[r]^(0.6){\widetilde{d\CD}} & 
    {T\SG/\SG} \ar[d]
    \\ {(Q\times Q)/\SG} \ar[r]^(0.6){\ti{\CD}} & {\SG/\SG}}
\end{equation*}

Let $\sigma:\RS \rightarrow (Q/\SG)\times (Q/\SG)$ be
$\sigma(v_0,r_1):= (p^{Q/\SG}(v_0),r_1)$; then, we have the vector
bundles $\sigma^*T(Q/\SG\times Q/\SG)$ over $\RS$ and $(\rho\times
id)^*\sigma^*T(Q/\SG\times Q/\SG)$ over $Q\times \SG\times
(Q/\SG)$. We define $\check{\chi}:(\rho\times
id)^*\sigma^*T(Q/\SG\times Q/\SG) \rightarrow T(Q\times Q)$ by
$\check{\chi}(q_0,w_0,r_1)(\delta r_0, \delta r_1):= (\HLc{q_0}(\delta
r_0), \HLc{q_1}(\delta r_1))$ for $q_1:=\ti{F}_1(q_0,w_0,r_1)$. 
As $\ti{F}_1$ and the horizontal lift are $\SG$-equivariant,
$\check{\chi}$ descends to a morphism of vector bundles $\chi$:
\begin{equation*}
  \xymatrix{ {\sigma^*T(Q/\SG\times Q/\SG)} \ar[d] \ar[r]^(.55){\chi} &
    {T(Q\times Q)/\SG} \ar[d] 
    \\ {\RS} \ar[r]^{\Psi_\CD}& {(Q\times Q)/\SG} }
\end{equation*}
It will be convenient to consider the composition
$\widetilde{d\CD}\circ \chi$ as a morphism of the vector bundles
$\sigma^*T(Q/\SG\times Q/\SG)\rightarrow (\widetilde{\CD}\circ
\Psi_\CD)^*(T\SG/\SG)$ over the same base space $\RS$. Even better,
since $\widetilde{\CD}\circ \Psi_\CD = p^{\SG/\SG}\circ p_1$ we have the
following notion.

\begin{definition}
  The \jdef{reduced mixed curvature} $\BMR$ is the morphism of bundles
  over $\RS$ $\BMR:= \widetilde{d\CD}\circ \chi :
  \sigma^*T(Q/\SG\times Q/\SG)\rightarrow (p^{\SG/\SG}\circ
  p_1)^*(T\SG/\SG)$. Explicitly,
  \begin{equation*}  
    \begin{split}
      \BMR(v_0,r_1)(\delta r_0,\delta r_1) := &
      [\BM(q_0,q_1)(\HLc{q_0}(\delta r_0),\HLc{q_1} (\delta r_1))] \\
      =& [ d\CD (q_0,q_1)(\HLc{q_0}(\delta r_0),\HLc{q_1}(\delta
      r_1))],
    \end{split}
  \end{equation*}
  where $v_0 = \rho(q_0,w_0)$, $q_1:=\ti{F}_1(q_0,w_0,r_1)$ and
  $[\cdot]$ denotes the equivalence class in $T\SG/\SG$.  Associated
  to the decomposition $\BM = \BMp + \BMm$ there is a decomposition
  $\BMR = \BMRp + \BMRm$ where $\BMRp : p_2^*T(Q/\SG) \rightarrow
  (p^{\SG/\SG}\circ p_1)^* (T\SG/\SG)$ and $\BMRm :(p^{Q/\SG}\circ
  p_1)^*T(Q/\SG) \rightarrow (p^{\SG/\SG}\circ p_1)^* (T\SG/\SG)$ are
  defined by
  \begin{equation*}
    \begin{split}
      \BMRp(v_0,r_1)(\delta r_1) :=& [\BMp(q_0,q_1)(\HLc{q_1}(\delta
      r_1))] = [D_2\CD(q_0,q_1)(\HLc{q_1}(\delta r_1))] \\
      \BMRm(v_0,r_1)(\delta r_0) :=& [\BMm(q_0,q_1)(\HLc{q_0}(\delta
      r_0))] = [D_1\CD(q_0,q_1)(\HLc{q_0}(\delta r_0))].
    \end{split}
  \end{equation*}
\end{definition}

Now we return to the discrete lagrangians introduced in
section~\ref{sec:reduced_lagrangians}. In particular, we associate a
new notion to the second term in the last equality
of~\eqref{eq:decomposition_reduced_extended_lagrangian}.

By Lemma~\ref{le:equivariance_of_DcL_d}, $D_2\check{L}_d :
p_2^*T\SG\rightarrow \R$ is a $\SG$-equivariant morphism of vector
bundles over $Q\times \SG\times (Q/\SG)$. Thus, it induces a morphism
of the corresponding quotient vector bundles
$\widetilde{D_2\check{L}_d}$:
\begin{equation*}
  \xymatrix{ {(p_2^*T\SG)/\SG} \ar[r]^(.6){\widetilde{D_2\check{L}_d}}
    \ar[d] & {\R} \ar[dl]\\
    {\RS} & {} }
\end{equation*}
Notice that $(p_2^*T\SG)/\SG \simeq (p^{\SG/\SG}\circ p_1)^*
T\SG/\SG$. Then, $\widetilde{D_2\check{L}_d} \circ \BMR$ is a well
defined morphism of vector bundles.

\begin{definition}\label{def:reduced_discrete_force}
  The \emph{reduced discrete force} is the morphism of vector bundles
  over $\RS$ defined by $\hat{F}_d := \widetilde{D_2\check{L}_d} \circ
  \BMR: \sigma^*(T(Q/\SG\times Q/\SG))\rightarrow \R $ (essentially a
  $1$-form over $\RS$). Concretely,
  \begin{equation} \label{eq:reduced_discrete_force-def}
    \hat{F}_d(v_0,r_1)(\delta r_0, \delta r_1) :=
    D_2\check{L}_d(q_0,w_0,r_1) d\CD(q_0,q_1)(\HLc{q_0}(\delta r_0),
    \HLc{q_1}(\delta r_1))
  \end{equation}
  where $v_0=\rho(q_0,w_0)$ and $q_1:=\ti{F}_1(q_0,w_0,r_1)$. Once
  again, using $T(Q/\SG\times Q/\SG)\simeq p_1^*T(Q/\SG)\oplus
  p_2^*T(Q/\SG)$, we define $\hat{F}_d^+:=\widetilde{D_2\check{L}_d}
  \circ \BMRp$ and $\hat{F}_d^-:=\widetilde{D_2\check{L}_d} \circ
  \BMRm$, so that $\hat{F}_d=\hat{F}_d^++\hat{F}_d^-$. Explicitly,
  \begin{equation} \label{eq:reduced_mixed_curvature_pm-def}
    \begin{split}
      \hat{F}_d^+(v_0,r_1)(\delta r_1) :=&\;
      D_2\check{L}_d(q_0,w_0,r_1) D_2\CD(q_0,q_1)(\HLc{q_1}(\delta r_1)) \\
      \hat{F}_d^-(v_0,r_1)(\delta r_0) :=& \;D_2\check{L}_d(q_0,w_0,r_1)
      D_1\CD(q_0,q_1)(\HLc{q_0} (\delta r_0)).
    \end{split}
  \end{equation}
\end{definition}

\begin{example}\label{ex:continued_example-reduced_force}
  Continuing the analysis of the discrete mechanical system introduced
  in Example~\ref{ex:continued_example-discrete_mechanical_system}, we
  compute its reduced discrete force with respect to $\CC$ and
  $\CDp{b}$. If $(v_0,r_1) = (r_0,w_0,r_1)\in (Q/\SG)\times \SG \times
  (Q/\SG)$ we can take $q_0 := (r_0,0)$ and
  $q_1:=l^Q_{w_0}(\HLds{q_0}(r_1)) = (r_1,w_0+b(r_1^2-r_0^2)/2)$ in
  the computation. Since
  \begin{equation*}
    D_2\check{L}_d((r_0,0),w_0,r_1) = 
    m \big(w_0+ b (r_1^2-r_0^2)/2\big) dw_0\big|_{w_0},
  \end{equation*}
  using the results of
  Example~\ref{ex:continued_example-mixed_curvature}, we have
  \begin{equation*}
    \begin{split}
      \hat{F}_d(r_0,w_0,r_1)(c_0\del_{r_0}\big|_{r_0},c_1\del_{r_1}\big|_{r_1})
      =& \;\Bigl(m \big(w_0+ b (r_1^2-r_0^2)/2\big)
      dw_0\big|_{w_0} \Bigr)\\ & \phantom{==} \big((1-b)(c_1 r_1- c_0
      r_0) \del_{w_0}\big|_{w_0}\big) \\=& \; m (1-b) (c_1 r_1- c_0 r_0)
      \big(w_0+ b (r_1^2-r_0^2)/2\big).
    \end{split}
  \end{equation*}
  Equivalently,
  \begin{equation*}
    \hat{F}_d(r_0,w_0,r_1) = m(1-b) \big(w_0+ b
    (r_1^2-r_0^2)/2\big) (r_1 dr_1-  r_0dr_0).
  \end{equation*}
\end{example}


\subsection{Reduced dynamics}
\label{sec:reduced_dynamics}

In this section we relate the dynamics of a symmetric discrete
mechanical system to the dynamics of a ``reduced'' system, defined
using a variational principle.

\begin{remark}\label{rem:group_notation}  
  The analysis of the symmetric system leads us to consider some objects
  defined on the symmetry group $\SG$, that is a Lie group. We review
  very briefly a notation that is convenient and used in this
  subject. When $w_0, w_1\in \SG$ and $\delta w_0\in T_{w_0}\SG$ we
  define
  \begin{equation*}
    w_1 \delta w_0 := dL_{w_1}(w_0)(\delta w_0)\in T_{w_1 w_0}\SG \stext{and}  
    \delta w_0 w_1 := dR_{w_1}(w_0)(\delta w_0)\in T_{w_0 w_1}\SG,
  \end{equation*}
  where $L_{w_1}$ and $R_{w_1}$ denote the multiplication by $w_1$ on
  the left and right, respectively. Analogously, when $\alpha_0 \in
  T_{w_0}^*\SG$,
  \begin{equation*}
    \begin{split}
      w_1 \alpha_0 :=& (dL_{w_1^{-1}}(w_1w_0))^*(\alpha_0)\in
      T_{w_1 w_0}^*\SG
      \stext{and}\\
      \alpha_0 w_1 :=& (dR_{w_1^{-1}}(w_0w_1))^*(\alpha_0)\in
      T^*_{w_0w_1}\SG.
    \end{split}
  \end{equation*}
  Last, we notice that if $\alpha_0\in T_{w_0}^*\SG$ and $\delta
  w_0\in T_{w_0}\SG $, the following identities hold
  \begin{equation*}
    \alpha_0(\delta w_0) = (w_1 \alpha_0)(w_1 \delta w_0) = (\alpha_0
    w_1)(\delta w_0 w_1).
  \end{equation*}
\end{remark}

\begin{theorem} \label{th:4_points-general} Let $q_\cdot$ be a
  discrete curve in $Q$, $r_k:=\pi(q_k)$, $w_k:=\CD(q_k,q_{k+1})$ and
  $v_k:=\rho(q_k,w_k)$ be the corresponding discrete curves in $Q/\SG$,
  $\SG$ and $\ti{\SG}$. Then, given a discrete mechanical system
  $(Q,L_d,\mathcal{D},\mathcal{D}_d)$ with symmetry group $\SG$, the
  following statements are equivalent.
  \begin{enumerate}
  \item \label{it:var_pple-general} $(q_k,q_{k+1})\in \mathcal{D}_d$
    for all $k$ and $q_\cdot$ satisfies the variational principle
    $dS_d(q_\cdot)(\delta q_\cdot) = 0$ for all vanishing end points
    variations $\delta q_\cdot$ such that $\delta q_k \in
    \mathcal{D}_{q_k}$ for all $k$.

  \item \label{it:eq_lda-general} $q_\cdot$ satisfies the
    Lagrange--D'Alembert equations~\eqref{eq:dla-eqs} for all $k$.

  \item \label{it:red_var_pple-general} $(v_k,r_{k+1})\in
    \hat{\mathcal{D}}_d := \Phi_{\CD}(\mathcal{D}_d/\SG)$ for all $k$
    and $d\hat{S}_d(r_\cdot,v_\cdot) (\delta r_\cdot, \delta v_\cdot)
    = 0$ for all variations $(\delta v_\cdot, \delta r_\cdot)$ with
    vanishing end points such that $\delta r_k\in
    \hat{\mathcal{D}}_{r_k} :=d\pi(q_k)(\mathcal{D}_{q_k})$ and
    \begin{equation} \label{eq:delta_vk-def}
      \begin{split}
        \delta v_k :=& \; d\rho(q_k,w_k)\big(\HLc{q_k}(\delta r_k),
        d\CD(q_k,q_{k+1})(\HLc{q_k}(\delta r_k),\HLc{q_{k+1}}(\delta
        r_{k+1}))\big)
        \\
        & + d\rho(q_k,w_k)\big((\xi_k)_Q(q_k),
        d\CD(q_k,q_{k+1})((\xi_k)_Q(q_k),(\xi_{k+1})_Q(q_{q+1}))\big),
      \end{split}
    \end{equation}
    where $(q_k,\xi_k)\in\jgsg^{\mathcal{D}}$.

  \item \label{it:red_lda_eq-general}
    $(v_k,r_{k+1})\in\hat{\mathcal{D}}_d$ for all $k$ and
    $(v_\cdot,r_\cdot)$ satisfies the following conditions for each
    fixed $(v_{k-1},r_k,v_k,r_{k+1})$.

    \begin{itemize}
    \item $\phi \in T_{r_k}^*(Q/\SG)$ defined by
      \begin{equation} \label{eq:red_hor_gen-Q}
        \begin{split}
          \phi :=& D_1 \check{L}_d(q_k,w_k,r_{k+1}) \circ \HLc{q_k}
          + D_3 \check{L}_d(q_{k-1},w_{k-1},r_{k}) \\
          & + \hat{F}_d^-(v_k,r_{k+1}) + \hat{F}_d^+(v_{k-1},r_k)
        \end{split}
      \end{equation}
      vanishes on $\hat{\mathcal{D}}_{r_k}$, \emph{i.e.,}
      \begin{equation} \label{eq:cond_red_hor_gen-Q}
        \phi\in\hat{\mathcal{D}}_{r_k}^\annihilator.
      \end{equation}

    \item $\psi \in \jgsg^*$ defined by
      \begin{equation} \label{eq:red_vert_gen-Q}
        \begin{split}
          \psi := & D_2\check{L}_d(q_{k-1},w_{k-1},r_{k}) w_{k-1}^{-1}
          - D_2\check{L}_d(q_k,w_k,r_{k+1}) w_k^{-1}
        \end{split}
      \end{equation}
      vanishes on $\jgsg^{\mathcal{D}}_{q_k}$, \emph{i.e.,}
      \begin{equation} \label{eq:cond_red_vert-gen-Q}
        \psi\in(\jgsg^{\mathcal{D}}_{q_k})^\annihilator.
      \end{equation}
    \end{itemize}
  \end{enumerate}
\end{theorem}

\begin{proof}
  First we study the equivalence of conditions $(q_k,q_{k+1}) \in
  \mathcal{D}_d$ and $(v_k,r_{k+1})\in\hat{\mathcal{D}}_d$. If
  $(q_k,q_{k+1})\in \mathcal{D}_d$, it follows that
  $\ti{\pi}(q_k,q_{k+1})\in \mathcal{D}_d/\SG$ and $(v_k,r_{k+1}) =
  \Phi_{\CD}(\ti{\pi}(q_k,q_{k+1})) \in
  \Phi_{\CD}(\mathcal{D}_d/\SG)$. Conversely, if $(v_k,r_{k+1}) =
  \Phi_{\CD}(\ti{\pi}(q_k,q_{k+1})) \in\hat{\mathcal{D}}_d$, we have
  that $\ti{\pi}(q_k,q_{k+1})\in \mathcal{D}_d/\SG$, so that
  $(q_k,q_{k+1}) \in \mathcal{D}_d$.

  Next we tackle the equivalence between the variational principles
  and equations.

  \ref{it:var_pple-general} $\jiff$ \ref{it:eq_lda-general}. Is a
  standard computation using calculus of variations.

  \ref{it:red_var_pple-general} $\imp$
  \ref{it:var_pple-general}. Given a variation $\delta q_\cdot$ as in
  the statement, we let $\delta r_k := d\pi(q_k)(\delta q_k)$ and
  $\delta v_k$ using~\eqref{eq:delta_vk-def} for $\xi_k
  :=\CC(q_k)(\delta q_k)$, so that $\delta q_k = \HLc{q_k}(\delta r_k)
  + (\xi_k)_Q(q_k)$. Using Lemma~\ref{le:dgamma_computation} to compute
  $d\Upsilon$,
  \begin{equation}\label{eq:3_imp_1}
    \begin{split}
      dS_d(q_\cdot)(\delta q_\cdot) =& \sum\nolimits_{k=0}^{N-1}
      dL_d(q_k,q_{k+1})(\delta q_k,\delta q_{k+1}) \\=&
      \sum\nolimits_{k=0}^{N-1} d\hat{L}_d(\Upsilon(q_k,q_{k+1}))
      d\Upsilon(\delta q_k,\delta q_{k+1}) \\=&
      \sum\nolimits_{k=0}^{N-1} d\hat{L}_d(v_k,r_{k+1})(\delta v_k,
      \delta r_{k+1}) = d\hat{S}_d(r_\cdot,v_\cdot)(\delta
      r_\cdot,\delta v_\cdot) = 0,
    \end{split}
  \end{equation}
  where the last equality holds because $(r_\cdot,v_\cdot)$ satisfies
  condition~\ref{it:red_var_pple-general}. Hence, $q_\cdot$ satisfies
  condition~\ref{it:var_pple-general}.

  \ref{it:var_pple-general} \imp \ref{it:red_var_pple-general}. Given
  a variation $(\delta r_\cdot, \delta v_\cdot)$ as in the statement,
  let $\delta q_k := \HLc{q_k}(\delta r_k) + (\xi_k)_Q(q_k)$. Then,
  using Lemma~\ref{le:dgamma_computation} and the explicit form
  of~\eqref{eq:delta_vk-def}, we have that $(\delta v_k, \delta
  r_{k+1}) = d\Upsilon(q_k,q_{k+1})(\delta q_k, \delta q_{k+1})$. A
  computation similar to~\eqref{eq:3_imp_1} shows that
  \begin{equation*}
    d\hat{S}_d(v_\cdot, r_\cdot)(\delta v_\cdot, \delta r_\cdot) = 
    dS_d(q_\cdot)(\delta q_\cdot) = 0.
  \end{equation*}
  since the variation $\delta q_\cdot$ satisfies
  condition~\ref{it:var_pple-general}. Hence,
  condition~\ref{it:red_var_pple-general} holds.

  \ref{it:red_var_pple-general} \jiff \ref{it:red_lda_eq-general}.
  For variations $(\delta v_\cdot, \delta r_\cdot)$, writing $\delta
  v_k = d\rho(q_k,w_k)(\delta q_k,\delta w_k)$,
  \begin{equation*}
    \begin{split}
      d\hat{S}_d(v_\cdot,r_\cdot)&(\delta v_\cdot, \delta r_\cdot) =
      \sum\nolimits_{k=0}^{N-1} d\hat{L}_d(v_k,r_{k+1})(\delta v_k,\delta
      r_{k+1}) \\=& \sum\nolimits_{k=0}^{N-1}
      d\check{L}_d(q_k,w_k,r_{k+1})(\delta q_k, \delta w_k, \delta
      r_{k+1}) \\=&
      \sum\nolimits_{k=0}^{N-1} \big( D_1 \check{L}_d(q_k,w_k,r_{k+1})(\delta
      q_k) + D_2 \check{L}_d(q_k,w_k,r_{k+1})(\delta w_k) \\ &+ D_3
      \check{L}_d(q_k,w_k,r_{k+1})(\delta r_{k+1})\big).
    \end{split}
  \end{equation*}
  Using $\CC$ to decompose $\delta q_k = \HLc{q_k}(\delta r_k) +
  (\xi_k)_Q(q_k)$ and rearranging indexes, for vanishing end point
  variations, we obtain
  \begin{equation}\label{eq:tiSd_computation}
    \begin{split}
      d\hat{S}_d(v_\cdot,r_\cdot)&(\delta v_\cdot, \delta r_\cdot) 
      = \sum_{k=0}^{N-1} \Bigl( D_1
      \check{L}_d(q_k,w_k,r_{k+1}) \circ \HLc{q_k} + D_3
      \check{L}_d(q_{k-1},w_{k-1},r_{k}) \\ & + D_2
      \check{L}_d(q_k,w_k,r_{k+1}) D_1\CD(q_k,q_{k+1})\circ \HLc{q_k} \\
      & + D_2 \check{L}_d(q_{k-1},w_{k-1},r_{k})
      D_2\CD(q_{k-1},q_{k})\circ \HLc{q_k}\Bigr)(\delta r_k) \\ & +
      \sum_{k=0}^{N-1}\bigl(D_1 \check{L}_d(q_k,w_k,r_{k+1})
      ((\xi_k)_Q(q_k)) \\ & \phantom{\sum} +
      D_2\check{L}_d(q_k,w_k,r_{k+1})
      d\CD(q_k,q_{k+1})((\xi_k)_Q(q_k),(\xi_{k+1})_Q(q_{k+1}))\bigr).
    \end{split}
  \end{equation}
  Since the variations $\delta r_\cdot$ are independent from those
  generated by the $\xi_\cdot$, we conclude that condition
  $d\hat{S}_d(r_\cdot,v_\cdot)(\delta r_\cdot,\delta v_\cdot)=0$ from
  point~\ref{it:red_var_pple-general} in the statement is equivalent
  to the vanishing of the first and second summations
  in~\eqref{eq:tiSd_computation} independently, for all vanishing end
  points variations $\delta r_\cdot$ with $\delta
  r_k\in\hat{\mathcal{D}}_{r_k}$ for all $k$ and for all $\xi_\cdot$
  with $(q_k,\xi_k)\in \jgsg^\mathcal{D}$ for all $k$ and
  $\xi_0=\xi_n=0$.

  Recalling that the variations $\delta r_\cdot$ are independent
  and~\eqref{eq:reduced_mixed_curvature_pm-def}, it is clear that
  the vanishing of the first summation in~\eqref{eq:tiSd_computation} is
  equivalent to the fact that $\phi$, defined
  in~\eqref{eq:red_hor_gen-Q}, satisfies
  condition~\eqref{eq:cond_red_hor_gen-Q}.

  Before we prove the equivalence of
  condition~\eqref{eq:cond_red_vert-gen-Q} and the
  vanishing of the second summation in~\eqref{eq:tiSd_computation} we need
  two auxiliary computations. On the one hand,
  \begin{equation*}
    \begin{split}
      d\CD(q_k,q_{k+1})((\xi_k)_Q(q_k),&(\xi_{k+1})_Q(q_{k+1})) \\=& \;
       \frac{d}{dt}\bigg|_{t=0}\big( \exp(t\xi_{k+1}) \CD(q_k,
      q_{k+1})\exp(-t\xi_k) \big) \\=& \; \xi_{k+1}w_k - w_k \xi_k,
    \end{split}
  \end{equation*}
  where the notation used in the last equality is the one introduced
  in Remark~\ref{rem:group_notation}. On the other hand, since
  $\check{L}_d(l^Q_g(q_k),w_k,r_{k+1}) =
  \check{L}_d(q_k,l^\SG_{g^{-1}}(w_k),r_{k+1})$,
  \begin{equation*}
    \begin{split}
      D_1\check{L}_d(q_k,w_k,r_{k+1})((\xi_k)_Q(q_k)) =& \;
      \frac{d}{dt}\bigg|_{t=0}
      \check{L}_d(l^Q_{\exp(t\xi_k)}(q_k),w_k,r_{k+1}) \\=& \;
      D_2\check{L}_d(q_k,w_k,r_{k+1})\frac{d}{dt}\bigg|_{t=0}
      (\exp(-t\xi_k)w_k\exp(t\xi_k))\\=& \;
      D_2\check{L}_d(q_k,w_k,r_{k+1}) (-\xi_kw_k + w_k\xi_k).
    \end{split}
  \end{equation*}
  Using the previous computations we see that the vanishing of the
  second summation in~\eqref{eq:tiSd_computation} is equivalent to
  \begin{equation*}
    \begin{split}
      0 =& \sum\nolimits_{k=1}^{N-1} \big( D_2\check{L}_d(q_k,w_k,r_{k+1})
      (\xi_{k+1}w_k-\xi_kw_k) \big) \\=& \sum\nolimits_{k=1}^{N-1} \big(
      D_2\check{L}_d(q_{k-1},w_{k-1},r_{k})w_{k-1}^{-1} -
      D_2\check{L}_d(q_k,w_k,r_{k+1}) w_k^{-1} \big)(\xi_k)
    \end{split}
  \end{equation*}
  for all $\xi_\cdot$ with $(q_k,\xi_k)\in \jgsg^\mathcal{D}$. It
  follows immediately that this last condition is equivalent to the
  fact that $\psi$, defined in~\eqref{eq:red_vert_gen-Q},
  meets condition~\eqref{eq:cond_red_vert-gen-Q}.
\end{proof}

From the proof of Theorem~\ref{th:4_points-general} we isolate the
following partial result.

\begin{lemma}
  \label{le:vertical_part_of_th_4_points_general} 
  With the notation of Theorem~\ref{th:4_points-general}, the
  following assertions are equivalent.
  \begin{enumerate}
  \item The discrete curve $q_\cdot$ satisfies $dS_d(q_\cdot)(\delta
    q_\cdot^\mathcal{S})=0$ for all vanishing end points variations
    such that $\delta q_k\in \mathcal{S}_{q_k}$ for all $k$.
  \item The curve $(v_\cdot, r_\cdot)$ satisfies
    condition~\eqref{eq:cond_red_vert-gen-Q} for
    $\psi$ defined by~\eqref{eq:red_vert_gen-Q}.
  \end{enumerate}
\end{lemma}

\begin{remark}
  We notice that the horizontal equations that appear in
  item~\ref{it:red_lda_eq-general} of
  Theorem~\ref{th:4_points-general} contain force terms that, because
  of~\eqref{eq:reduced_discrete_force-def} are a composition of
  derivatives of the reduced lagrangian and a term involving $\CC$ and
  $\CD$. Since in the continuous setting, the forces that appear in
  the analogous equation involve the reduced curvature of the
  nonholonomic connection, we chose to call the corresponding term in
  the discrete setting the reduced mixed curvature. In this case, we
  recover the standard result that the vanishing of the mixed
  curvature implies the vanishing of the discrete forces in the
  reduced system. Other than that, we do not have a good reason to
  call $\BM$ or $\BMR$ ``curvatures''.
\end{remark}

\begin{remark}
  For classical mechanical systems, the Lagrange--D'Alembert
  equations~\eqref{eq:generalized_lagrange_dalambert} are second order
  differential equations, while their reduced counterparts that appear
  in Theorem~\ref{th:generalized_reduction-continuous} are second
  order in the $Q/\SG$ variables but only first order in the
  $\ti{\jgsg}$ ones. Analogously, for discrete mechanical systems, the
  discrete Lagrange--D'Alembert equations~\eqref{eq:dla-eqs} are
  recurrence equations of second order while the reduced equations
  obtained in item~\ref{it:red_lda_eq-general} of
  Theorem~\ref{th:4_points-general} are second order in the $Q/\SG$
  variables but only first order in the $\ti{\SG}$ ones.
\end{remark}

\begin{remark}
  Theorem~\ref{th:4_points-general} is similar in spirit to the
  reduction theorems for classical mechanical systems, like
  Theorem~\ref{th:generalized_reduction-continuous}. Still there is a
  noticeable technical difference between both types of results. In
  the continuous case, even in the generalized context, the choice of
  a connection on the principal bundle $\pi$ serves the dual purpose
  of constructing a model for the reduced space via the diffeomorphism
  $\alpha_\CC$ and determines a horizontal / vertical splitting of the
  variational principle and equations of motion. In the discrete
  context a continuous connection serves the same splitting purpose
  but a discrete connection $\CD$ is used to construct the model
  reduced space. It would be interesting to see if there is any
  advantage in using two different connections in the reduction of
  continuous mechanical systems.
\end{remark}

\begin{example}\label{ex:continued_example-reduced_equations}
  Continuing the analysis of the discrete mechanical system introduced
  in Example~\ref{ex:continued_example-discrete_mechanical_system}, we
  apply Theorem~\ref{th:4_points-general} to find the equations of
  motion of the reduced system. The reduced variational constraint is
  \begin{equation*}
    \begin{split}
      \hat{\mathcal{D}}_{r_k} =& d\pi((r_k,0))(\mathcal{D}_{(r_k,0)})
      = d\pi((r_k,0))(\langle \del_x\big|_{(r_k,0)} - r_k
      \del_y\big|_{(r_k,0)} \rangle) \\=& \langle \del
      r\big|_{r_k}\rangle = T_{r_k}(Q/\SG),
    \end{split}
  \end{equation*}
  and, since
  \begin{equation*}
    (r_0,w_0,r_1) \in \hat{\mathcal{D}}_d \jiff 
    \ti{\Phi}_{\CDp{b}}^{-1}((r_0,0),w_0,r_1) \in \mathcal{D}_d \jiff 
    w_0 = (1-b)(r_1^2-r_0^2)/2,
  \end{equation*}
  the reduced kinematic constraints are
  \begin{equation*}
    \hat{\mathcal{D}}_d = \{(r_0,w_0,r_1)\in (Q/\SG)\times \SG\times (Q/\SG) : 
    w_0 = (1-b)(r_1^2-r_0^2)/2\}.
  \end{equation*}

  Next we compute the reduced equations of motion. Notice that
  $\mathcal{S}=\{0\}$ for this system, so that
  $\jgsg^\mathcal{D}=\{0\}$ and
  condition~\eqref{eq:cond_red_vert-gen-Q} is
  trivially satisfied, that is, there are no vertical equations.

  In order to find the horizontal equations we compute
  \begin{gather*}
    D_1\check{L}_d((r_k,h_k),w_k,r_{k+1}) = -m\big( (r_{k+1}-r_{k}) +
    (w_{k}+b(r_{k+1}^2-r_{k}^2)/2 ) b r_{k}\big)
    dx_k\big|_{(r_k,h_k)}\\
    D_3\check{L}_d((r_{k-1},h_{k-1}),w_{k-1},r_{k}) = m\big(
    (r_{k}-r_{k-1}) + (w_{k-1}+b(r_{k}^2-r_{k-1}^2)/2 ) b r_{k}\big)
    dr_k\big|_{r_k},
  \end{gather*}
  recall from Example~\ref{ex:continued_example-mixed_curvature}
  that
  \begin{equation*}
    \HLc{(r_k,h_k)} = \big( \del_x\big|_{(r_k,h_k)} + 
    r_k\del_y\big|_{(r_k,h_k)} \big) \otimes dr_k\big|_{r_k},
  \end{equation*}
  and, from Example~\ref{ex:continued_example-reduced_force} the
  expression of the reduced forces.
  \begin{equation*}
    \begin{split}
      \hat{F}_d^-(r_k,w_k,r_{k+1}) =& -m(1-b) \big(w_k+ b
      (r_{k+1}^2-r_{k}^2)/2\big) r_k dr_k\big|_{r_k}\\
      \hat{F}_d^+(r_{k-1},w_{k-1},r_{k}) =& \; m(1-b) \big(w_{k-1}+ b
      (r_{k}^2-r_{k-1}^2)/2\big) r_k dr_k\big|_{r_k}.
    \end{split}
  \end{equation*}
  Then, $\phi = m U_k dr_k\big|_{r_k}$ for 
  \begin{equation*}
    \begin{split}
      U_k :=& 
      -\big((r_{k+1}-r_{k})-(r_{k}-r_{k-1}) \big) \\ & 
      -r_k
      \big(w_k-w_{k-1} + b
      ((r_{k+1}^2-r_{k}^2)-(r_{k}^2-r_{k-1}^2))/2\big),
    \end{split}
  \end{equation*}
  and, since $\hat{\mathcal{D}} = T(Q/\SG)$,
  condition~\eqref{eq:cond_red_hor_gen-Q} says
  that $U_k=0$, which is the horizontal equation of motion for the
  reduced system. Thus, the reduced evolution is determined by the
  system of equations
  \begin{equation}\label{eq:continued_example-red_evol_eqs-2}
    U_k = 0 \stext{ and } w_k = (1-b)(r_{k+1}^2-r_{k}^2)/2.
  \end{equation}
  Using the second equation to eliminate the $w$ dependence in the
  first expression we obtain
  \begin{equation}\label{eq:continued_example-reduced_equations-r}
    0 = \big((r_{k+1}-r_{k})-(r_{k}-r_{k-1}) \big)+ 
    r_k ((r_{k+1}^2-r_{k}^2)-(r_{k}^2-r_{k-1}^2))/2.
  \end{equation}
  From this equation of degree two in $r_{k+1}$ the evolution of this
  variable is obtained and, using the second equation
  of~\eqref{eq:continued_example-red_evol_eqs-2} the dynamics of $w_k$
  is determined.
\end{example}

\begin{remark}
  Notice that, depending on the value of the parameter $b$, the
  reduced system constructed in
  Example~\ref{ex:continued_example-reduced_equations} is forced or
  not. Besides making the reduced system unforced, the value $b=1$
  also gives a very simple dynamics to the $w_\cdot$ variables. This
  will be a characteristic that we explore in more detail in
  Section~\ref{sec:reduced_equations_of_motion-chaplygin} when we
  consider systems of Chaplygin type.
\end{remark}

Back in the general setting,
conditions~\eqref{eq:cond_red_hor_gen-Q}
and~\eqref{eq:cond_red_vert-gen-Q} in
Theorem~\ref{th:4_points-general} establish the equations of motion of
the reduced system. However, they are explicitly written in terms of
$q_\cdot$. We will see in
Section~\ref{sec:intrinsic_version_of_the_reduced_equations_of_motion}
that, in fact, they can be defined in an intrinsic manner, in terms of
objects defined on the reduced space.


\section{Intrinsic version of the reduced equations of motion}
\label{sec:intrinsic_version_of_the_reduced_equations_of_motion}

In this section we write the horizontal and vertical equations in
terms of the reduced system. In fact, the equations of motion will be
given as conditions on morphisms of vector bundles on the second order
reduced manifold $\RSsec := \RS\times_{Q/\SG}\RS$, where the fibered
product is taken over the maps $p_2:\RS\rightarrow Q/\SG$ and
$p^{Q/\SG}\circ p_1:\RS\rightarrow Q/\SG$.

It is convenient to consider the space $\RSsecC:=Q\times
\SG\times (Q/\SG) \times_{Q/\SG}\RS$ with the $\SG$-action
$l^{\RSsecC}_g(q_0,w_0,r_1,v_1,r_2):=(l^Q_g(q_0),l^\SG_g(w_0),r_1,v_1,r_2)$, so
that $\RSsecC/\SG=\RSsec$. Additionally, we define the maps
\begin{equation*}
  F_1:\RSsecC\rightarrow Q \stext{by} 
  F_1(q_0,w_0,r_1,v_1,r_2) := \ti{F}_1(q_0,w_0,r_1),
\end{equation*}
and $F_2:\RSsecC\rightarrow \SG$ by
\begin{equation*}
  F_2(q_0,w_0,r_1,v_1,r_2):= l^\SG_{\tau(\ti{F}_1(q_0,w_0,r_1),\ti{q}_1)}(\ti{w}_1) 
  \stext{ if } v_1 = \rho(\ti{q}_1,\ti{w}_1),
\end{equation*}
where $\tau:Q\times_{Q/\SG} Q\rightarrow\SG$ is defined by $\tau(q,q')=g$
if $l^Q_g(q')=q$. It is easy to check that $F_1$ and $F_2$ are
$\SG$-equivariant.


\subsection{Horizontal equations}
\label{sec:horizontal_equations-intrinsic}

In order to give an intrinsic meaning
to~\eqref{eq:cond_red_hor_gen-Q}, we consider the
following commutative diagram,
\begin{equation}\label{eq:diagram_intrinsic_horizontal_eqs_incomplete}
  \xymatrix{{} & {\R} \ar[dl] & {p_3^*T(Q/\SG)} \ar[dll] 
    \ar[l]_(.6){\check{\phi}} & 
    {T(Q/\SG)} \ar[d]\\
    {\RSsecC}  \ar[rrr]^{p_3}&  {} & {} & {Q/\SG} }
\end{equation}
where $\check{\phi}$ is defined by
\begin{equation*}
  \begin{split}
    \check{\phi}(q_0,w_0,r_1,v_1,r_2,\delta r_1) :=& \big(D_1
    \check{L}_d(q_1,w_1,r_2)\circ \HLc{q_1} +
    D_3\check{L}_d(q_0,w_0,r_1)\\ & + \hat{F}_d^-(v_1,r_{2}) +
    \hat{F}_d^+(\rho(q_0,w_0),r_1)\big)(\delta r_1)
  \end{split}
\end{equation*}
for $q_1:=\ti{F}_1(q_0,w_0,r_1)$ and
$w_1:=F_2(q_0,w_0,r_1,v_1,r_2)$. Also, $\SG$ acts on $p_3^*T(Q/\SG)$
by $l^{p_3^*T(Q/\SG)}_g(q_0,w_0,r_1,v_1,r_2,\delta r_1) :=
(l^Q_g(q_0),l^\SG_g(w_0),r_1,v_1,r_2,\delta r_1)$ and trivially on
$\R$.

\begin{lemma}\label{le:check_phi_is_SG_equivariant}
  The map $\check{\phi}$ is a $\SG$-equivariant morphism of vector
  bundles.
\end{lemma}
\begin{proof}
  From the explicit definition, $\check{\phi}$ is a morphism of vector
  bundles. 

  In order to check the $\SG$-equivariance, for
  $(q_0,w_0,r_1,v_1,r_2)\in \RSsecC$, we let $q_1:=\ti{F}_1(q_0,w_0,r_1)$
  and $w_1:=F_2(q_0,w_0,r_1,v_1,r_2)$. Then, writing
  \begin{equation*}
    \begin{split}
      \check{\phi}_1(q_0,w_0,r_1,v_1,r_2,\delta r_1) :=& \;\big(D_1
      \check{L}_d(q_1,w_1,r_2)\circ \HLc{q_1} +
      D_3\check{L}_d(q_0,w_0,r_1)\big)(\delta r_1) \\
      \check{\phi}_2(q_0,w_0,r_1,v_1,r_2,\delta r_1) :=& \;
      \big(\hat{F}_d^-(v_1,r_{2}) +
      \hat{F}_d^+(\rho(q_0,w_0),r_1)\big)(\delta r_1),
    \end{split}
  \end{equation*}
  we have $\check{\phi} =\check{\phi}_1 +\check{\phi}_2$. From the
  definition, it is immediate that $\check{\phi}_2$ is
  $\SG$-equivariant, so we concentrate on $\check{\phi}_1$.
  \begin{equation}\label{eq:equivariance_c_phi_1-aux}
    \begin{split}
      \check{\phi}_1&(l^{p_3^*T(Q/\SG)}_g(q_0,w_0,r_1,v_1,r_2,\delta r_1)) = 
      \check{\phi}_1(l^Q_g(q_0),l^\SG_g(w_0),r_1,v_1,r_2,\delta r_1)
      \\&= D_1 \check{L}_d(l^Q_g(q_1),l^\SG_g(w_1),r_2)\circ
      \HLc{l^Q_g(q_1)}(\delta r_1) +
      D_3\check{L}_d(l^Q_g(q_0),l^\SG_g(w_0),r_1)(\delta r_1),
    \end{split}
  \end{equation}
  where we used the $\SG$-equivariance of $\ti{F}_1$ and $F_2$ in the
  first term of the last identity. By
  Lemma~\ref{le:equivariance_of_DcL_d}, we have $D_1 \check{L}_d
  (l^Q_g (q_1),l^\SG_g(w_1),r_2) l^{TQ}_g = D_1 \check{L}_d
  (q_1,w_1,r_2)$, so that
  \begin{equation*}
    \begin{split}
      D_1 \check{L}_d (l^Q_g (q_1),l^\SG_g(w_1),r_2) \circ
      \HLc{l^Q_g(q_1)} = & \;D_1 \check{L}_d (l^Q_g
      (q_1),l^\SG_g(w_1),r_2) \circ l_g ^{TQ} \circ \HLc{q_1}\\
       = & \; D_1 \check{L}_d (q_1,w_1,r_2) \circ \HLc{q_1}.
    \end{split}
  \end{equation*}
  Also, from $\check{L}_d (l^Q_g(q_0),l^\SG_g(w_0),r_1) = \check{L}_d
  (q_0,w_0,r_1)$, we obtain
  \begin{equation*}
    D_3\check{L}_d(l^Q_g(q_0),l^\SG_g(w_0),r_1) = D_3\check{L}_d(q_0,w_0,r_1). 
  \end{equation*}
  Replacing the last two identities back
  in~\eqref{eq:equivariance_c_phi_1-aux}, the $\SG$-equivariance of
  $\check{\phi}_1$ follows.
\end{proof}

By Lemma~\ref{le:check_phi_is_SG_equivariant}, $\check{\phi}$ defines
a morphism of vector bundles
$\bar{\phi}:(p_3^*T(Q/\SG))/\SG\rightarrow \R$. Since
$(p_3^*T(Q/\SG))/\SG=p_2^*T(Q/\SG)$, where $p_2:\RSsec\rightarrow
Q/\SG$ is the projection, we have $\bar{\phi}:p_2^*T(Q/\SG)\rightarrow
\R$. Concretely, if $(q_0,w_0)\in \rho^{-1}(v_0)$,
\begin{equation*}
  \begin{split}
    \bar{\phi}(v_0,r_1,v_1,r_2,\delta r_1) =
    \check{\phi}(q_0,w_0,r_1,v_1,r_2, \delta r_1).
  \end{split}
\end{equation*}

We extend
diagram~\eqref{eq:diagram_intrinsic_horizontal_eqs_incomplete} to the
following commutative diagram, where $\hat{\mathcal{D}}\subset
T(Q/\SG)$ is the subbundle introduced in
Theorem~\ref{th:4_points-general}.
\begin{equation*}
  \xymatrix{{} & {\R} \ar[dl] & {p_3^*T(Q/\SG)} \ar[dll] 
    \ar[l]_(.6){\check{\phi}} & 
    {\hat{\mathcal{D}}} \ar@{^{(}->}[r]  & {T(Q/\SG)} \ar[d]\\
    {\RSsecC} \ar[dd] \ar[rrrr]^{p_3}& {} & {} & {} & {Q/\SG} \\
    {} & {\R} \ar[dl] & {p_2^*T(Q/\SG)} \ar[dll] \ar[l]_(.6){\bar{\phi}} & 
    {p_2^*\hat{\mathcal{D}}} \ar@{_{(}->}[l] & {}\\
    {\RSsec} \ar@/_1.3cm/[rrrruu]_{p_2} & {} & {} & {} & {}}
\end{equation*}

\begin{proposition}\label{prop:horizontal_condition_of_motion-intrinsic}
  Condition~\eqref{eq:cond_red_hor_gen-Q} in
  Theorem~\ref{th:4_points-general} is equivalent to
  \begin{equation}\label{eq:cond_red_hor_gen-intrinsic}
    \bar{\phi}\mid_{p_2^*(\hat{\mathcal D})} = 0.
  \end{equation}
\end{proposition}
\begin{proof}
  For $v_k$ and $r_k$ defined as in Theorem~\ref{th:4_points-general},
  the explicit form of $\bar{\phi}$ coincides with that of $\phi$
  defined in~\eqref{eq:red_hor_gen-Q} and the
  vanishing
  condition~\eqref{eq:cond_red_hor_gen-Q}
  coincides
  with~\eqref{eq:cond_red_hor_gen-intrinsic}.
\end{proof}


\subsection{Vertical equations}
\label{sec:vertical_equations-intrinsic}

Consider the commutative diagram
\begin{equation}\label{eq:diagram_intrinsic_vertical_eqs_incomplete}
  \xymatrix{ {\R} \ar[dr] & {F_1^*(Q\times\jgsg)} \ar[d] 
    \ar[l]_(.6){\check{\psi}} & {F_1^* \jgsg^\mathcal{D}} \ar@{_{(}->}[l] \ar[dl]
    & {\jgsg^\mathcal{D}} \ar@{^{(}->}[r] 
    & {Q\times \jgsg} \ar[d]\\
    {} & {\RSsecC} \ar[rrr]_{F_1} & {} & {} & {Q}}
\end{equation}
where 
\begin{equation*}
  \check{\psi}(q_0,w_0,r_1,v_1,r_2,\xi_1) := 
  \big(D_2\check{L}_d(q_{0},w_{0},r_{1}) w_{0}^{-1}
  - D_2\check{L}_d(q_1,w_1,r_{2}) w_1^{-1}\big) (\xi_1)
\end{equation*}
for $q_1:=\ti{F}_1(q_0,w_0,r_1)$ and
$w_1:=F_2(q_0,w_0,r_1,v_1,r_2)$. 

Notice that $F_1^*(Q\times \jgsg) = \RSsecC\times \jgsg$. In addition,
$\SG$ acts on $\R$ trivially and on $\RSsecC\times \jgsg$ by
$l^{\RSsecC\times \jgsg}_g(q_0,w_0,r_1,v_1,r_2,\xi_1) :=
(l^Q_g(q_0),l^\SG_g(w_0),r_1,v_1,r_2,l^{T\SG}_g(\xi_1))$.

\begin{lemma}\label{le:check_psi_is_SG_equivariant}
  The map $\check{\psi}$ is a $\SG$-equivariant morphism of vector
  bundles.
\end{lemma}
\begin{proof}
  From the explicit definition, $\check{\psi}$ is a morphism of vector
  bundles. 

  In order to check the $\SG$-equivariance, for
  $(q_0,w_0,r_1,v_1,r_2)\in \RSsecC$, we let
  $q_1:=\ti{F}_1(q_0,w_0,r_1)$ and $w_1:=F_2(q_0,w_0,r_1,v_1,r_2)$.
  Since
  \begin{equation*}
    \check{\psi}(q_0,w_0,r_1,v_1,r_2,\xi_1) = 
    D_2\check{L}_d(q_{0},w_{0},r_{1})(\xi_1 w_{0})
    - D_2\check{L}_d(q_1,w_1,r_{2}) (\xi_1 w_1),
  \end{equation*}
  using the $\SG$-equivariance of $\ti{F}_1$ and $F_2$, we have
  \begin{equation}\label{eq:check_psi_SG_equivariance-aux-0}
    \begin{split}
      \check{\psi}(l^{\RSsecC\times \jgsg}_g(q_0,w_0,r_1,v_1,r_2,\xi_1)) = &
      \check{\psi}(l^Q_g(q_0),l^\SG_g(w_0),r_1,v_1,r_2,l^{T\SG}_g(\xi_1))
      \\=& D_2\check{L}_d(l^{Q}_g(q_{0}),
      l^\SG_g(w_{0}),r_{1})(l^{T\SG}_g(\xi_1) l^\SG_g(w_{0})) \\ &-
      D_2\check{L}_d(l^Q_g(q_1), l^{\SG}_g(w_1),r_{2})
      (l^{T\SG}_g(\xi_1) l^\SG_g(w_1)).
    \end{split}
  \end{equation}
  On one hand, recalling the notation of
  Remark~\ref{rem:group_notation}, we have
  \begin{equation}\label{eq:check_psi_SG_equivariance-aux-1}
    l^{T\SG}_g(\xi_1) l^\SG_g(w_{0}) = g \xi_1 g^{-1} g w_0 g^{-1} = 
    g \xi_1 w_0 g^{-1} = dl^{\SG}_g(w_0)(\xi_1 w_0).
  \end{equation}
  On the other, the $\SG$-equivariance of $D_2\check{L}_d$ proved in
  Lemma~\ref{le:equivariance_of_DcL_d} together
  with~\eqref{eq:check_psi_SG_equivariance-aux-1}, applied
  to~\eqref{eq:check_psi_SG_equivariance-aux-0} lead to the
  $\SG$-equivariance of $\check{\psi}$.
\end{proof}

By Lemma~\ref{le:check_psi_is_SG_equivariant}, $\check{\psi}$ defines
a morphism of vector bundles $\bar{\psi}:(\RSsecC\times\jgsg)/\SG
\rightarrow \R$ over $\RSsec$. Explicitly, for any
$(q_0,w_0)\in\rho^{-1}(v_0)$ and $\xi_1\in\jgsg$,
\begin{equation}\label{eq:bar_psi-explicit}
  \begin{split}
    \bar{\psi}([(q_0,w_0,r_1,v_1,r_2,\xi_1)]) =&
    \check{\psi}(q_0,w_0,r_1,v_1,r_2,\xi_1) \\=&
    \big(D_2\check{L}_d(q_{0},w_{0},r_{1}) w_{0}^{-1} -
    D_2\check{L}_d(q_1,w_1,r_{2}) w_1^{-1}\big) (\xi_1),
  \end{split}
\end{equation}
where $q_1:=\ti{F}_1(q_0,w_0,r_1)$ and
$w_1:=F_2(q_0,w_0,r_1,v_1,r_2)$.

A simple computation shows that, because $\mathcal{D}\subset TQ$ is
$\SG$-invariant, $F_1^*\jgsg^\mathcal{D}\subset \RSsecC\times \jgsg$
is $\SG$-invariant. Therefore $(F_1^*\jgsg^\mathcal{D})/\SG \subset
(\RSsecC\times\jgsg)/\SG$ is a vector subbundle.

We collect the different objects in the following commutative diagram,
which is the quotient of (part
of)~\eqref{eq:diagram_intrinsic_vertical_eqs_incomplete}
\begin{equation*}
  \xymatrix{ {\R} \ar[dr] & {(\RSsecC\times \jgsg)/\SG} \ar[d] 
    \ar[l]_(.65){\bar{\psi}} & 
    {(F_1^* \jgsg^\mathcal{D})/\SG} \ar@{_{(}->}[l] \ar[dl]\\ 
    {} & {\RSsec} & {} }
\end{equation*}

\begin{proposition}\label{prop:vertical_condition_of_motion-intrinsic}
  Condition~\eqref{eq:cond_red_vert-gen-Q} in
  Theorem~\ref{th:4_points-general} is equivalent to
  \begin{equation}\label{eq:cond_red_vert_gen-intrinsic}
    \bar{\psi}\mid_{(F_1^*\jgsg^\mathcal{D})/\SG}\; = 0.
  \end{equation}
\end{proposition}
\begin{proof}
  For $v_k$ and $r_k$ defined as in Theorem~\ref{th:4_points-general},
  the explicit form of $\bar{\psi}$ coincides with that of $\psi$
  defined in~\eqref{eq:red_vert_gen-Q} and the vanishing
  condition~\eqref{eq:cond_red_vert-gen-Q} coincides
  with~\eqref{eq:cond_red_vert_gen-intrinsic}.
\end{proof}

\begin{corollary}\label{cor:4_points_general-intrinsic}
  Any one of the four equivalent conditions of
  Theorem~\ref{th:4_points-general} is equivalent to $(v_k,r_{k+1})
  \in \hat{\mathcal D}_d$ and
  conditions~\eqref{eq:cond_red_hor_gen-intrinsic}
  and~\eqref{eq:cond_red_vert_gen-intrinsic} are
  met for all $k$.
\end{corollary}


\section{Reconstruction}
\label{sec:reconstruction}

Given a discrete curve $q_\cdot$ in $Q$ and its image
$(v_\cdot,r_\cdot)$ in $\RS$, Theorem~\ref{th:4_points-general}
establishes an equivalence between $q_\cdot$ being a trajectory of the
original system (\emph{i.e.}, a solution of the original dynamics) and
$(v_\cdot,r_\cdot)$ being a solution of the reduced dynamics.  In this
section we study the reconstruction problem, that is, given a curve
$(v_\cdot,r_\cdot)$ in $\RS$ that satisfies adequate conditions
(constraints and equations of motion), is it possible to find a
trajectory $q_\cdot$ of the original system that projects to
$(v_\cdot,r_\cdot)$?

Consider the following construction. Given a discrete curve
$(v_\cdot,r_\cdot) \in\RS$ and $q_k \in Q$ (one value of $k$) such
that $\pi(q_k) = r_k = p^{Q/\SG}(v_k)$, if
$v_k=\rho(\tilde{q}_k,\tilde{w}_k) $ we define
\begin{equation}\label{eq:reconstruction-recursive_formula}
  u_k := l^\SG_{\tau(q_k,\ti{q}_k)}(\ti{w}_k) \in \SG \stext{ and } 
  q_{k+1} := \ti{F}_1(q_k,u_k,r_{k+1}).
\end{equation}

Since $l^\SG_{\tau(q_k,l^Q_g(\tilde{q_k}))}(l^\SG_{g}(\tilde{w_k})) =
l^\SG_{\tau(q_k,\tilde{q_k})g^{-1}}(l^\SG_{g}(\tilde{w_k}))=
l^\SG_{\tau(q_k,\tilde{q_k})}(\tilde{w_k})$, the previous construction
is independent of the chosen representatives of $v_k$. Therefore,
given $q_0\in Q$ with $\pi(q_0)=r_0$, applying the construction
iteratively defines a unique discrete curve $q_\cdot$ in $Q$.  The
following result establishes the properties of $q_\cdot$.

\begin{theorem}\label{th:reconstruction-general}
  Consider $(\bar{q}_0,\bar{q}_1) \in \mathcal{D}_d$ and
  $(v_\cdot,r_\cdot) \in \RS$ such that $\pi(\bar{q}_0)=r_0$,
  $\pi(\bar{q}_1)=r_1$,
  $v_0=\rho(\bar{q}_0,\CD(\bar{q}_0,\bar{q}_1))$. If
  $(v_\cdot,r_\cdot)$ satisfies $(v_k,r_{k+1}) \in \hat{\mathcal D}_d$
  and
  conditions~\eqref{eq:cond_red_hor_gen-intrinsic}
  and~\eqref{eq:cond_red_vert_gen-intrinsic} for
  all $k$, then the discrete curve $q_\cdot$ constructed
  by~\eqref{eq:reconstruction-recursive_formula} from $\bar{q}_0$ is a
  trajectory of the system on $Q$ whose image by $\Upsilon$ is the curve
  $(v_\cdot,r_\cdot)$ and satisfies $q_0=\bar{q}_0$ and $q_1=\bar{q}_1$.
\end{theorem}

\begin{proof}
  The curve $q_\cdot$ is a lifting of the curve $(v_\cdot,r_\cdot)$ to
  $Q$; that is, $r_k = \pi(q_k)$ and $v_k =
  \rho(q_k,\CD(q_k,q_{k+1}))$. Indeed,
  from~\eqref{eq:reconstruction-recursive_formula}, $\pi(q_{k}) =
  \pi(l^Q_{w_{k-1}} (\HLds{q_{k-1}}(r_{k}))) =
  \pi(\HLds{q_{k-1}}(r_{k})) = r_{k}$. On the other hand,
  \begin{equation*}
    \CD(q_k,q_{k+1}) = \CD(q_k,l^Q_{u_k} (\HLds{q_k}(r_{k+1}))) =
    u_k \CD\underbrace{(q_k,\HLds{q_k}(r_{k+1}))}_{\in
      Hor_{\CD}} = u_k,
  \end{equation*}
  so that
  \begin{equation*}
    \begin{split}
      \rho(q_k,\CD(q_k,q_{k+1})) = & \;\rho(q_k,u_k) = \rho(q_k,
      l^\SG_{\tau(q_k,\ti{q}_k)}(\ti{w}_k)) \\=& \;
      \rho(l^Q_{\tau(\ti{q}_k,q_k)}(q_k), \ti{w}_k) =
      \rho(\ti{q}_k,\ti{w}_k)=v_k,
    \end{split}
  \end{equation*}
  which completes the argument, that is $\Upsilon(q_k,q_{k+1}) =
  (v_k,r_{k+1})$ for all $k$.

  The lifted curve satisfies the first initial condition because, by
  construction, $r_0$ is lifted to $\bar{q}_0$. Since
  $v_0=\rho(\bar{q}_0,\CD(\bar{q}_0,\bar{q}_1))$, according
  to~\eqref{eq:reconstruction-recursive_formula},
  $u_0=\CD(\bar{q}_0,\bar{q}_1)$ so that $r_1$ is lifted to
  $q_1=l^Q_{\CD(q_0,q_1)}(\HLds{q_0}(r_1))$.
  Using~\eqref{eq:horizontal_lift_using_Ad}, we conclude
  that $q_1=\bar{q}_1$.

  Next, we see that $(q_k,q_{k+1})\in\mathcal{D}_d$ for all $k$. By
  the $\SG$-invariance of $\mathcal{D}_d$,
  \begin{equation*}
    \begin{split}
      (q_k,q_{k+1})\in \mathcal{D}_d \; \jiff & \; \ti{\pi}(q_k,q_{k+1})\in
      \mathcal{D}_d/\SG \; \jiff \;\Phi_{\CD}(\ti{\pi}(q_k,q_{k+1})) \in
      \Phi_{\CD}(\mathcal{D}_d/\SG) \\\jiff & \;(v_k,r_{k+1}) \in
      \hat{\mathcal{D}}_d,
    \end{split}
  \end{equation*}
  which holds by hypothesis for all $k$.

  The only thing left to do is to check that $q_\cdot$ is a trajectory
  of the discrete mechanical system. By hypothesis,
  $(v_\cdot,r_\cdot)$ satisfies
  conditions~\eqref{eq:cond_red_hor_gen-intrinsic}
  and~\eqref{eq:cond_red_vert_gen-intrinsic}. Then,
  by Corollary~\ref{cor:4_points_general-intrinsic}
  conditions~\eqref{eq:cond_red_hor_gen-Q}
  and~\eqref{eq:cond_red_vert-gen-Q} hold. But,
  since the relationship among $q_\cdot$, $v_\cdot$, $r_\cdot$ and
  $w_\cdot$ is precisely that of the statement of
  Theorem~\ref{th:4_points-general} and we proved that
  condition~\ref{it:red_lda_eq-general} holds, we conclude that
  $q_\cdot$ satisfies condition~\ref{it:var_pple-general} in
  Theorem~\ref{th:4_points-general}, hence it is a trajectory of the
  system on $Q$.
\end{proof}

\begin{example}\label{ex:continued_example-reconstruction}
  The last step to complete our analysis of the system
  $(Q,L_d,\mathcal{D},\mathcal{D}_d)$ introduced in
  Example~\ref{ex:continued_example-discrete_mechanical_system} is to
  consider the reconstruction of the evolution of the original system
  given a trajectory $(r_\cdot,w_\cdot)$ of the reduced system
  compatible with some initial data $(\bar{q}_0,\bar{q}_1)\in
  \mathcal{D}_d$.

  According to Theorem~\ref{th:reconstruction-general} we
  use~\eqref{eq:reconstruction-recursive_formula} to construct the
  trajectory $q_\cdot$. Since $v_k = (r_k,w_k) = \rho((r_k,0),w_k)$
  and $\tau((x_k,y_k),(r_k,0)) = y_k$ (remember that $x_k =
  \pi(x_k,y_k) = r_k$),
  \begin{gather*}
    u_k = l^\SG_{y_k}(w_k) = w_k,\\
    q_{k+1} = (x_{k+1},y_{k+1}) =
    l^Q_{w_k}(\HLds{(x_k,y_k)}(r_{k+1})) =
    (r_{k+1},y_k+b(r_{k+1}^2-x_k^2)/2+w_k)
  \end{gather*}
  that expresses $q_\cdot$ in terms of the known data
  $(r_\cdot,w_\cdot)$. Simplifying we obtain
  \begin{equation*}
    (x_k,y_k) = (r_k,\bar{y}_0+(r_k^2-r_0^2)/2) \stext{ for all } k,
  \end{equation*}
  with $r_0=\bar{x}_0,\, r_1=\bar{x}_1.$ Notice that the resulting
  trajectory is independent of the parameter $b$ chosen to do the
  reduction, as it should be.
\end{example}

\begin{remark}
  The reduction and reconstruction techniques developed so far can be
  applied in the case where the configuration space is the symmetry
  group (acting by left multiplication), a problem that has already
  been studied
  in~\cite{ar:fedorov_zenkov-discrete_nonholonomic_ll_systems_on_lie_groups,ar:fedorov_zenkov-dynamics_of_the_discrete_chaplygin_sleigh,ar:mclachlan_perlmutter-integrators_for_nonholonomic_mechanical_systems}. In
  this context, Theorems~\ref{th:4_points-general}
  and~\ref{th:reconstruction-general} allow us to re derive Theorem 3
  and Corollary 4
  of~\cite{ar:mclachlan_perlmutter-integrators_for_nonholonomic_mechanical_systems}.
\end{remark}


\section{Nonholonomic discrete momentum}
\label{sec:nonholonomic_discrete_momentum}

For discrete or continuous holonomic mechanical systems the presence
of continuous symmetries automatically leads to the existence of
conserved quantities (momenta), due to Noether's Theorem. In the
nonholonomic case, this is no longer true, essentially due to the
behavior of constraint forces. Instead, in this case, one obtains an
equation that describes the evolution of momenta over the trajectory
of the system. Below we discuss the relationship of the discrete
momentum evolution equation and the dynamics of a discrete mechanical
system with symmetries. 

\begin{definition}
  Given $(Q,L_d,\mathcal{D},\mathcal{D}_d)$ with symmetry group $\SG$,
  the \jdef{nonholonomic discrete momentum map} is the application
  $J_d : Q\times Q \rightarrow (\jgsg^\mathcal{D})^*$ defined by
  \begin{equation}\label{eq:non_holonomic_discrete_momentum-def}
    J_d(q_0,q_1)(q_0,\xi) := -D_1L_d(q_0,q_1)\xi_Q(q_0).
  \end{equation}
  for all $(q_0,\xi) \in \jgsg^\mathcal{D}$.
\end{definition}

Given any section $\ti{\xi} \in \Gamma(\jgsg^\mathcal{D})$, define the
map $(J_d)_{\ti{\xi}}:Q\times Q \rightarrow \R$ by
\begin{equation*}
  (J_d)_{\ti{\xi}}(q_0,q_1) := J_d(q_0,q_1)(q_0,\ti{\xi}(q_0)).
\end{equation*}

Straightforward computations give the following result.
\begin{lemma}
  For $\xi \in \jgsg$ and $(q_0,q_1) \in Q \times Q$, 
  \begin{equation}
    \label{eq:D1=-D2}
    D_1 L_d (q_0,q_1)(\xi_Q (q_0)) = -D_2 L_d
    (q_0,q_1)(\xi_Q(q_1)).
  \end{equation}
  Also, if $(q_0,w_0,r_1) = \ti{\Phi}_{\CD}(q_0,q_1) = (q_0,{\mathcal
    A}_d(q_0,q_1),\pi(q_1))$, then
  \begin{equation}\label{eq:D1_in_Q-D_2_in_RS}
    -D_1L_d(q_0,q_1)(\xi_Q (q_0)) = (D_2 \check{L}_d(q_0,w_0,r_1)w_0^{-1})(\xi).
  \end{equation}
  In particular, $J_d(q_0,q_1) = D_2\check{L}_d
  (q_0,w_0,r_1)w_0^{-1}$.
\end{lemma}

\begin{remark}
  By definition, when $Q\times Q$ and $\jgsg^\mathcal{D}$ are seen as
  bundles over $Q$ with respect to the projection on the first
  variable, $J_d$ is a bundle mapping (that is, it maps fibers to
  fibers). J. Cort\'es defines $J_d^{nh}$
  in~\cite[pp. 154-5]{bo:cortes-non_holonomic} in a slightly different
  way: $J_d^{nh}(q_0,q_1)(q_1,\xi) := D_2L_d(q_0,q_1)\xi_Q(q_1)$,
  which is a bundle map when the projection on the second variable is
  considered in $Q\times Q$. In any case, due to~\eqref{eq:D1=-D2},
  both maps are, essentially, the same. Also $(J_d)_{\ti{\xi}} =
  (J_d^{nh})_{\ti{\xi}}$, for all $\ti{\xi}\in
  \Gamma(\jgsg^\mathcal{D})$.
\end{remark}

A trajectory of a discrete mechanical system is determined by the
Discrete Lagrange--D'Alembert Principle
(Definition~\ref{def:discrete_lagrange_dalembert_principle}). When the
variational constraints are decomposed by $\mathcal D = \mathcal S
\oplus \mathcal H$, it is possible to decompose all admissible
variations into horizontal and vertical variations, in the sense that
they belong to $\mathcal{H}$ or $\mathcal{S}$. It is also possible to
decompose the variational principle accordingly. The following result,
whose proof is obvious, makes a precise statement.

\begin{proposition}
  Let $q_\cdot$ be a discrete curve in $Q$. Then, the following
  conditions are equivalent.
  \begin{enumerate}
  \item $q_\cdot$ satisfies the variational part of the discrete
    Lagrange--D'Alembert principle
    (Definition~\ref{def:discrete_lagrange_dalembert_principle}).
  \item $q_\cdot$ satisfies 
    \begin{equation*}
      dS_d(q_\cdot)(\delta q_\cdot^\mathcal{S}) = 0 \stext{ and } 
      dS_d(q_\cdot)(\delta q_\cdot^\mathcal{H}) = 0
    \end{equation*}
    for all pairs of variations vanishing at the end points $\delta
    q_\cdot^\mathcal{S}$ and $\delta q_\cdot^\mathcal{H}$ in
    $\mathcal{S}$ and $\mathcal{H}$ respectively.
  \end{enumerate}
\end{proposition}

The next result relates the vertical variational principle, the
evolution of the nonholonomic discrete momentum and
condition~\eqref{eq:cond_red_vert_gen-intrinsic}.

\begin{theorem}\label{th:vertical_conditions_equivalence}
  In the context of the theorem~\ref{th:4_points-general}, if
  $q_\cdot$ is a discrete curve in $Q$, let $r_k = \pi(q_k), w_k =
  \CD(q_k,q_{k+1})$ and $v_k = \rho(q_k,w_k)$. Then, the following
  conditions are equivalent.
  \begin{enumerate}
  \item \label{it:discrete_non_holonomic_momentum-var_pple}
    $q_\cdot$ satisfies the vertical variational principle. That is,
    $dS_d(q_\cdot)(\delta q_\cdot^\mathcal{S}) =0$ for all vanishing
    end point variation $\delta q_k^\mathcal{S} \in \mathcal{S}_{q_k}$
    for all $k$.
  \item \label{it:discrete_non_holonomic_momentum-eqs_in_r}
    Condition~\eqref{eq:cond_red_vert_gen-intrinsic}
    is satisfied.
  \item \label{it:discrete_non_holonomic_momentum-momentum_eq} For all
    sections $\ti{\xi}\in\Gamma(\jgsg^\mathcal{D})$,
    \begin{equation}
      \label{eq:discrete_non_holonomic_momentum_evolution_eq}
      (J_d)_{\ti{\xi}}(q_k,q_{k+1}) - (J_d)_{\ti{\xi}}(q_{k-1},q_{k}) = 
      -D_1L_d(q_{k-1},q_k)( \ti{\xi}(q_k)-\ti{\xi}(q_{k-1}) )_Q(q_{k-1}).
    \end{equation}
  \end{enumerate}
\end{theorem}

\begin{proof}
  \begin{itemize}
  \item \ref{it:discrete_non_holonomic_momentum-var_pple} \jiff
    \ref{it:discrete_non_holonomic_momentum-eqs_in_r}.  This is
    Proposition~\ref{prop:vertical_condition_of_motion-intrinsic}
    and Lemma~\ref{le:vertical_part_of_th_4_points_general}.
  \item \ref{it:discrete_non_holonomic_momentum-eqs_in_r} \jiff
    \ref{it:discrete_non_holonomic_momentum-momentum_eq}. For $\ti{\xi}\in
    \Gamma(\jgsg^\mathcal{D})$, using the definition of
    $(J_d)_{\ti{\xi}}$,
    equation~\eqref{eq:discrete_non_holonomic_momentum_evolution_eq}
    is equivalent to
    \begin{equation}
      \label{eq:simplified_discrete_non_holonomic_momentum_eq}
      -D_1L_d(q_k,q_{k+1})(\ti{\xi}(q_k))_Q(q_k) = 
      -D_1L_d(q_{k-1},q_k)( \ti{\xi}(q_k))_Q(q_{k-1}).
    \end{equation}
    Using now~\eqref{eq:D1_in_Q-D_2_in_RS} with $\xi=\ti{\xi}(q_k)$
    and recalling the explicit formula~\eqref{eq:bar_psi-explicit}
    for $\bar{\psi}$ we
    obtain~\eqref{eq:cond_red_vert_gen-intrinsic}.
  \end{itemize}
\end{proof}

\begin{remark}
  During the proof of
  Theorem~\ref{th:vertical_conditions_equivalence}, we saw that
  the discrete nonholonomic momentum evolution equation was
  equivalent
  to~\eqref{eq:simplified_discrete_non_holonomic_momentum_eq},
  that is simpler to use in practice and, by the same Theorem is
  equivalent to any of the other conditions in the statement.
\end{remark}

\begin{remark}
  The discrete nonholonomic evolution
  equation~\eqref{eq:discrete_non_holonomic_momentum_evolution_eq} was
  first obtained by Cort\'es and Mart\'inez
  in~\cite[Thm. 5.3]{ar:cortes_martinez-non_holonomic_integrators}, where they
  prove that any solution of the discrete Lagrange--D'Alembert
  equation
  satisfies~\eqref{eq:discrete_non_holonomic_momentum_evolution_eq}. In
  the context of groupoids, Iglesias et al. obtain the evolution
  equation of the discrete nonholonomic momentum map
  in~\cite[Thm. 3.20]{ar:iglesias_marrero_martin_martines-discrete_nonholonomic_lagrangian_systems_on_lie_groupoids}
\end{remark}


\section{Reduced equations of motion: trivial bundle case}
\label{sec:reduced_equations_of_motion_for_trivial_bundles}

In this section we consider the case where $Q:=R\times \SG$, for a
manifold $R$ and a Lie group $\SG$ that acts on $Q$ by
$l^Q_g(r,h):=(r,gh)$. In this case, $\pi:Q\rightarrow Q/\SG$ is the
trivial principal bundle $p_1:R\times \SG\rightarrow R$. The goal is
to give an explicit description of the reduced system as well as the
corresponding equations of motion.

Let $\alpha:\ti{\SG}\rightarrow R\times \SG$ be given by
$\alpha(\rho((r_0,h_0),w_0)):=(r_0,h_0^{-1}w_0h_0)$, then $\alpha$ is
a diffeomorphism with inverse
$\beta(r_0,\vartheta_0):=\rho((r_0,e),\vartheta_0)$. We have the following
commutative diagram
\begin{equation*}
  \xymatrix{{Q\times \SG \times (Q/\SG)} \ar[r]^{\check{\alpha}\times id} 
    \ar[d]_{\rho\times id} & 
    {R\times \SG\times\SG\times R} \ar[d]^{\rho^t\times id} 
    \ar@/^/[l]^{\check{\beta}\times id}\\
    {\RS} \ar[r]^{\alpha\times id} & {R\times\SG\times R} 
    \ar@/^/[l]^{\beta\times id}}
\end{equation*}
where $\rho^t(r_0,h_0,\vartheta_0):=(r_0,h_0^{-1}\vartheta_0h_0)$,
$\check{\alpha}((r_0,h_0),\vartheta_0):=(r_0,h_0,\vartheta_0)$ and
$\check{\beta}=\check{\alpha}^{-1}$. In addition, a useful ingredient
is the section $s:\ti{G}\rightarrow Q\times \SG$ given by
$s(\rho((r_0,h_0),\vartheta_0)) := ((r_0,e),h_0^{-1} \vartheta_0 h_0)$.

In the current context, by~\eqref{eq:Ad_GxG}, $\CD$ can be written as
$\CD((r_0,h_0),(r_1,h_1)) = h_1 \CDp{t}(r_0,r_1)h_0^{-1}$, where
$\CDp{t}(r_0,r_1) := \CD((r_0,e),(r_1,e))$.

Using $\alpha$ the reduced constraint manifold is
$\hat{\mathcal{D}}_d^t := (\alpha\times id)(\hat{\mathcal{D}}_d)$, so
that the reduced kinematic constraint condition becomes
$(r_0,\vartheta_0,r_1)\in \hat{\mathcal{D}}_d^t$ if and only if
$\big((r_0,e),(r_1,\vartheta_0\CDp{t}(r_0,r_1)^{-1})\big)\in \mathcal{D}_d$. If
$\mathcal{D}_d$ is described by equations $\phi_b(q_0,q_1)=0$ for all
$b$, the reduced constraint condition becomes
$\phi_b((r_0,e),(r_1,\vartheta_0\CDp{t}(r_0,r_1)^{-1})) = 0$ for all $b$.

We also have the induced lagrangians $\check{L}_d^t:=\check{L}_d\circ
(\check{\beta}\times id)$ and $\hat{L}_d^t:=
\hat{L}_d\circ(\beta\times id)$, so that
$\hat{L}_d^t(r_0,\vartheta_0,r_1) =
\check{L}_d((r_0,e),\vartheta_0,r_1)$.

Next we characterize the second order reduced manifolds $\RSsecC$ and
$\RSsec$, that are needed to set the equations of motion, according to
Corollary~\ref{cor:4_points_general-intrinsic}. It is clear that
\begin{equation*}
  \RSsecC \simeq R\times G \times G \times R \times G \times R, 
\end{equation*}
with $((r_0,h_0),w_0,r_1,\rho((r_1,h_1),w_1),r_2) \mapsto
(r_0,h_0,w_0,r_1,h_1^{-1}w_1h_1,r_2)$, and
\begin{equation*}
  \RSsec \simeq R \times G \times R \times G \times R,
\end{equation*}
with $(\rho((r_0,h_0),w_0),r_1,\rho((r_1,h_1),w_1),r_2) \mapsto
(r_0,,h_0^{-1}w_0h_0,r_1,h_1^{-1}w_1h_1,r_2)$. Next, we have to
characterize the bundles over $\RSsec$ where the morphisms
$\bar{\phi}$ and $\bar{\psi}$ are defined. Instead, it is easier to
notice that the section $s$ provides a diffeomorphism of $\RSsec$ with
its image, so that we can view $\RSsec$ as $R\times \{e\} \times G
\times R \times G \times R$. The advantage of this approach is that
instead of having to work with quotient bundles we have to work with
the restriction of bundles on $\RSsecC$ to the image of $s$.

We start by obtaining trivialized versions of the maps $\ti{F}_1$,
$F_1$ and $F_2$. They are
\begin{equation*}
  \ti{F}_1^t(r_0,h_0,\vartheta_0,r_1):=\ti{F}_1((r_0,h_0),\vartheta_0,r_1) = 
  (r_1,\vartheta_0 h_0 (\CDp{t}(r_0,r_1))^{-1}),
\end{equation*}
\begin{equation*}
  \begin{split}
    F_1^t(r_0,h_0,\vartheta_0,r_1,\vartheta_1,r_2):= &
    \; F_1((r_0,h_0),\vartheta_0,r_1,\rho((r_1,e),\vartheta_1),r_2) \\=&
     \; (r_1,\vartheta_0 h_0
    (\CDp{t}(r_0,r_1))^{-1}),
  \end{split}
\end{equation*}
and
\begin{equation*}
  \begin{split}
    F_2^t(r_0,h_0,\vartheta_0,r_1,\vartheta_1,r_2):=& \;
    F_2((r_0,h_0),\vartheta_0,r_1,\rho((r_1,e),\vartheta_1),r_2)\\=&
    \; \vartheta_0 h_0 (\CDp{t}(r_0,r_1))^{-1}) \vartheta_1
    \CDp{t}(r_0,r_1) h_0^{-1} \vartheta_0^{-1}.
  \end{split}
\end{equation*}

\begin{remark}
  The case when the symmetry group is abelian has some specially nice
  features. For one thing, even when $\pi:Q\rightarrow Q/\SG$ is not
  trivial, there is a diffeomorphism $\alpha:\ti{G}\rightarrow
  (Q/\SG)\times \SG$ that is given by
  $\alpha(\rho(q_0,w_0)):=(\pi(q_0),w_0)$. When, in addition,
  $\pi:Q\rightarrow Q/\SG$ is trivial, this diffeomorphism coincides
  with the one introduced at the beginning of this section.
\end{remark}


\subsection{Horizontal equations}
\label{sec:horizontal_equations-trivial_bundle}

Here we write down the morphism $\bar{\phi}:(\RSsec\times p_2^*TR)
\rightarrow \R$ or, rather, its realization under the identification
of $\RSsec$ with the image of the section $s$ in $\RSsecC$. Let
$\check{\phi}^t$ be the pullback of $\check{\phi}$ to $R\times \SG
\times \SG \times R \times \SG \times R$. Explicitly,
\begin{equation*}
  \begin{split}
    \check{\phi}^t(r_0,h_0,\vartheta_0,r_1,\vartheta_1,r_2,\delta r_1)
    :=& \;
    \check{\phi}((r_0,h_0),\vartheta_0,r_1,\rho((r_1,e),\vartheta_1),r_2,\delta
    r_1) \\=& \; \big(D_1 \check{L}_d((r_1,h_1),w_1,r_2)\circ
    \HLc{q_1} + D_3\check{L}_d((r_0,h_0),\vartheta_0,r_1)\\ & \; +
    (\hat{F}_d^-)^t(r_1,\vartheta_1,r_{2}) +
    (\hat{F}_d^+)^t(r_0,\vartheta_0,r_1)\big)(\delta r_1),
  \end{split}
\end{equation*}
where $h_1:=\vartheta_0 h_0 (\CDp{t}(r_0,r_1))^{-1}$, $q_1:=(r_1,h_1)$,
$w_1:=F_2^t(r_0,h_0,\vartheta_0,r_1,\vartheta_1,r_2)$, and
\begin{equation*}
  (\hat{F}_d^\pm)^t(r_1,\vartheta_1,r_{2}):= \hat{F}_d^\pm(\rho((r_1,e),\vartheta_1),r_2).
\end{equation*}
From Definition~\ref{def:reduced_discrete_force}, we see that
$\hat{F}_d^t: p_1^*TR\oplus p_3^*TR\rightarrow \R$ is a morphism of
bundles. Explicitly,
\begin{equation*}
  \begin{split}
    \hat{F}_d^t(r_0,\vartheta_0,r_1)&(\delta r_0,\delta r_1) :=\\ &
     D_2\check{L}_d((r_0,e),\vartheta_0,r_1)
    d\CD((r_0,e),(r_1,h_1))(\HLc{(r_0,e)}(\delta r_0),
    \HLc{(r_1,h_1)}(\delta r_1)),
  \end{split}
\end{equation*}
with $h_1$ as above. Even more, we can write
\begin{equation*}
  \begin{split}
    d\CD((&r_0,h_0),(r_1,h_1))(\delta r_0,\delta h_0,\delta r_1,\delta
    h_1) := \; h_1 D_1\CDp{t}(r_0,r_1)(\delta r_0) h_0^{-1} \\ & -h_1
    \CDp{t}(r_0,r_1) h_0^{-1}(\delta h_0) h_0^{-1} + h_1
    D_2\CDp{t}(r_0,r_1)(\delta r_1) h_0^{-1} + \delta h_1
    \CDp{t}(r_0,r_1) h_0^{-1},
  \end{split}
\end{equation*}
and we notice that this expression is written only in terms of
$\CDp{t}$.

We can also write the horizontal lift explicitly. In order to do so,
we view the vector bundle $\mathcal{W}\oplus \mathcal{H}\subset TQ$ as
the graph of a bundle map $M:p_1^*TR \rightarrow p_2^*T\SG$ over
$Q=R\times \SG$. Therefore,
\begin{equation*}
  \HLc{(r,h)}(\delta r) := (\delta r, M(r,h)(\delta r)) \in T_{(r,h)}Q.
\end{equation*}

Finally, we have
\begin{equation}\label{eq:bar_phi_t}
  \begin{split}
    \bar{\phi}^t(r_0,\vartheta_0,r_1,\vartheta_1,r_2,\delta r_1) :=& \;
    \check{\phi}^t(r_0,e,\vartheta_0,r_1,\vartheta_1,r_2,\delta r_1)\\
    = \; &\Bigl( D_1\hat{L}_d^t(r_1,\vartheta_1,r_2) +
    D_3\hat{L}_d^t(r_0,\vartheta_0,r_1)\\ & \;+
    D_2\hat{L}_d^t(r_1,\vartheta_1,r_2) \circ (\vartheta_1 M(r_1,e)
    -M(r_1,e) \vartheta_1) \\ & \; +
    (\hat{F}_d^-)^t(r_1,\vartheta_1,r_{2}) +
    (\hat{F}_d^+)^t(r_0,\vartheta_0,r_1)\Bigr)(\delta r_1).
  \end{split}
\end{equation}

On the other hand, we have that the representation of
$p_2^*\hat{\mathcal{D}}$ over $R\times \SG\times R\times \SG\times R$
is the vector bundle $\hat{\mathcal{D}}^t := p_3^*\hat{\mathcal{D}}$.

Finally,
condition~\eqref{eq:cond_red_hor_gen-intrinsic}
becomes
\begin{equation}
  \label{eq:cond_red_hor_gen-intrinsic-trivial_bundle}
  \bar{\phi}^t(r_0,\vartheta_0,r_1,\vartheta_1,r_2,\delta r_1)= 0
\end{equation}
for all $\delta r_1$ such that $\HLc{(r_1,e)}(\delta r_1)\in
\mathcal{D}_{(r_1,e)}$. 

Alternatively, if $\langle
\omega_1,\ldots,\omega_K\rangle = \mathcal{D}^\annihilator$,
condition~\eqref{eq:cond_red_hor_gen-intrinsic-trivial_bundle}
becomes
\begin{equation*}
  \bar{\phi}^t(r_0,\vartheta_0,r_1,\vartheta_1,r_2,\delta r_1) =
  \sum\nolimits_{a=1}^K \lambda_{a} \omega^a(r_1,e)(\delta r_1,M(r_1,e)(\delta r_1)) 
\end{equation*}
for all $\delta r_1\in T_{r_1}R$ and where $\lambda_{a}\in\R$ are
unknown.


\subsection{Vertical equations}
\label{sec:vertical_equations-trivial_bundle}

Here we write down the morphism $\bar{\psi}:(\RSsec\times\jgsg)/\SG
\rightarrow \R$ or, rather, its realization under the identification
of $\RSsec$ with the image of the section $s$ in $\RSsecC$. Let
$\check{\psi}^t$ be the pullback of $\check{\psi}$ to $R\times \SG
\times \SG \times R \times \SG \times R$. Explicitly,
\begin{equation*}
  \begin{split}
    \check{\psi}^t(r_0,h_0,\vartheta_0,r_1,&\vartheta_1,r_2,\xi_1) :=
    \check{\psi}((r_0,h_0),\vartheta_0,r_1,\rho((r_1,e),\vartheta_1),r_2,\xi_1)
    \\=& \big(D_2\check{L}_d((r_0,h_0),\vartheta_{0},r_{1}) w_{0}^{-1}
    - D_2\check{L}_d((r_1,h_1),w_1,) w_1^{-1}\big) (\xi_1),
  \end{split}
\end{equation*}
where $h_1:=\vartheta_0 h_0 (\CDp{t}(r_0,r_1))^{-1}$ and
$w_1:=F_2^t(r_0,h_0,\vartheta_0,r_1,\vartheta_1,r_2)$.

Since, $\bar{\psi}^t := \check{\psi}^t|_{R\times \{e\}\times G\times
  R\times G\times R}$,
\begin{equation}\label{eq:bar_psi_t}
  \begin{split}
    \bar{\psi}^t(r_0,\vartheta_0,&r_1,\vartheta_1,r_2,\xi_1) =
    \check{\psi}^t(r_0,e,\vartheta_0,r_1,\vartheta_1,r_2,\xi_1)\\=& \;
    \big(D_2\check{L}_d((r_0,e),\vartheta_{0},r_{1})
    \vartheta_{0}^{-1} - D_2\check{L}_d((r_1,h_1),w_1,r_2)
    w_1^{-1}\big) (\xi_1) \\=&\;
    \big(D_2\hat{L}_d^t(r_0,\vartheta_{0},r_{1}) \vartheta_{0}^{-1} -
    h_1 D_2\hat{L}_d^t(r_1,\vartheta_1,r_2) \vartheta_1^{-1}
    h_1^{-1}\big) (\xi_1),
  \end{split}
\end{equation}
where 
\begin{equation}\label{eq:trivial_bundle-vertical-h1-0}
  h_1:=\vartheta_0 (\CDp{t}(r_0,r_1))^{-1} \stext{and}
  w_1:=F_2^t(r_0,e,\vartheta_0,r_1,\vartheta_1,r_2).
\end{equation}

The other required ingredient is the pullback of
$F_1^*\jgsg^\mathcal{D}$, that is, $(F_1^t)^*\jgsg^\mathcal{D}$. Then,
\begin{equation*}
  \begin{split}
    (F_1^t)^*\jgsg^\mathcal{D} =&
    \;\big\{(r_0,h_0,\vartheta_0,r_1,\vartheta_1,r_2,\xi_1) \in
    \RSsecC\times \jgsg : ((r_1,h_1),(0,\xi_1 h_1))\in
    \mathcal{D}_{(r_1,h_1)}\big\} \\=& \;
    \big\{(r_0,h_0,\vartheta_0,r_1,\vartheta_1,r_2,\xi_1) \in
    \RSsecC\times \jgsg : ((r_1,e),(0,h_1^{-1} \xi_1 h_1))\in
    \mathcal{D}_{(r_1,e)}\big\},
  \end{split}
\end{equation*}
where $h_1$ is as in~\eqref{eq:trivial_bundle-vertical-h1-0}.
Therefore,
condition~\eqref{eq:cond_red_vert_gen-intrinsic},
becomes
\begin{equation}
  \label{eq:cond_red_vert_gen-intrinsic-trivial_bundle}
  \bar{\psi}^t(r_0,\vartheta_0,r_1,\vartheta_1,r_2,\xi_1) = 0 
\end{equation}
for all $\xi_1\in\jgsg$ such that $(0,h_1^{-1}\xi_1
h_1)\in\jgsg^\mathcal{D}_{(r_1,e)}$ where $h_1$ is defined
by~\eqref{eq:trivial_bundle-vertical-h1-0}. 

Alternatively, if $\langle \omega_1,\ldots, \omega_K
\rangle = \mathcal{D}^\annihilator$,
condition~\eqref{eq:cond_red_vert_gen-intrinsic-trivial_bundle}
becomes
\begin{equation}
  \label{eq:cond_red_vert_gen-intrinsic-trivial_bundle-annihilator}
  \bar{\psi}^t(r_0,\vartheta_0,r_1,\vartheta_1,r_2,\xi_1) =
  \sum\nolimits_{a=1}^K \lambda_{a} \omega^b(r_1,e)(0,h_1^{-1}\xi_1h_1) 
  \stext{ for all }\xi_1\in \jgsg,
\end{equation}
where $h_1$ is defined by~\eqref{eq:trivial_bundle-vertical-h1-0} and
$\lambda_{a}\in\R$ are unknown.


\subsection{Reconstruction}
\label{sec:trivial_bundle-reconstruction}

Given initial conditions on $Q\times Q$, the reconstruction process in
the current case goes as follows. Let
$\bar{q}_0=(\bar{r}_0,\bar{h}_0)$ and
$\bar{q}_1=(\bar{r}_1,\bar{h}_1)$ with
$(\bar{q}_0,\bar{q}_1)\in\mathcal{D}_d$, and a discrete curve
$(r_\cdot,\vartheta_\cdot)\in R\times \SG$. If $(r_\cdot,\vartheta_\cdot)$
satisfies the constraint $\hat{\mathcal{D}}_d^t$, the
equations~\eqref{eq:cond_red_hor_gen-intrinsic-trivial_bundle}
and~\eqref{eq:cond_red_vert_gen-intrinsic-trivial_bundle},
and $r_0=\bar{r}_0$, $r_1=\bar{r}_1$ and
$\vartheta_0=\bar{h}_0^{-1}\bar{h}_1\CDp{t}(\bar{r}_0,\bar{r}_1)$, by
Theorem~\ref{th:reconstruction-general}, the lifted curve $q_\cdot$
constructed from $\bar{q}_0$
using~\eqref{eq:reconstruction-recursive_formula} is a trajectory of
the original mechanical system on $Q\times Q$ with $q_0 = \bar{q}_0$
and $q_1 = \bar{q}_1$. In the trivial bundle case, if $q_k=(r_k,h_k)$,
\begin{equation}\label{eq:trivial_bundle-recosntruction_formula}
  q_{k+1} = \big(r_{k+1},h_k\vartheta_k\CDp{t}(r_k,r_{k+1})^{-1}\big) 
  \stext{ for all $k$.}
\end{equation}


\subsection{An example of reduction in the trivial bundle case}
\label{sec:an_example_of_reduction_in_the_trivial_bunlde_case}

In this section we specialize the previous discussion to the discrete
mechanical system $(Q,L_d,\mathcal{D},\mathcal{D}_d)$ where
$Q:=\R^2\times S^1\times S^1$ with coordinates $(x,y,\theta,\phi)$
---we consider $S^1:=\R/2\pi\Z$, and identify operations in $S^1$ with
regular addition in $\R^1$---,
\begin{equation*}
  L_d(q_0,q_1):=\frac{m}{2}\big((q_1^x-q_0^x)^2+(q_1^y-q_0^y)^2\big) + 
  \frac{I}{2}(q_1^\theta - q_0^\theta)^2 +\frac{J}{2}(q_1^\phi-q_0^\phi)^2,
\end{equation*}
\begin{equation*}
  \begin{split}
    \mathcal{D}_q:=& \langle \del_\phi, \del_\theta+ A \cos(q^\phi)
    \del_x + A\sin(q^\phi)\del_y \rangle = \langle
    \omega_1(q),\omega_2(q) \rangle^\annihilator \subset T_qQ,
  \end{split}
\end{equation*}
where $\omega_1(q):= dx -r\cos(q^\phi)d\theta$ and
$\omega_2(q):=dy-r\sin(q^\phi) d\theta$, and
\begin{equation*}
  \begin{split}
    \mathcal{D}_d:= \{(q_0,q_1)\in Q\times Q : &q_1^x-q_0^x = A
    (q_1^\theta-q_0^\theta)(\cos((q_0^\phi+q_1^\phi)/2) \text{ and } \\
    &q_1^y-q_0^y = A
    (q_1^\theta-q_0^\theta)(\sin((q_0^\phi+q_1^\phi)/2)\},
  \end{split}
\end{equation*}
for $m$, $A$, $I$ and $J$ positive constants. This system can be
obtained as a discretization of the classical mechanical system formed
by a vertical disk of radius $A$, mass $m$ with inertia momenta $I$ and
$J$, rolling without slipping on a horizontal plane. 

We consider the Lie group $\SG:=\R^2\times S^1$ acting on $Q$ by
\begin{equation*}
  l^Q_g(q) := (q^x+g^x,q^y+g^y,q^\theta+g^\theta,q^\phi).
\end{equation*}
The corresponding lifted action is $l^{TQ}_g(q,v) = (l^Q_g(q),v)$.  We
are in the trivial bundle case analyzed above because
$\pi:Q\rightarrow Q/\SG$ is the trivial principal bundle
$p_2:\SG\times S^1\rightarrow S^1$ with structure group $\SG$, which
is a symmetry group of the system.

We consider various subbundles of $TQ$: the vertical bundle
$\mathcal{V}^\SG_q=\langle \del_x, \del_y, \del_\theta \rangle$, the
intersection bundle $\mathcal{S}_q= \langle
\del_\theta+A\cos(q^\phi)\del_x + A\sin(q^\phi)\del_y\rangle$ and the
corresponding complements $\mathcal{H}_q:=\langle \del_\phi\rangle$,
$\mathcal{U}_q:=\langle \del_x, \del_y\rangle$ and
$\mathcal{W}:=\{0\}$.  Following
Definition~\ref{def:non_holonomic_connection}, the previous bundles
induce a nonholonomic connection $\CC$ whose horizontal space is
$Hor_\CC:=\mathcal{H}$. The corresponding horizontal lift
$\HLc{q}:T_{r}S^1\rightarrow \mathcal{H}_q \subset T_qQ$ for $q\in
\pi^{-1}(r)$ is $\HLc{q}(\dot{r}\del_r) = \dot{r}\del_\phi$; notice
that, in this case, the bundle map $M$ that describes the horizontal
bundle as a graph is identically $0$.

By Proposition~\ref{prop:discrete_connection_vs_one_form}, since
\begin{equation*}
  \CD(q_0,q_1) = (q_1^x,q_1^y,q_1^\theta) (q_0^x,q_0^y,q_0^\theta)^{-1} = 
  (q_1^x-q_0^x, q_1^y-q_0^y,q_1^\theta-q_0^\theta) \in \SG
\end{equation*}
satisfies~\eqref{eq:Ad_GxG}, there is a discrete connection with level
$\gamma(q)=e$ whose discrete connection $1$-form is $\CD$. Notice that
our choice of $\CD$ has $\CDp{t}(r_0,r_1)=e$ for all $r_0,r_1\in
S^1$. The corresponding discrete horizontal lift is $\HLd{q_0}(r_1) =
(q_0,(q_0^x,q_0^y,q_0^\theta,r_1))$.

The reduced space in this case is $\RS\simeq S^1\times \SG\times S^1$
and the reduced second order manifold is $\RSsec\simeq S^1\times
\SG\times S^1\times \SG\times S^1$. The trivialized reduced lagrangian is 
\begin{equation*}
  \hat{L}_d^t(r_0,\vartheta_0,r_1) =
  \frac{m}{2}\big((\vartheta^x_0)^2+
  (\vartheta^y_0)^2\big) 
  + \frac{I}{2}(\vartheta^\theta_0)^2 + \frac{J}{2}(r_1-r_0)^2.
\end{equation*}

Next we describe the horizontal and vertical equations of the reduced
system, but first we notice that since $\CD$ only depends on the group
variables, it turns out that $d\CD$ annihilates horizontal vectors and
the mixed curvature $\BM$ vanishes. Consequently, there are no reduced
discrete forces on the system.

From~\eqref{eq:bar_phi_t}, we derive
\begin{equation*}
  \bar{\phi}^t(r_{k-1},\vartheta_{k-1},r_k,\vartheta_k,r_{k+1},\delta r_k) = 
  \big(J \big((r_k-r_{k-1})- (r_{k+1}-r_k)\big) dr|_{r_k}\big)(\delta r_k),
\end{equation*}
and, noticing that $\mathcal{W}=\{0\}$ implies that
$\hat{\mathcal{D}}^t = TS^1$, the horizontal equation is
\begin{equation}\label{eq:trivial_bundle-example-horizontal_equation}
  (r_k-r_{k-1}) - (r_{k+1}-r_k)=0.
\end{equation}

From~\eqref{eq:bar_psi_t}, we obtain
\begin{equation*}
  \begin{split}
    \bar{\psi}^t(r_{k-1},\vartheta_{k-1},r_k,\vartheta_k,r_{k+1},\xi_k) =& \;\big(
    m(\vartheta^x_{k-1}-\vartheta^x_k) dx|_{e^x} + m ( \vartheta^y_{k-1}-\vartheta^y_k)
    dy|_{e^y} \\ & \;+ I (\vartheta^\theta_{k-1}-\vartheta^\theta_k)
    d\theta|_{e^\theta} \big)(\xi_k),
  \end{split}
\end{equation*}
where $e=(e^x,e^y,e^\theta)$ is the identity element in $\SG$. Since
the right hand side
of~\eqref{eq:cond_red_vert_gen-intrinsic-trivial_bundle-annihilator}
is $\lambda_1 (dx|_{e^x}-A\cos(r_k) d\theta|_{e^\theta}) + \lambda_2
(dy|_{e^y} - A\sin(r_k) d\theta|_{e^\theta})$, we obtain the following
vertical equations
\begin{equation*}
  \begin{split}
    m (\vartheta^x_{k-1}-\vartheta^x_k) =& \; \lambda_1,\\
    m (\vartheta^y_{k-1}-\vartheta^y_k) =& \; \lambda_2,\\
    I (\vartheta^\theta_{k-1}-\vartheta^\theta_k) =& -\lambda_1 A\cos(r_k)
    -\lambda_2 A\sin(r_k).
  \end{split}
\end{equation*}
On the other hand, the reduced kinematic constraint equations are
\begin{equation}\label{eq:trivial_bundle-example-kinematic_constraint_equation}
  \begin{split}
    \vartheta_k^x = & \; A\vartheta_k^\theta \cos((r_k+r_{k+1})/2),\\
    \vartheta_k^y = & \; A\vartheta_k^\theta \sin((r_k+r_{k+1})/2).
  \end{split}
\end{equation}

Next, we find the reduced dynamics by solving the reduced
equations. From~\eqref{eq:trivial_bundle-example-horizontal_equation}
\begin{equation}\label{eq:trivial_bundle-example-r_solution}
  r_k = (r_1-r_0) k + r_0 \stext{ for all } k\in\N\cup\{0\}.
\end{equation}
From the vertical equations we obtain
\begin{equation*}
  I (\vartheta^\theta_{k-1}-\vartheta^\theta_k) = - m (\vartheta^x_{k-1}-\vartheta^x_k)A\cos(r_k) 
  - m (\vartheta^y_{k-1}-\vartheta^y_k) A\sin(r_k).
\end{equation*}
Plugging~\eqref{eq:trivial_bundle-example-kinematic_constraint_equation}
into this last equation (for $k-1$ and $k$) and simplifying we get
\begin{equation*}
  \vartheta^\theta_k = \frac{I + m A^2 \cos((r_k-r_{k-1})/2)}{I+m A^2 \cos((r_{k+1}-r_k)/2)} \vartheta^\theta_{k-1}.
\end{equation*}
Using the horizontal
equation~\eqref{eq:trivial_bundle-example-horizontal_equation} we
conclude that $\vartheta^\theta_k = \vartheta^\theta_0$ for all $k$. Thus,
\begin{equation*}
  \vartheta^x_k = A \vartheta^\theta_0 \cos\bigg((r_1-r_0)k+\frac{r_1+r_0}{2}\bigg) \stext{ and } 
  \vartheta^y_k = A \vartheta^\theta_0 \sin\bigg((r_1-r_0)k+\frac{r_1+r_0}{2}\bigg).
\end{equation*}
This last expression completes the description of the reduced system's
dynamics.

The reconstruction of the trajectories of the original system in
$Q\times Q$ is done as described in
Section~\ref{sec:trivial_bundle-reconstruction}. By~\eqref{eq:trivial_bundle-recosntruction_formula},
the reconstructed trajectory $q_\cdot$ is
\begin{equation}\label{eq:trivial_bundle-example-reconstruction_iterative_q}
  q_{k+1} = \big((q^x_{k+1},q^y_{k+1},q^\theta_{k+1}),q^\phi_{k+1}\big) = 
  \big((q^x_{k},q^y_{k},q^\theta_{k})\vartheta_k,r_{k+1}\big),
\end{equation}
so that $(q^x_{k+1},q^y_{k+1},q^\theta_{k+1}) =
(q^x_{0},q^y_{0},q^\theta_{0}) \prod_{j=0}^{k} \vartheta_j$ or,
\begin{align*}
  q^x_{k+1} = & \; q^x_0 + A (q^\theta_1 -q^\theta_0) \cos((q^\phi_1-q^\phi_0)/2) \frac{\sin((q^\phi_1-q^\phi_0)k+q^\phi_0) - \sin(q^\phi_0)}{\sin(q^\phi_1-q^\phi_0)},\\
  q^y_{k+1} = & \; q^y_0 + A (q^\theta_1 -q^\theta_0) \cos((q^\phi_1-q^\phi_0)/2) \frac{\cos(q^\phi_0)-\cos((q^\phi_1-q^\phi_0)k+r_0)} {\sin(q^\phi_1-q^\phi_0)},\\
  q^\theta_{k+1} = & \; (q^\theta_1- q^\theta_0) k +q^\theta_0.\\
  \intertext{Also,
    from~\eqref{eq:trivial_bundle-example-reconstruction_iterative_q}
    and~\eqref{eq:trivial_bundle-example-r_solution},} q^\phi_{k+1} =
  & \;(q_1^\phi-q_0^\phi)(k+1)+ q_0^\phi.
\end{align*}


\section{Reduced equations of motion: Chaplygin case}
\label{sec:reduced_equations_of_motion-chaplygin}

In this section we specialize Theorems~\ref{th:4_points-general}
and~\ref{th:reconstruction-general} to the case where the original
system is a discrete mechanical system with Chaplygin type symmetries.
In this case we go beyond the result of equivalence between the
discrete mechanical system in $Q$ and a dynamical system in $\RS$ to
obtain an equivalence between discrete mechanical systems on $Q$
and $Q/\SG$.

\begin{definition}\label{def:chaplygin_symmetry}
  A symmetry group $\SG$ of $(Q,L_d,\mathcal{D},\mathcal{D}_d)$ is a
  \jdef{Chaplygin type symmetry group} if it satisfies the following
  conditions.
  \begin{enumerate}
  \item $TQ=\mathcal{V}^{\SG}\oplus \mathcal{D}$ and
  \item \label{it:chaplygin-discrete_condition} $\mathcal{D}_{d}$
    defines an affine discrete connection $\CD$ on the principal
    bundle $\pi:Q\rightarrow Q/\SG$.
  \end{enumerate}
\end{definition}

Notice that, by definition, the condition $TQ=\mathcal{V}^{\SG}\oplus
\mathcal{D}$ is equivalent to the fact that the
decomposition~\eqref{eq:TQ_decomposition} has the form
\begin{equation}\label{eq:descomposicion_TQ-chaplygin}
  TQ = \underbrace{\{0\}}_{\mathcal{W}} \oplus 
  \underbrace{\mathcal{V}^\SG}_{\mathcal{U}}\oplus 
  \underbrace{\{0\}}_{\mathcal{S}} \oplus \underbrace{\mathcal{D}}_{\mathcal{H}}.
\end{equation}
As in the general case, this decomposition defines a connection $\CC$
on $Q\rightarrow Q/\SG$, whose horizontal space is
$\mathcal{H}=\mathcal{D}$. Condition~\ref{it:chaplygin-discrete_condition}
in Definition~\ref{def:chaplygin_symmetry} requires that, for every
$(q_0,q_1)$ there exists a unique $g\in\SG$ such that
$(q_0,l^Q_{g^{-1}}(q_1))\in\mathcal{D}_d$.

\begin{example}
  $\SG:=\R$ is a Chaplygin type symmetry group of the discrete
  mechanical system $(Q,L_d,\mathcal{D},\mathcal{D}_d)$ defined in
  Example~\ref{ex:continued_example-discrete_mechanical_system}. The
  first condition is clear from the decomposition of $TQ$ that appears
  in Example~\ref{ex:continued_example-mixed_curvature}; the second
  point corresponds to the discrete connection $\CDp{b}$ 
  considered in
  Example~\ref{ex:continued_example-affine_discrete_connection} when
  $b=1$.
\end{example}


\subsection{An inclusion}
\label{sec:an_inclusion}

Let $\mathcal{Y}:(Q/\SG)\times (Q/\SG) \rightarrow\RS$ be defined by
\begin{equation*}
  \mathcal{Y}(r_0,r_1) := (\rho(q_0,e),r_1),
\end{equation*}
where $q_0\in\pi^{-1}(r_0)$ and $e$ is the identity of $\SG$.

\begin{lemma}
  The application $\mathcal{Y}$ is well defined. Even more, if $s$ is
  a local section of $\pi: Q\rightarrow Q/\SG$, then
  $\mathcal{Y}(r_0,r_1)=(\rho(s(r_0),e),r_1)$.
\end{lemma}

\begin{proof}
  $\mathcal{Y}$ does not depend on the choice of $q_0\in\pi^{-1}(r_0)$
  by the $\SG$-invariance of $\rho$. The expression of $\mathcal{Y}$
  in terms of $s$ holds because $s(r_0)\in \pi^{-1}(r_0)$. This
  formula also shows that $\mathcal{Y}$ is a smooth map.
\end{proof}
Using $\mathcal{Y}$ we transport the existing structure on $\RS$ to
$(Q/\SG)\times (Q/\SG)$. More precisely, we define the (forced,
unconstrained) discrete mechanical system
$(Q/\SG,\breve{L}_d,\breve{F}_d)$ where
$\breve{L}_d:=\mathcal{Y}^*(\hat{L}_d)=\hat{L}_d\circ \mathcal{Y}$ and
the discrete force $\breve{F}_d := \mathcal{Y}^*(\hat{F}_d)$. We also
define $\breve{\mathcal{D}}:=\hat{\mathcal{D}}$ as a subbundle of
$T(Q/\SG)$, but notice that, $\breve{\mathcal{D}}=T(Q/\SG)$ in the
Chaplygin case, due to~\eqref{eq:descomposicion_TQ-chaplygin}.

\begin{lemma}\label{le:relation_breveLd_checkL_d}
  If $\CC$ is a connection on the principal bundle $\pi:Q\rightarrow
  Q/\SG$ and $q_0\in \pi^{-1}(r_0)$, then
  \begin{equation}\label{eq:relation_breveLd_checkL_d}
    d\breve{L}_d(r_0,r_1)(\delta r_0,\delta r_1) =
    D_1\check{L}_d(q_0,e,r_1)(\HLc{q_0}(\delta r_0)) +
    D_3\check{L}_d(q_0,e,r_1)(\delta r_1).
  \end{equation}
\end{lemma}

\begin{proof}
  If $s$ is a local section of $Q\rightarrow Q/\SG$ with $q_0=s(r_0)$
  then $\breve{L}_d(r_0,r_1) = \hat{L}_d(\rho(s(r_0),e),r_1) =
  \check{L}_d(s(r_0),e,r_1)$, so that
  \begin{equation}\label{eq:relation_breveLd_checkL_d-aux0}
    \begin{split}
      d\breve{L}_d(r_0,r_1)(\delta r_0,&\delta r_1) =
      (d\check{L}_d)(s(r_0),e,r_1)(ds(r_0)(\delta r_0),0,\delta r_1)
      \\=& D_1\check{L}_d(s(r_0),e,r_1)(ds(r_0)(\delta r_0)) +
      D_3\check{L}_d(s(r_0),e,r_1)(\delta r_1).
    \end{split}
  \end{equation}

  Since $\pi\circ s = id_{Q/\SG}$, we have
  \begin{equation*}
    \begin{split}
      \delta r_0 = d\pi(s(r_0))ds(r_0)(\delta r_0) = d\pi(s(r_0))
      Hor_{\CC}(ds(r_0)(\delta r_0)),
    \end{split}
  \end{equation*}
  thus $Hor_{\CC}(ds(r_0)(\delta r_0))=\HLc{s(r_0)}(\delta r_0)$.
  Also, as $D_1\check{L}_d(s(r_0),e,r_1)$ vanishes on vertical vectors
  by the $\SG$-invariance of $\check{L}_d$,
  \begin{equation*}
    \begin{split}
      D_1\check{L}_d(s(r_0),e,r_1)(ds(r_0)(\delta r_0)) =&
      D_1\check{L}_d(s(r_0),e,r_1) (Hor_{\CC}(ds(r_0)(\delta r_0))) \\=&
      D_1\check{L}_d(s(r_0),e,r_1)(\HLc{s(r_0)}(\delta r_0))
    \end{split}
  \end{equation*}
  which, replaced in~\eqref{eq:relation_breveLd_checkL_d-aux0} leads
  to~\eqref{eq:relation_breveLd_checkL_d}.
\end{proof}

\begin{lemma}\label{le:relation_hatSd_hatS_d}
  Let $(v_\cdot,r_\cdot)$ be a discrete curve in $\RS$ such that
  $v_k=\rho(q_k,e)$ for some $q_k\in Q$ and all $k$. Then
  \begin{equation*}
    d\hat{S}_d(v_\cdot,r_\cdot)(\delta v_\cdot,\delta r_\cdot) = 
    d\breve{S}_d(r_\cdot)(\delta r_\cdot) + 
    \sum_{k=1}^{N-1} \big( \breve{F}_d^-(r_k,r_{k+1}) + 
    \breve{F}_d^+(r_{k-1},r_k) \big)(\delta r_k)
  \end{equation*}
  for all vanishing end points variations $\delta r_\cdot$ and
  \begin{equation*}
    \delta v_k = d\rho(q_k,e)\big( \HLc{q_k}(\delta r_k),
    d\CD(q_k,q_{k+1})(\HLc{q_k}(\delta r_k), \HLc{q_{k+1}}(\delta r_{k+1}))\big)
  \end{equation*}
  for all $k$.
\end{lemma}

\begin{proof}
  By definition, using the explicit form of $\delta v$,
  \begin{equation*}
    \begin{split}
      d\hat{S}_d(v_\cdot,r_\cdot)(\delta v_\cdot,\delta r_\cdot) =&
      \sum\nolimits_{k=0}^{N-1} \big(
      D_1\check{L}_d(q_k,e,r_{k+1})(\HLc{q_k}(\delta r_k)) \\&+
      D_2\check{L}_d(q_k,e,r_{k+1})(d\CD(q_k,q_{k+1})(\HLc{q_k}(\delta
      r_k), \HLc{q_{k+1}}(\delta r_{k+1}))) \\ &+
      D_3\check{L}_d(q_k,e,r_{k+1})(\delta r_{k+1})\big).
    \end{split}
  \end{equation*}
  Using Lemma~\ref{le:relation_breveLd_checkL_d} and the decomposition
  $d\CD = D_1\CD + D_2\CD$ we obtain
  \begin{equation*}
    \begin{split}
      d\hat{S}_d(v_\cdot,r_\cdot)(\delta v_\cdot,\delta r_\cdot) = &
      \; d\breve{S}(r_\cdot)(\delta r_\cdot) \\&+
      \sum\nolimits_{k=0}^{N-1} \big(\hat{F}_d^-(q_k,e,r_{k+1})(\delta
      r_k) + \hat{F}_d^+(q_k,e,r_{k+1})(\delta r_{k+1})\big),
    \end{split}
  \end{equation*}
  and the result follows.
\end{proof}

\begin{lemma}\label{le:relation_hor_eqs_hat_breve}
  Let $(v_\cdot,r_\cdot)$ be a discrete curve in $\RS$ such that
  $v_k=\rho(q_k,e)$ for some $q_k\in Q$ and all $k$. Then
  $(v_\cdot,r_\cdot)$ satisfies
  condition~\eqref{eq:cond_red_hor_gen-intrinsic}
  if and only if
  \begin{equation}\label{eq:relation_hor_eqs_hat_breve}
    D_1\breve{L}_d(r_k,r_{k+1}) + D_2\breve{L}_d(r_{k-1},r_k)
    + \breve{F}_d^+(r_{k-1},r_k)+\breve{F}_d^-(r_k,r_{k+1}) 
    \in (\breve{\mathcal{D}}_{r_k})^\annihilator
  \end{equation}
  for all $k$.
\end{lemma}

\begin{proof}
  Considering the form of $v_k$, the vanishing
  condition~\eqref{eq:cond_red_hor_gen-intrinsic}
  immediately translates, via $\mathcal{Y}$ to the vanishing
  condition~\eqref{eq:relation_hor_eqs_hat_breve}.
\end{proof}


\subsection{Reduced dynamics}
\label{sec:chaplygin_reduction}

\begin{theorem}\label{th:4_points_chaplygin}
  Let $\SG$ be a Chaplygin symmetry group of
  $(Q,L_d,\mathcal{D},\mathcal{D}_d)$, $q_\cdot$ be a discrete curve
  in $Q$, and $r_{k}:=\pi(q_{k})$ be the corresponding curve in
  $Q/\SG$. If $(q_k,q_{k+1})\in \mathcal{D}_d$ for all $k$, the
  following statements are equivalent.
  \begin{enumerate}
  \item \label{it:var_pple-chaplygin} $q_\cdot$ satisfies the
    variational principle $dS_d(q)(\delta q_\cdot) = 0$ for all
    vanishing end points variations $\delta q_\cdot$ such that $\delta
    q_k \in \mathcal{D}_{q_k}$ for all $k$.
  \item \label{it:eq_lda-chaplygin} $q_\cdot$ satisfies the
    Lagrange--D'Alembert equations~\eqref{eq:dla-eqs}.
  \item \label{it:red_var_pple-chaplygin} $r_\cdot$ satisfies the
    variational principle
    \begin{equation*}
      d\breve{S}_d(r_\cdot) (\delta r_\cdot) =
      -\sum_{k=1}^{N-1} (\breve{F}_d^-(r_k,r_{k+1}) + 
      \breve{F}_d^+(r_{k-1},r_k))(\delta r_k)
    \end{equation*}
    for all vanishing end points variations $\delta r_\cdot$ with
    $\delta r_k \in T_{r_k}(Q/\SG)$ and where
    \begin{equation*}
      \breve{S}_d(r_\cdot):=\sum_{k=0}^{N-1}\breve{L}_d(r_k,r_{k+1}).      
    \end{equation*}
  \item \label{it:red_eq_lda-chaplygin} $r_\cdot$ satisfies the
    following equation, for all $k$.
    \begin{equation*}
      D_1\breve{L}_d(r_k,r_{k+1}) + D_2\breve{L}_d(r_{k-1},r_k) = 
      -(\breve{F}_d^+(r_{k-1},r_k)+\breve{F}_d^-(r_k,r_{k+1})).
    \end{equation*}
  \end{enumerate}
\end{theorem}

\begin{proof}
  Since $\mathcal{S}=\{0\}$, by Theorem~\ref{th:4_points-general}, we
  see that each one of the points~\ref{it:var_pple-chaplygin}
  and~\ref{it:eq_lda-chaplygin} of the present result is equivalent
  to any one of
  \begin{roman-enumerate}
  \item \label{it:chaplygin_red_intermediate-3}
    Item~\ref{it:red_var_pple-general} in
    Theorem~\ref{th:4_points-general} holds for all horizontal
    variations, ($\xi_k=0$ for all $k$).
  \item \label{it:chaplygin_red_intermediate-4} For all $k$
    $(v_k,r_{k+1})\in\hat{\mathcal{D}}_d$ and
    condition~\eqref{eq:cond_red_hor_gen-Q}
    holds.
  \end{roman-enumerate}

  As $Hor_\CD = \mathcal{D}_d$, $(q_k,q_{k+1})\in\mathcal{D}_d$ if and
  only if $w_k=\CD(q_k,q_{k+1})=e$. Then,
  \begin{equation*}
    (v_k,r_{k+1}) \in\hat{\mathcal{D}}_d \jiff (\rho(q_k,e),r_{k+1})
    \in \Phi_\CD(\mathcal{D}_d/\SG) \jiff 
    (q_k,\HLds{q_k}(r_{k+1})) \in \mathcal{D}_d,
  \end{equation*}
  which always holds by definition of $\CD$. Therefore
  $\hat{\mathcal{D}}_d=\RS$, so that we can drop the reduced kinematic
  constraint condition from~\ref{it:chaplygin_red_intermediate-3}
  and~\ref{it:chaplygin_red_intermediate-4}.

  By Lemma~\ref{le:relation_hatSd_hatS_d} the variational principle
  that appears in~\ref{it:chaplygin_red_intermediate-3} is equivalent
  to the one in point~\ref{it:red_var_pple-chaplygin} of the
  statement. Similarly, as $\breve{\mathcal{D}}=T(Q/\SG)$,
  Lemma~\ref{le:relation_hor_eqs_hat_breve} and
  Proposition~\ref{prop:horizontal_condition_of_motion-intrinsic} show
  that~\ref{it:chaplygin_red_intermediate-4} is equivalent to the one
  in point~\ref{it:red_eq_lda-chaplygin} of the statement.
\end{proof}

\begin{remark}
  Notice that the reduction of a discrete mechanical system with
  Chaplygin type symmetry results in an unconstrained discrete
  mechanical system but with external forces given by $\breve{F}^{\pm
  }$ in the previous theorem. A similar analysis has been done in the
  more general groupoid setting
  in~\cite{ar:iglesias_marrero_martin_martines-discrete_nonholonomic_lagrangian_systems_on_lie_groupoids}.
\end{remark}

\begin{example}\label{ex:continued_example-chaplygin_reduction}
  As we noted above, $\SG=\R$ is a Chaplygin type symmetry group of
  the discrete mechanical system introduced in
  Example~\ref{ex:continued_example-discrete_mechanical_system}. We
  notice that the construction described in this section corresponds
  to performing the reduction using the discrete connection $\CDp{1}$
  from Example~\ref{ex:continued_example-affine_discrete_connection}.
  Furthermore, the inclusion $\mathcal{Y}$ from
  section~\ref{sec:an_inclusion} was essentially already present in
  the analysis of the reduction in
  Example~\ref{ex:continued_example-reduced_equations} in the form of
  the section $(Q/\SG)\times (Q/\SG) \rightarrow \ti{\SG}$ given by
  $s(r_0,r_1):=\rho((r_0,0),r_1)$.

  In any case, the discrete unconstrained mechanical system associated
  to the reduced system by Theorem~\ref{th:4_points_chaplygin} is
  \begin{gather*}
    \breve{Q}:=Q/\SG = \R,\stext{ with coordinate } r,\\
    \breve{L}_d(r_0,r_1) = \frac{m}{2} \big( (r_1-r_0)^2 +
    (r_1^2-r_0^2)^2/4\big)
  \end{gather*}
  with no forces since, by
  Example~\ref{ex:continued_example-mixed_curvature}, the mixed
  curvature vanishes in the case $b=1$, hence there are no reduced
  forces. Finally, the discrete Euler--Lagrange equation that
  determines the evolution of the $(\breve{Q},\breve{L}_d)$ system is
  \begin{equation}\label{eq:continued_example-reduced_equations_chaplygin-r}
    (r_{k+1}-r_k)-(r_k-r_{k-1}) + 
    r_k \big((r_{k+1}^2-r_k^2)-(r_k^2-r_{k-1}^2)\big)/2 =0,
  \end{equation}
  that is, precisely,~\eqref{eq:continued_example-reduced_equations-r}.
\end{example}

Last, we adapt Theorem~\ref{th:reconstruction-general} to the
reconstruction in the present setting.
\begin{theorem}\label{th:reconstruction-chaplygin}
  Let $\SG$ be a Chaplygin symmetry group of
  $(Q,L_d,\mathcal{D},\mathcal{D}_d)$. Let $(\bar{q}_0,\bar{q}_1)\in
  \mathcal{D}_d$ and $r_\cdot$ be a discrete trajectory of the forced
  discrete mechanical system
  $(Q/\SG,\breve{L}_d,\breve{\mathcal{D}},\breve{\mathcal{D}}_d,\breve{F}_d)$
  such that $\pi(\bar{q}_j)=r_j$ for $j=0,1$. Define the discrete
  curve $q_\cdot$ in $Q$ inductively by $q_0=\bar{q}_0$ and $q_{k+1} =
  \HLds{q_k}(r_{k+1})$ for all $k\in\NZ$ (here the horizontal lift is
  associated to the affine discrete connection whose horizontal space
  is $\mathcal{D}_d$). Then $q_\cdot$ is a trajectory of the original
  discrete mechanical system with $q_j=\bar{q}_j$ for $j=0,1$.
\end{theorem}

\begin{proof}
  The trajectory $r_\cdot$ defines a trajectory $(v_k,r_{k+1}) =
  \mathcal{Y}(r_k,r_{k+1})$ in $\RS$. Indeed, by
  Lemma~\ref{le:relation_hatSd_hatS_d} both curves satisfy
  simultaneously the corresponding variational principles. Also, since
  $\bar{q}_0\in \pi^{-1}(r_0)$, $v_0 = \rho(\bar{q}_0,e) =
  \rho(\bar{q}_0,\CD(\bar{q}_0,\bar{q}_1))$. In general,
  $v_k=\rho(\ti{q}_k,e)$ for some $\ti{q}_k\in \pi^{-1}(r_k)$.

  Then, since $(\bar{q}_0,\bar{q}_1)\in\mathcal{D}_d$, by
  Theorem~\ref{th:reconstruction-general}, there is a trajectory
  $q_\cdot$ in $Q$ such that $\pi(q_k)=r_k$ for all $k$ and $q_0 =
  \bar{q}_0$, $q_1=\bar{q}_1$. Furthermore, $q_{k+1}$ satisfies the
  formula in the statement because, $q_{k+1}$
  satisfies~\eqref{eq:reconstruction-recursive_formula} with $u_k=e$
  for all $k$.
\end{proof}

\begin{example}\label{ex:continued_example-chaplygin_reconstruction}
  We use Theorem~\ref{th:reconstruction-chaplygin} to reconstruct the
  trajectories of the system reduced in
  Example~\ref{ex:continued_example-chaplygin_reduction}. Applying the
  recursive formula
  and~\eqref{eq:continued_example-discrete_connection_formulas} we
  have
  \begin{equation*}
    q_{k+1} = (x_{k+1},y_{k+1}) = \HLds{q_k}(r_{k+1}) = 
    (r_{k+1},y_k+(r_{k+1}^2-r_k^2)/2),
  \end{equation*}
  or, simplifying, $(x_k,y_k) = (r_k,y_0+(r_k^2-r_0^2)/2)$, agreeing
  with the result obtained in
  Example~\ref{ex:continued_example-reconstruction}.
\end{example}


\section{Reduced equations of motion: horizontal symmetries case}
\label{sec:reduced_equations_of_motion-horizontal}

In this section we specialize theorem \ref{th:4_points-general} to the
case where the original system is a mechanical system with horizontal
symmetries. In the same way as in the Chaplygin case discussed in
Section~\ref{sec:reduced_equations_of_motion-chaplygin}, in this case
we go beyond the result of equivalence between the mechanical
system in $Q$ and a dynamical system in $\RS$ to obtain an equivalence
between the discrete mechanical system in $Q$ and another one in
$Q/\SG$.

\begin{definition}
  Let $M$ be a symmetry group of
  $(Q,L_d,\mathcal{D},\mathcal{D}_d)$. A closed subgroup $\SG\subset
  M$ is said to be a \jdef{horizontal symmetry subgroup} for
  $(Q,L_d,\mathcal{D},\mathcal{D}_d)$ if
  \begin{equation} \label{eq:horizontal_symmetry-continuous_condition}
    \mathcal{V}^M(q)\cap \mathcal{D}_q= \mathcal{V}^\SG(q) \stext{for
      all} q\in Q.
  \end{equation}
\end{definition}

From now on we will forget the group $M$ and consider the action of
$\SG$ on the system. It is in this context that we specialize
Theorem~\ref{th:4_points-general}. Due to
condition~\eqref{eq:horizontal_symmetry-continuous_condition}, we have
that $\mathcal{V}^\SG = \mathcal{S}$, the
decomposition~\eqref{eq:TQ_decomposition} of $TQ$ becomes
\begin{equation*}
  TQ = \mathcal{W}\oplus \underbrace{\{0\}}_{\mathcal{U}} \oplus \mathcal{S} 
  \oplus \mathcal{H},
\end{equation*}
for any complementary subbundles $\mathcal{H}$ of $\mathcal{S}$ in
$\mathcal{D}$ and $\mathcal{W}$ of $\mathcal{D}$ in
$TQ$. Fixing one such decomposition we define a connection
$\CC$ on the principal bundle $\pi:Q\rightarrow Q/\SG$ requiring that
$Hor_\CC = \mathcal{H}$.

In the context of this section the discrete nonholonomic momentum map
$J_{d}$ defined by~\eqref{eq:non_holonomic_discrete_momentum-def} has
some special properties, which are studied next.

\begin{lemma}
  If $\SG$ is horizontal symmetry subgroup the following statements are
  true.
  \begin{enumerate}
  \item \label{it:hor_sym-triv_gD} $\jgsg^\mathcal{D} = Q\times \jgsg$
    and $(\jgsg^\mathcal{D})^* \simeq Q\times \jgsg^*$.
  \item \label{it:hor_sym-Jd} Composing the nonholonomic momentum
    map $J_{d}$ defined in~\eqref{eq:non_holonomic_discrete_momentum-def}
    with the projection onto the second variable defines a momentum
    application $J_d:Q\times Q\rightarrow \jgsg^*$. Explicitly, for
    $\xi \in \jgsg$,
    \begin{equation*}
      J_d(q_0,q_1)\xi := -D_1 L_d(q_0,q_1)\xi_Q(q_0).
    \end{equation*}
  \item \label{it:th_iff_cons_Jd} Any of the equivalent conditions of
    Theorem~\ref{th:vertical_conditions_equivalence} is equivalent
    to the condition that $J_{d}$ is constant on the trajectory
    $q_\cdot$.
  \end{enumerate}
\end{lemma}

\begin{proof}
  Item~\ref{it:hor_sym-triv_gD} is a direct consequence of
  $\mathcal{V}^\SG\subset \mathcal{D}$, while item~\ref{it:hor_sym-Jd}
  is clear from item~\ref{it:hor_sym-triv_gD}.

  Last we check item~\ref{it:th_iff_cons_Jd}.  Let $q_\cdot$ be a
  discrete curve in $Q$. Assume that
  equation~\eqref{eq:discrete_non_holonomic_momentum_evolution_eq}
  holds on $q_\cdot$ for any section $\ti{\xi}$. Then, since any
  $\xi\in\jgsg$ defines a (constant) section,
  evaluating~\eqref{eq:discrete_non_holonomic_momentum_evolution_eq}
  on this section yields $\big( J_d(q_k,q_{k+1})
  -J_d(q_{k-1},q_k)\big) \xi = 0$. Thus, since $\xi\in\jgsg$ is
  arbitrary, $J_d$ is conserved on $q_\cdot$.

  Conversely, if $J_d$ is constant on $q_\cdot$,
  equation~\eqref{eq:discrete_non_holonomic_momentum_evolution_eq}
  holds for constant sections of $\jgsg^\mathcal{D}$. But, it can be
  readily checked that
  if~\eqref{eq:discrete_non_holonomic_momentum_evolution_eq} holds
  for a section $\ti{\xi}$, it also holds for the section $f\ti{\xi}$,
  for arbitrary $f:Q\rightarrow\R$. Then, since in our setup every
  section is a linear combination of constant sections with variable
  coefficients, we conclude
  that~\eqref{eq:discrete_non_holonomic_momentum_evolution_eq}
  holds for all sections.
\end{proof}

Below we construct an affine discrete connection adapted to the
present geometry. Later we use that connection to specialize
Theorem~\ref{th:4_points-general} to the horizontal setting.


\subsection{Affine discrete connection for horizontal symmetries}
\label{sec:affine_discrete_connection_for_horizontal_symmetries}

We recall two well known results.
\begin{lemma}
  If $\xi\in \jgsg$, $g\in \SG$ and $q\in Q$, then $\xi_Q(l^Q_g(q)) =
  l^{TQ}_g((Ad_{g^{-1}}(\xi))_Q(q))$.
\end{lemma}

\begin{lemma}\label{le:Jd_equivariance}
  If $g\in \SG$ and $q_0,q_1\in Q$, then
  \begin{equation*}
    J_d(l^{Q\times Q}_g(q_0,q_1)) = Ad_{g^{-1}}^*(J_d(q_0,q_1)).
  \end{equation*}
\end{lemma}

When $F:Q\times Q\rightarrow \R$ is a smooth map,
$D_1F:p_1^*TQ\rightarrow \R$ is defined by $D_1F =
dF|_{p_1^*TQ}$. Noticing that $p_1^*TQ \simeq TQ\times Q$, we have
$D_2(D_1F):p_2^*TQ\rightarrow \R$ (where $p_2:TQ\times Q\rightarrow Q$
is the projection). As $p_2^*TQ \simeq TQ\times TQ$, it is customary
to consider $D_2D_1L_d:TQ\times TQ\rightarrow\R$. It is easy to check
that $D_2D_1L_d$ is bilinear in the tangent vectors.

\begin{definition}
  Let $L_d:Q\times Q\rightarrow\R$ be a discrete lagrangian. We say
  that $L_d$ is \jdef{regular} at $(q_0,q_1)\in Q\times Q$ if the
  bilinear mapping $D_2 D_1 L_d(q_0,q_1) : T_{q_0}Q\times T_{q_1}Q
  \rightarrow \R$ is nondegenerate, \emph{i.e.}, if $X_0\in T_{q_0}
  Q$ satisfies $D_2 D_1 L_d(q_0,q_1) (X_0,X_1) = 0$ for all $X_1\in
  T_{q_1} Q$, then $X_0=0$. In coordinates, the regularity condition
  becomes that the matrix $\frac{\del^2 L_d(q_0,q_1)}{\del q_0 \del
    q_1}$ be invertible.
\end{definition}

\begin{definition}
  Let $\SG$ be a symmetry group of
  $(Q,L_d,\mathcal{D},\mathcal{D}_d)$. We say that $L_d$ is
  \jdef{$\SG$-regular} at $(q_0,q_1)\in Q\times Q$ if the restriction
  of the bilinear form $D_2 D_1 L_d (q_0,q_1): T_{q_0}Q\times T_{q_1}Q
  \rightarrow \R$ to $\mathcal{V}^\SG(q_0)\times \mathcal{V}^\SG(q_1)$
  is nondegenerate.
\end{definition}

The notions of regularity introduced above have already been
considered by other authors
(\cite{ar:marsden_west-discrete_mechanics_and_variational_integrators},~\cite{ar:cortes_martinez-non_holonomic_integrators},~\cite{ar:mclachlan_perlmutter-integrators_for_nonholonomic_mechanical_systems}).
In order to study the relationship between the regularities of a
lagrangian and the fact that the discrete momentum values be regular
we begin by recalling the following fact.

\begin{lemma}
  Let $\pi :Q\rightarrow Q/\SG$ a principal bundle and $\{v_1,\ldots,
  v_k\}\subset \jgsg$ a linearly independent subset. Then, if $q\in Q$,
  $\{(v_1)_Q(q),\ldots, (v_k)_Q(q)\}\subset T_qQ$ is a linearly
  independent subset. Furthermore, if the first set is a basis, then
  the second one is a basis too.
\end{lemma}

\begin{proposition}\label{prop:mu_reg_value_of_Jd}
  Let $\SG$ be a horizontal symmetry group of
  $(Q,L_d,\mathcal{D},\mathcal{D}_d)$ with regular $L_d$, discrete
  momentum mapping $J_d:Q\times Q\rightarrow\jgsg^*$ and
  $\mu\in\jgsg^*$. Then,
  \begin{itemize}
  \item $\mu$ is a regular value of $J_d$ and, consequently,
    $\mathcal{J}_\mu:=J_d^{-1}(\mu)\subset Q\times Q$ is a
    submanifold.
  \item if for $q_0\in Q$ we let $J_d^{q_0}:Q\rightarrow \jgsg^*$ by
    $J_d^{q_0}(q_1) := J_d(q_0,q_1)$, $\mu$ is a regular value of
    $J_d^{q_0}$ and, consequently, $(J_d^{q_0})^{-1}(\mu)\subset Q$ is
    a submanifold that, if not empty, has dimension $\dim Q-\dim \SG$.
  \end{itemize}
\end{proposition}

\begin{proof}
  We have to prove that for all $(q_0,q_1)\in\mathcal{J}_\mu$, the map
  \begin{equation*}
    dJ_d(q_0,q_1) : T_{(q_0,q_1)} (Q\times Q) \rightarrow T_\mu\jgsg^*
    \simeq \jgsg^*
  \end{equation*}
  is onto. Let $\{e_1,\ldots, e_r\}$ be a basis of $\jgsg$ and
  $\{e_1^*,\ldots,e_r^*\}$ its dual basis. Hence
  \begin{equation*}
    J_d(q_0,q_1) = \sum\nolimits_{j=1}^r \phi_j(q_0,q_1) e_j^*,
  \end{equation*}
  where $\phi_j(q_0,q_1) = J_d(q_0,q_1)(e_j) =
  -D_1L_d(q_0,q_1)(e_j)_Q(q_0)$. Then, 
  \begin{equation*}
      dJ_d(q_0,q_1)(X_0,X_1) = \sum\nolimits_{j=1}^r
      (D_1\phi_j(q_0,q_1)(X_0)+D_2\phi_j(q_0,q_1)(X_1)) e_j^*.
  \end{equation*}
  
  Given $\psi = \sum_{j=1}^r a_j e_j^*\in\jgsg^*$, by the regularity
  of $L_d$ there is $X_1\in T_{q_1}Q$ such that
  \begin{equation*}
    D_2\phi_j(q_0,q_1)(X_1) = -D_2 D_1 L_d(q_0,q_1)((e_j)_Q(q_0),X_1) = a_j
  \end{equation*}
  for all $j$. Then,
  \begin{equation*}
    dJ_d(q_0,q_1)(0,X_1) = \sum\nolimits_{j=1}^r
    D_2\phi_j(q_0,q_1)(X_1) e_j^* = \sum\nolimits_{j=1}^r a_j e_j^* = \psi,
  \end{equation*}
  so that $dJ_d(q_0,q_1)$ is onto; thus $\mu$ is a regular value of
  $J_d$ and standard results allow us to conclude that
  $\mathcal{J}_\mu\subset Q\times Q$ is a submanifold.

  We see in the previous computation that, for $\psi\in \jgsg^*$
  \begin{equation*}
    dJ_d^{q_0}(q_1)(X_1) = D_2 J_d(q_0,q_1)(X_1) = \psi, 
  \end{equation*}
  so that $\mu$ is also a regular value of $J_d^{q_0}$, hence,
  $(J_d^{q_0})^{-1}(\mu)\subset Q$ is also a submanifold. The
  dimension of $(J_d^{q_0})^{-1}(\mu)$ can be computed noticing that,
  being $\mu$ a regular value of $J_d^{q_0}$, if $q_1\in
  (J_d^{q_0})^{-1}(\mu)$, we have
  \begin{equation*}
    \begin{split}
      \dim((J_d^{q_0})^{-1}(\mu)) =& \dim(\ker(dJ_d^{q_0}(q_1))) =
      \dim(T_{q_1}Q) - \dim(\im(dJ_d^{q_0}(q_1))) \\ &= \dim(Q)- \dim(\SG).
    \end{split}
  \end{equation*}
\end{proof}

\begin{proposition}\label{prop:Jd_transversality}
  Let $\SG$ be a horizontal symmetry group of
  $(Q,L_d,\mathcal{D},\mathcal{D}_d)$ with $L_d$ regular and
  $\SG$-regular. Then, for all $q_0\in Q$ and $\mu\in\jgsg^*$, if
  $q_1\in l^Q_\SG(\{q_0\}) \cap (J_d^{q_0})^{-1}(\mu)$, then $T_{q_1}Q
  = T_{q_1}l^Q_\SG(\{q_0\}) \oplus T_{q_1}(J_d^{q_0})^{-1}(\mu)$.
\end{proposition}

\begin{proof}
  Assume that $(J_d^{q_0})^{-1}(\mu)\neq\emptyset$ since, otherwise,
  the statement is valid. As $\dim (J_d^{q_0})^{-1}(\mu) = \dim Q -
  \dim \SG$, it suffices to see that $T_{q_1}l^Q_\SG(\{q_0\}) \cap
  T_{q_1}(J_d^{q_0})^{-1}(\mu) =\{0\}$ for all $q_1$ as in the
  statement. If $X_1$ is in this last intersection,
  \begin{equation*}
    0 = dJ_d^{q_0}(q_1)(X_1) = 
    \sum\nolimits_{j=1}^r - D_2 D_1L_d(q_0,q_1)((e_j)_Q(q_0),X_1) e_j^*, 
  \end{equation*}
  so that $D_2 D_1L_d(q_0,q_1)((e_j)_Q(q_0),X_1)=0$, for all $j$. By
  the $\SG$-regularity of $L_d$ and being
  $\{(e_1)_Q(q_0)\ldots,(e_k)_Q(q_0)\}$ a basis of
  $T_{q_1}l^Q_\SG(\{q_0\})$, we have that $X_1=0$.
\end{proof}

\begin{definition}
  Let $\SG$ be a horizontal symmetry group of
  $(Q,L_d,\mathcal{D},\mathcal{D}_d)$ with $L_d$ regular and
  $\SG$-regular. Given $\mu\in\jgsg^*$ we say that $\SG$ is a group of
  \jdef{$\mu$-good symmetries} if, in addition, for each $q\in Q$
  there is a unique $g\in \SG$ such that $J_d^{q}(l^Q_g(q)) = \mu$. In
  this case, we define $\gamma:Q\rightarrow \SG$ by $\gamma(q) :=
  g$. 
\end{definition}

It is possible to extend the previous notion to systems where there
are more than one $g\in \SG$ with the required property but, in this
case, the action must be accompanied by a smooth unique determination
of $g$.

\begin{proposition}\label{prop:mu_good_imp_aff_conn}
  Let $\SG$ be a group of $\mu$-good symmetries of
  $(Q,L_d,\mathcal{D},\mathcal{D}_d)$ for some $\mu\in \jgsg^*$. If
  $\mu$ satisfies $Ad^*_{g}(\mu)=\mu$ for all $g\in \SG$, then
  $\mathcal{J}_\mu\subset Q\times Q$ defines an affine discrete
  connection in $\pi:Q\rightarrow Q/\SG$ of level $\gamma$, given by
  the $\mu$-goodness of $\SG$.
\end{proposition}

\begin{proof}
  By Proposition~\ref{prop:mu_reg_value_of_Jd}, $\mathcal{J}_\mu$ is a
  submanifold of $Q\times Q$ and, by Lemma~\ref{le:Jd_equivariance},
  it is also $\SG$-invariant. According to the definition of $\gamma$,
  $J_d(q,l^Q_{\gamma(q)}(q))=\mu$, so that $\Gamma\subset
  \mathcal{J}_\mu$. Last, by
  Proposition~\ref{prop:mu_reg_value_of_Jd},
  $\mathcal{J}_\mu(q)\subset Q$ is a submanifold and, using
  Proposition~\ref{prop:Jd_transversality}, we conclude that
  $\mathcal{J}_\mu$ defines an affine discrete connection. Notice that
  $\gamma$ is $\SG$-equivariant by Lemma~\ref{le:Jd_equivariance} and
  the condition $Ad^*_{g}(\mu)=\mu$ for all $g\in \SG$.
\end{proof}

\begin{remark}\label{rem:full_stabilizer_is_no_condition}
  In the context of Proposition~\ref{prop:mu_good_imp_aff_conn}, the
  condition that $Ad^*_g(\mu)=\mu$ for all $g\in G$ can be avoided by
  considering the (probably smaller) symmetry group $G_\mu:=\{g\in
  G:Ad^*_g(\mu) = \mu\}$ instead of $G$. Indeed, $G_\mu$ is a Lie
  group and a symmetry group of
  $(Q,L_d,\mathcal{D},\mathcal{D}_d)$. Since $G$
  satisfies~\eqref{eq:horizontal_symmetry-continuous_condition},
  $\mathcal{V}^G(q)\subset \mathcal{D}_q$ for all $q$, so that
  $\mathcal{V}^{G\mu}(q)\subset \mathcal{D}_q$ for all
  $q$. Hence,~\eqref{eq:horizontal_symmetry-continuous_condition}
  holds if we put $G_\mu$ instead of $M$ and $G$ and run the reduction
  arguments in this setting. This same remark applies to the
  statements made in what remains of this section.
\end{remark}


\subsection{Reduction}
\label{sec:reduction-horizontal}

In this section we assume that $\SG$ is a group of $\mu$-good
symmetries of a regular and $\SG$-regular system for some
$\mu\in\jgsg^*$. In the previous section we saw that if $\mu$
satisfies $Ad_{g}^*(\mu)=\mu$ for all $g\in \SG$, we can define an
affine discrete connection $\CD$ such that $Hor_{\CD} =
\mathcal{J}_\mu$. Below, we use this connection to specialize
Theorem~\ref{th:4_points-general} to the context of horizontal
symmetries.

As in the Chaplygin case studied in Section~\ref{sec:an_inclusion},
using $\mathcal{Y}$ we define a forced discrete mechanical system on
$Q/\SG$ with discrete lagrangian
$\breve{L}_d(r_0,r_1):=\mathcal{Y}^*(\hat{L}_d) = \hat{L}_d\circ
\mathcal{Y}$, variational constraints $\breve{\mathcal{D}} :=
d\pi(\mathcal{D})$, kinematic constraints $\breve{\mathcal{D}}_d :=
\mathcal{Y}^{-1}(\hat{\mathcal{D}}_d)$ and forces $\breve{F} :=
\mathcal{Y}^*(\hat{F})$.

\begin{theorem}\label{th:4_points_horizontal}
  Let $\SG$ be a $\mu$-good symmetry group of
  $(Q,L_d,\mathcal{D},\mathcal{D}_d)$ for some $\mu\in\jgsg^*$. Assume
  that $Ad_g^*(\mu)=\mu$ for all $g\in \SG$. Let $q_\cdot$ be a
  discrete curve in $Q$ and $r_k:=\pi(q_k)$ the corresponding discrete
  curve in $Q/\SG$. Then, if $J_d(q_k,q_{k+1})=\mu$ for all $k$, the
  following statements are equivalent.
  \begin{enumerate}
  \item \label{it:var_pple-horizontal} $(q_k,q_{k+1})\in
    \mathcal{D}_d$ for all $k$ and $q_\cdot$ satisfies the variational
    principle $dS_d(q)(\delta q_\cdot) = 0$ for all vanishing end
    points variation $\delta q_\cdot$ such that $\delta q_k \in
    \mathcal{D}_{q_k}$ for all $k$.
  \item \label{it:eq_lda-horizontal} $(q_k,q_{k+1})\in\mathcal{D}_d$
    for all $k$ and $q_\cdot$ satisfies the discrete
    Lagrange--D'Alembert equations~\eqref{eq:dla-eqs} for all $k$.
  \item \label{it:red_var_pple-horizontal} $r_\cdot$ satisfies the
    variational principle
    \begin{equation*}
      d\breve{S}_d(r_\cdot) (\delta r_\cdot) =
      -\sum\nolimits_{k=1}^{N-1} (\breve{F}_d^-(r_k,r_{k+1}) + 
      \breve{F}_d^+(r_{k-1},r_k))(\delta r_k)
    \end{equation*}
    for all vanishing end points variations such that $\delta r_k\in
    \breve{\mathcal{D}}_{r_k}$ for all $k$. In addition,
    $(r_k,r_{k+1})\in \breve{\mathcal{D}}_d$ for all $k$.
  \item \label{it:red_eq_lda-horizontal} $r_\cdot$
    satisfies~\eqref{eq:relation_hor_eqs_hat_breve}. In addition,
    $(r_k,r_{k+1})\in \breve{\mathcal{D}}_d$ for all $k$.
  \end{enumerate}
\end{theorem}

\begin{proof}
  In the present context $J_d$ is conserved over $q_\cdot$. Hence,
  equation~\eqref{eq:discrete_non_holonomic_momentum_evolution_eq}
  is satisfied and, by
  Theorem~\ref{th:vertical_conditions_equivalence} the vertical
  equations in Theorem~\ref{th:4_points-general} are satisfied. Then,
  by Theorem~\ref{th:4_points-general},
  points~\ref{it:var_pple-horizontal} or~\ref{it:eq_lda-horizontal}
  of the present result are equivalent to either one of
  \begin{roman-enumerate}
  \item \label{it:hor_red_intermediate-3} Item~\ref{it:red_var_pple-general} in
    Theorem~\ref{th:4_points-general} holds for all horizontal
    variations, ($\xi_k=0$ for all $k$).
  \item \label{it:hor_red_intermediate-4} For all $k$
    $(v_k,r_{k+1})\in\hat{\mathcal{D}}_d$ and
    condition~\eqref{eq:cond_red_hor_gen-Q}
    holds.
  \end{roman-enumerate}
  Notice that $w_k=\CD(q_k,q_{k+1}) = e$ due to the conservation of
  $J_d$ and the choice of $\CD$. Then,
  $(v_k,r_{k+1})=(\rho(q_k,e),r_{k+1}) = \mathcal{Y}(r_k,r_{k+1}) \in
  \hat{\mathcal{D}}_d$ if and only if $(r_k,r_{k+1})\in
  \mathcal{Y}^{-1}(\hat{\mathcal{D}}_d) = \breve{\mathcal{D}}_d$. Thus
  the constraint conditions that appear in
  items~\ref{it:hor_red_intermediate-3}
  or~\ref{it:hor_red_intermediate-4} are equivalent to the ones that
  appear in~\ref{it:red_var_pple-horizontal}
  or~\ref{it:red_eq_lda-horizontal}.

  By Lemma~\ref{le:relation_hatSd_hatS_d}, the variational principle
  that appears in item~\ref{it:hor_red_intermediate-3} is equivalent
  to the variational principle in
  point~\ref{it:red_var_pple-horizontal}. Similarly, by
  Lemma~\ref{le:relation_hor_eqs_hat_breve}, the equations that appear
  in item~\ref{it:hor_red_intermediate-4} are equivalent to the ones
  that appear in~\ref{it:red_eq_lda-horizontal}.
\end{proof}

The following reconstruction result is the analogue of
Theorem~\ref{th:reconstruction-general} which is valid in the current
setting.

\begin{theorem}\label{th:reconstruction-horizontal}
  Let $\SG$ be a horizontal symmetry group of
  $(Q,L_d,\mathcal{D},\mathcal{D}_d)$ and $(\bar{q}_0,\bar{q}_1)\in
  \mathcal{D}_d$. Let $\mu =
  J_d(\bar{q}_0,\bar{q}_1)\in\jgsg^*$. Assume that $Ad_g^*(\mu)=\mu$
  for all $g\in\SG$ and that $\SG$ is $\mu$-good.

  Let $r_\cdot$ be a discrete trajectory of
  $(Q/\SG,\breve{L}_d,\breve{\mathcal{D}},\breve{\mathcal{D}}_d,\breve{F}_d)$
  such that $\pi(\bar{q}_j)=r_j$ for $j=0,1$. Define the discrete
  curve $q_\cdot$ in $Q$ inductively by $q_0=\bar{q}_0$, and $q_{k+1}
  = \HLds{q_k}(r_{k+1})$ for all $k$ (here the horizontal lift is
  associated to the affine discrete connection whose horizontal space
  is $\mathcal{J}_\mu$). Then $q_\cdot$ is a trajectory of the
  original discrete mechanical system that satisfies $q_0=\bar{q}_0$
  and $q_1=\bar{q}_1$.
\end{theorem}

\begin{proof}
  See the proof of Theorem~\ref{th:reconstruction-chaplygin}.
\end{proof}


\subsection{An example of reduction of horizontal symmetries}
\label{sec:an_example_of_reduction_of_horizontal_symmetries}

In this section we apply the analysis developed in the previous
sections to a system that exhibits horizontal symmetries. The system
is a discretization of the nonholonomic free particle considered by
J. Cort\'es on page 100 of~\cite{bo:cortes-non_holonomic}. More
precisely, the system $(Q,L_d,\mathcal{D},\mathcal{D}_d)$ has
$Q:=\R^3$ and
\begin{equation*}
  L_d(q_0,q_1) := \frac{m}{2} ((q^x_1-q^x_0)^2 + (q^y_1-q^y_0)^2 + 
  (q^z_1-q^z_0)^2),
\end{equation*}
where $q:=(q^x,q^y,q^z)$. The constraints are
\begin{equation*}
  {\mathcal D}_q := \big\{\dot{x}\del_x + \dot{y} \del_y +
  \dot{z}\del_z \in T_q Q: \dot{y} = q^x \dot{x}\big\} = 
  \big\langle (\del_x\big|_q + q^x\del_y\big|_q) , \del_z\big|_q \big\rangle ,
\end{equation*}
\begin{equation*}
  \mathcal{D}_d := 
  \big\{ (q_0,q_1) \in Q\times Q : q^y_1 - q^y_0 = ((q^x_1)^2-(q^x_0)^2)/2\big\}.
\end{equation*}

The group $M:=\R^2$ acts on $Q$ by $l^Q_g(q) :=
(q^x,q^y+g^y,q^z+g^z)$, where $g:=(g^y,g^z)\in M$. The Lie algebra of
$M$ is $\jgm$ that we identify with $\R^2$. The corresponding lifted
action is $l^{TQ}_g(q,v) = (l^Q_g(q),v)$.  Therefore, $L_d$,
$\mathcal{D}$ and $\mathcal{D}_d$ are $M$-invariant.

The vertical space for the action is $\mathcal{V}^M_q = \langle
\del_y\big|_q, \del_z\big|_q\rangle \subset T_qQ$, so that 
\begin{equation*}
  \mathcal{S}_q := \mathcal{D}_q \cap \mathcal{V}^M_q = 
  \langle \del_z|_q \rangle.
\end{equation*}
Thus, the closed subgroup $\SG:=\{0\}\times \R\subset M$, satisfies
$\mathcal{D}_q \cap \mathcal{V}^M_q = \mathcal{S}_q =
\mathcal{V}^\SG_q$, for all $q\in Q$. Hence, $\SG$ is a horizontal
symmetry (sub)group of the system.  We have
$\jgsg=\{0\}\times\R\subset \R^2 =\jgm$.

Clearly, $L_d$ is regular as well as $\SG$-regular. Also, being $\SG$
abelian, $Ad_g^*(\mu)=\mu$ for all $g\in\SG$ and $\mu\in\jgsg^*$. A
simple computation shows that $J_d(q_0,q_1) = m(q^z_1-q^z_0) 1^*$,
where $1^*$ denotes the basis of $\jgsg^*$ that is dual to
$(0,1)$. Then, if $q\in Q$ and $\mu=\mu^z 1^*$, $J_d^q(l_g^Q(q)) =
\mu$ has a unique solution $g:=(0,\frac{1}{m}\mu^z)\in \SG$ and $\SG$
is a group of $\mu$-good symmetries. Hence
Theorem~\ref{th:4_points_horizontal} applies to identify the reduced
system with a forced discrete mechanical system on $Q/\SG\simeq
\R^2$. Below we give an explicit description of this reduced system.

As in the previous sections, we use the affine discrete connection
$\CD$ whose horizontal space is, for a fixed $\mu\in \jgsg^*$,
$\mathcal{J}_\mu = \{(q_0,q_1)\in Q\times Q : m(q^z_1-q^z_0) =
\mu^z\}$. Equivalently, $\CD(q_0,q_1) =
(0,q^z_1-q^z_0-\frac{1}{m}\mu^z)\in \SG$. The discrete horizontal lift
of $\CD$ is $\HLds{q}(r) = (r',r'',\frac{1}{m}\mu^z+q^z_0)$.

As an intermediate step we have to describe the reduced system on
$\RS$. Since $\pi:Q\rightarrow Q/\SG$ is a trivial principal bundle with
structure group $\SG$ we apply the description of the reduced system
given in
Section~\ref{sec:reduced_equations_of_motion_for_trivial_bundles}. The resulting system is
\begin{equation*}
  \begin{cases}
    \RS \simeq (Q/\SG)\times \SG\times (Q/\SG) = \R^2 \times \R \times \R^2,\\
    \hat{L}_d^t(r_k,\vartheta_k,r_{k+1}) = \frac{m}{2}\big( (r_{k+1}'-r_k')^2 + (r_{k+1}''-r_k'')^2 +
    (\frac{1}{m}\mu^z+\vartheta_k)^2 \big),\\
\hat{\mathcal{D}}^t_r = \langle \del_{r'}|_r + 
  r'\del_{r''}|_r \rangle \subset T_r(Q/\SG),\\
\hat{\mathcal{D}}^t_d = \{(r_k,\vartheta_k,r_{k+1})\in 
  (Q/\SG)\times \SG\times (Q/\SG) : r_{k+1}''-r_k'' = 
  \big((r_{k+1}')^2-(r_k')^2\big)/2\}.
  \end{cases}
\end{equation*}
In order to complete the description of the reduced system we need to
compute the reduced forces $\hat{F}_d$. We fix a splitting
\begin{equation*}
  T_qQ = \underbrace{\langle \del_y|_q\rangle}_{=\mathcal{W}_q} \oplus 
  \underbrace{\{0\}}_{=\mathcal{U}_q} \oplus 
  \underbrace{\langle \del_z|_q\rangle}_{=\mathcal{S}_q} \oplus 
  \underbrace{\langle \del_x|_q + q^x \del_y|_q \rangle}_{=\mathcal{H}_q},
\end{equation*}
which determines the nonholonomic connection $\CC$. In particular,
notice that since $d\CD(q_0,q_1) = (0,dz|_{q_1}-dz|_{q_0})$ and
$Hor_{\CC}(\delta q) \in \langle \del_x\big|_q +
q^x\del_y\big|_q\rangle$ the mixed curvature $\BM$ vanishes, and so
does the reduced force $\hat{F}_d$.

The inclusion $\mathcal{Y}$ is
$\mathcal{Y}(r_k,r_{k+1}):=(r_k,e,r_{k+1})$ so that the reduced
mechanical system
$(\breve{Q},\breve{L}_d,\breve{\mathcal{D}},\breve{\mathcal{D}}_d)$
has $\breve{Q}=Q/\SG$ and
\begin{equation}
  \label{eq:discrete_mechanical_system_from_reduced_system-horizontal}
  \begin{cases}
    \breve{L}_d(r_k,r_{k+1}) = \frac{m}{2} \big(
    (r_{k+1}'-r_k')^2 + (r_{k+1}'' - r_k'')^2 + (\frac{1}{m}\mu^z)^2 \big),\\
    \breve{\mathcal{D}}_r 
    = \langle \del_{r'}\big|_r + r'\del_{r''}\big|_r\rangle,\\
    \breve{\mathcal{D}}_d = \{
    (r_{k},r_{k+1})\in (Q/\SG)\times (Q/\SG) :r_{k+1}''-r_k'' =
    \frac{(r_{k+1}')^2-(r_k')^2}{2}\}.
  \end{cases}
\end{equation}
The discrete Lagrange--D'Alembert equations~\eqref{eq:dla-eqs} for
this system are
\begin{equation*}
  \begin{cases}
    -(r_{k+1}'-r_{k}') + (r_k'-r_{k-1}') = \lambda_k r_k',\\
    -(r_{k+1}''-r_{k}'') + (r_k''-r_{k-1}'') = -\lambda_k,\\
    r_{k+1}''-r_k'' = ((r_{k+1}')^2-(r_k')^2)/2.
  \end{cases}
\end{equation*}
Those equations can be easily solved to determine the evolution of the
reduced system, $r_\cdot$. The corresponding trajectory $q_\cdot$ of
the original system is obtained recursively, according to
Theorem~\ref{th:reconstruction-horizontal}, by, for a given $q_k$,
defining
\begin{equation*}
  \begin{split}
    q_{k+1} =& (q^x_{k+1},q^y_{k+1},q^z_{k+1}) =\HLds{q_k}(r_{k+1}) =
    (r_{k+1}', r_{k+1}'',\frac{1}{m}\mu^z+q^z_k) \\=&
    (r_{k+1}',r_{k+1}'',q^z_1-q^z_0+q^z_k),
  \end{split}
\end{equation*}
or,
\begin{equation*}
  (q^x_k,q^y_k,q^z_k) = (r_k',r_k'', k(q^z_1-q^z_0)+q^z_0).
\end{equation*}

Interestingly, the discrete mechanical system
$(\breve{Q},\breve{L}_d,\breve{\mathcal{D}},\breve{\mathcal{D}}_d)$
defined
by~\eqref{eq:discrete_mechanical_system_from_reduced_system-horizontal}
still has a residual symmetry group. The group $M/\SG \simeq \R$ acts
on $Q/\SG$ via the action induced by $l^Q$ on $Q/\SG$. This is
precisely the system and symmetry group whose reduction and
reconstruction as a Chaplygin system was discussed in
Examples~\ref{ex:continued_example-chaplygin_reduction}
and~\ref{ex:continued_example-chaplygin_reconstruction}. Using those
results we get the trajectories of the original system
\begin{equation*}
  (q^x_k,q^y_k,q^z_k) = (r_k, (r_k^2-r_0^2)/2+y_0, k(q^z_1-q^z_0)+q^z_0),
\end{equation*}
where $r_k$ is determined
from~\eqref{eq:continued_example-reduced_equations_chaplygin-r} with
$r_0=q^x_0$ and $r_1=q^x_1$.

\begin{remark}
  The presence of this second reduction step by $M/\SG$ is an example
  of reduction by stages for discrete mechanical systems, a topic that
  will be discussed elsewhere.
\end{remark}




\def\cprime{$'$} 
\def\polhk#1{\setbox0=\hbox{#1}{\ooalign{\hidewidth\lower1.5ex\hbox{`}\hidewidth\crcr\unhbox0}}}
\providecommand{\bysame}{\leavevmode\hbox to3em{\hrulefill}\thinspace}
\providecommand{\MR}{\relax\ifhmode\unskip\space\fi MR }
\providecommand{\MRhref}[2]{%
  \href{http://www.ams.org/mathscinet-getitem?mr=#1}{#2}
}
\providecommand{\href}[2]{#2}


\end{document}